\documentclass{amsart}

\usepackage{amssymb}

\usepackage[all]{xy}

\usepackage{hyperref}

 \usepackage[dvipsnames,usenames]{xcolor}
 \hypersetup{
 colorlinks=true,
 linkcolor={blue!50!black},
 citecolor={green!50!black},
 urlcolor={red!80!black}
 }

\makeatletter
\renewcommand\thesubsection   {\thesection.\@arabic\c@subsection}
\renewcommand\thesubsubsection{\thesubsection .\@arabic\c@subsubsection}
\renewcommand\theparagraph    {\thesubsubsection.\@arabic\c@paragraph}
\renewcommand\thesubparagraph {\theparagraph.\@arabic\c@subparagraph}
\makeatother

\newtheorem{thm}{Theorem}

\newtheorem{lem}[thm]{Lemma}
\newtheorem{cor}[thm]{Corollary}

\newtheorem{prop}[thm]{Proposition}

\newtheorem{conj}[thm]{Conjecture}
   \newtheorem{thm-defn}[thm]{Theorem-Definition}

\theoremstyle{definition}
\newtheorem{defn}[thm]{Definition}

\newtheorem{say}[thm]{}
\newtheorem{exmp}[thm]{Example}

\newtheorem{ques}[thm]{Question}    

\newtheorem{rem}[thm]{Remark}          

\newtheorem*{ack}{Acknowledgments}       

\newtheorem{defn-thm}[thm]{Definition--Theorem}  
\newtheorem{defn-lem}[thm]{Definition--Lemma}  

\newtheorem{temp-defn}[thm]{Temporary Definition}

\newtheorem{variants}[thm]{Variants}

\theoremstyle{remark}

\newtheorem*{cor-no-num}{Corollary}

\let \cedilla =\c
\renewcommand{\c}[0]{{\mathbb C}}  

\renewcommand{\o}[0]{{\mathcal O}} 
\newcommand{\z}[0]{{\mathbb Z}}
\newcommand{\n}[0]{{\mathbb N}}

\renewcommand{\a}[0]{{\mathbb A}}

\newcommand{\p}[0]{{\mathbb P}}
\newcommand{\f}[0]{{\mathbb F}}
\newcommand{\q}[0]{{\mathbb Q}}

\newcommand{\map}[0]{\dasharrow}
\newcommand{\qtq}[1]{\quad\mbox{#1}\quad}
\newcommand{\spec}[0]{\operatorname{Spec}}
\newcommand{\pic}[0]{\operatorname{Pic}}

\newcommand{\mult}[0]{\operatorname{mult}}

\newcommand{\supp}[0]{\operatorname{Supp}}    
\newcommand{\red}[0]{\operatorname{red}}    
\newcommand{\codim}[0]{\operatorname{codim}}

\newcommand{\Hom}[0]{\operatorname{Hom}}

\newcommand{\aut}[0]{\operatorname{Aut}}

\newcommand{\sing}[0]{\operatorname{Sing}}

\newcommand{\tors}[0]{\operatorname{tors}}

\newcommand{\chr}[0]{\operatorname{char}}

\newcommand{\cl}[0]{\operatorname{Cl}}

\newcommand{\len}[0]{\operatorname{length}}

\newcommand{\ann}[0]{\operatorname{Ann}}
\newcommand{\diag}[0]{\operatorname{diag}}
\newcommand{\ord}[0]{\operatorname{ord}}

\newcommand{\mor}[0]{\operatorname{Mor}} 
 
\newcommand{\tsum}[0]{\textstyle{\sum}}
\newcommand{\tprod}[0]{\textstyle{\prod}}

\newcommand{\cdiv}[0]{\operatorname{CDiv}}

\def\into{\DOTSB\lhook\joinrel\to}

\def\loccoh#1.#2.#3.#4.{H^{#1}_{#2}(#3,#4)}

\DeclareMathAlphabet{\mathchanc}{OT1}{pzc}%
                                {m}{it}

 \newcommand{\ass}[0]{\operatorname{Ass}}

\newcommand{\emb}[0]{\operatorname{emb}}

\newcommand{\pure}[0]{\operatorname{pure}}

\newcommand{\grass}[0]{\operatorname{Gr}}

\newcommand{\chow}[0]{\operatorname{Chow}}

\newcommand{\chh}[0]{\operatorname{Ch}}

\newcommand{\dsupp}[0]{\operatorname{DSupp}}
\newcommand{\ssupp}[0]{\operatorname{SSupp}}

\newcommand{\flag}[0]{\operatorname{Flag}}
\newcommand{\icorr}[0]{\operatorname{Inc}}
\newcommand{\fcm}[0]{\operatorname{FlatCM}}
\newcommand{\floc}[0]{\operatorname{Flat}}

\newcommand{\kdiv}[0]{\operatorname{KDiv}}

\newcommand{\vtors}[0]{\operatorname{v-tors}}
\newcommand{\vpure}[0]{\operatorname{v-pure}}

\newcommand{\mcl}[0]{\operatorname{MCl}}

\newcommand{\chull}[0]{\operatorname{CHull}}

\begin{document}

\bibliographystyle{amsalpha}
\hfill\today

 \title{Families of divisors}
 \author{J\'anos Koll\'ar}

\begin{abstract}  We establish a new moduli theory for divisors,  that interpolates between the Hilbert scheme and the Cayley-Chow variety.
This completes the last step in the construction of a good moduli theory for stable pairs  $(X,\Delta)$.
\end{abstract} 
 \maketitle

\tableofcontents

A persistent problem in the moduli theory of pairs $(X,\Delta)$ is that, while the underlying varieties form  flat families, the divisorial parts $\Delta$ do not. Neither of the two main traditional methods of parametrizing varieties or schemes  gives the right answer for  the divisorial part.
\begin{itemize}
\item Hilbert schemes take into account embedded points, but we need to ignore them entirely.
\item Cayley-Chow varieties work well only over seminormal  schemes, but we wish to have a theory over arbitrary base schemes.
\end{itemize}

Our aim is to develop a theory that interpolates between these two, managing to keep from both of them the properties that we need. 
The definition of Mumford divisors codifies basic properties of the divisorial part of families of stable pairs, though with a  new name introduced in \cite{k-mumf}.
The main new result is 
Definition~\ref{K-flat.defn.1}. The rest of the paper is then devoted to proving that it has all the hoped-for properties.

 \begin{defn}[Mumford divisors]\label{fam.divs.ch3}   
Let $f:X\to S$ be a flat morphism with $S_2$ fibers of pure dimension $n$.
A  subscheme $D\subset X$  is a  {\it relative  Mumford   divisor}  or a {\it generically flat family of  Mumford divisors} over $S$   if  
there is an 
 open subset $U\subset X$  such that 
\begin{enumerate}
\item  $\codim_{X_s}(X_s\setminus U)\geq 2$ for every $s\in S$,
\item $D|_U$ is a relative Cartier divisor,
\item $D$ is the closure of $D|_U$ and
\item $X_s$ is smooth at generic points of $D_s$ for every $s\in S$.
\end{enumerate}
If  $U\subset X$  denotes the largest open set with the properties (1--2),
then $U$ is  called the {\it Cartier locus} and  $Z:=X\setminus U$   the {\it non-Cartier locus} of $D$.

Let $q:W\to S$ be any morphism.
We have a fiber product diagram
$$
\begin{array}{ccc}
X_W & \stackrel{q_X}{\to} &  X\\
f_W\downarrow\hphantom{f_W} && \hphantom{f}\downarrow f\\
W & \stackrel{q}{\to} &   \hphantom{.}S.
\end{array}
\eqno{(\ref{fam.divs.ch3}.5)}
$$
 Then   $q_X^*\bigl(D|_U\bigr)$ is a   relative Cartier divisor
on  $U_W:=q_X^{-1}(U)$ and its closure
 is a relative Mumford divisor, called the {\it divisorial pull-back} of $D$ by $q$. It is denoted by $q^{[*]}_XD$ or simply by $D_W$.
If $q$ is flat then  $q^{[*]}_XD=q^{*}_XD=W\times_SD$.
\end{defn}

\begin{defn}[K-flatness]\label{K-flat.defn.1} Let $f:X\to S$ be a flat, projective morphism with $S_2$ fibers of pure dimension $n$.   A relative Mumford divisor $D\subset X$   
 is {\it K-flat} over $S$ iff one the following---increasingly more general---conditions hold.
\begin{enumerate} 
\item ($S$ local with infinite residue field)  For  every finite morphism  $\pi: X\to \p^n_S$,  $\pi_*D\subset \p^n_S$ is a relative Cartier divisor.
\item ($S$ local) $q^*D$ is K-flat over $S'$ for some (equivalently every)
flat, local morphism $q:S'\to S$, where $S'$ has infinite residue field.
\item ($S$ arbitrary)  $D$ is K-flat over every  localization of  $S$.
\end{enumerate}

Let us start with some comments on the definition.
\begin{enumerate}\setcounter{enumi}{3} 
\item  Here K stands for the first syllable of Cayley.  
We  use C-flat for a closely related (possibly equivalent) notion; see (\ref{CK.flat.verss.defns}).
\item The definition of $\pi_*D$ is not always obvious; in essence Section~\ref{divisorial.support.sect} is entirely devoted to establishing it.
However, $\pi_*D$ equals the scheme-theoretic image of $D$ if    $\red D\to \red(\pi(D))$ is birational and 
 $\pi$ is \'etale at every generic point of every fiber $D_s$ (\ref{dsupp.car.nd.dl.4}.2). It is sufficient to check condition (1) for such morphisms  $\pi:X\to \p^n_S$.
\item If $S$ is not local then there may not be any finite morphisms  $\pi: X\to \p^n_S$ (\ref{no.morph.to.prod.exmp}); this is one reason for the 3 step definition.
\item The infinite residue field extensions in (2) are necessary in some cases; see (\ref{Cn.CC.flat.ord.1.exmp}.9).
\item  The definition of 
K-flatness  is global in nature, but we show that it is in fact \'etale local on $X$ (\ref{loc.C.flat.C.flat.lem}).  
\end{enumerate} 
\end{defn}

\subsection*{Good properties of K-flatness}{\ }

K-flat families have several good properties. Some of them are needed for
the moduli theory of stable pairs, but others, for example  (\ref{bertini,intro.thm}--\ref{K.flat.push.fwd}), come as  bonus.

\begin{say}[Functoriality] Being K-flat is preserved by arbitrary base changes
and it descends from faithfully flat base changes  (\ref{CC.flat.4.vers.thm.cor.1}).
Thus we get the functor $KDiv(X/S)$ of K-flat,  relative Mumford divisors 
on $X/S$.  If we have a fixed relatively ample divisor $H$ on $X$, 
then  $KDiv_d(X/S)$ denotes the functor of K-flat,  relative Mumford divisors
of degree $d$.

We have a  disjoint union decomposition  $KDiv(X/S)=\cup_d KDiv_d(X/S)$.
\end{say}

The main result is the following, to be proved in (\ref{CMDIV.exits.thm}).

\begin{thm} \label{KDIV.exits.thm} Let $f:X\to S$ be a flat, projective morphism with $S_2$ fibers of pure dimension $n$.   Then the functor $KDiv_d(X/S)$ of K-flat,  relative Mumford divisors of degree $d$
 is representable by  a separated $S$-scheme of finite type  $\kdiv_d(X/S)$.
\end{thm}

{\it Complement \ref{KDIV.exits.thm}.1.} 
 $\kdiv_d(X/S)$ is proper  over $S$ (and non-empty) iff the fibers of $f$ are normal.
This is, however, not a problem for the moduli of stable pairs.

\begin{say}[Comparison with flatness]\label{comp.w.flat.i}  K-flatness is a generalization of flatness  and it is equivalent to it for smooth morphisms. 
\begin{enumerate}
\item If $D$ is K-flat and $f$ is smooth at $x\in X$, then $f|_D$ is flat at $x$. Equivalently,  $D$ is a relative Cartier divisor in a neighborhood of $x$; see (\ref{mumf.d.Cf.cor.2}). 
\item  If $f|_D$ is flat at $x$  then $f$ is K-flat at $x$; see  (\ref{mumf.d.Cf.cor.3}). 
\end{enumerate}
In particular,  the notion of K-flatness gives something new only at the  points  where $f$ is not smooth and $f|_D$ is not flat. 
\end{say}

\begin{say}[Reduced base schemes] If $S$ is reduced then every relative Mumford divisor is K-flat; see (\ref{dsupp.red.sm.lem}), which in turn  follows directly from \cite[4.36]{k-modbook}. In retrospect, this is the reason 
why the moduli theory of pairs could be developed over reduced base schemes
without this notion in \cite[Chap.4]{k-modbook}. 

So in practice the main task is to understand K-flatness  over  Artin base schemes $S$.  This takes care of the general  case since
 $D\subset X$ is K-flat over $S$ iff  $D_A\subset X_A$ is K-flat over $A$ for every Artin subscheme $A\subset S$; see (\ref{over.art.enough.ptop}).
\end{say}

\begin{thm}[Bertini theorems, up and down]\label{bertini,intro.thm}  Assume that $n\geq 2$ and let $|H|$ be a basepoint-free linear system on $X$. Then $D$ is K-flat iff  $D|_H$ is K-flat for general
$H\in |H|$. 
\end{thm}

This is established by combining  (\ref{bert.sec.lem.2}--\ref{bert.sec.lem.3}) with (\ref{K-flat.=.sC-flat.thm}). 
As a consequence, K-flatness is really a question about families of surfaces and curves on them. There are  similar theorems for families of stable pairs, see \cite[Chaps.2 and 5]{k-modbook} or the original papers
\cite{k-gl1, bha-dej, k-gl2}.

This reduction to surfaces  is very helpful conceptually, but also computationally since we have rather complete lists of singularities of log canonical surface pairs $(X, \Delta)$, at least when the coefficients of $\Delta$ are not too small.

Another variant of the phenomenon, that higher codimension points sometimes do not matter much, is the 
Hironaka-type flatness theorem  \cite[10.68]{k-modbook}, which is 
 a generalization of 
\cite{MR0102519}; see also  \cite[III.9.11]{hartsh}.

\begin{say}[K-flatness does not depend on $X$] \label{K.flat.indep.X}
It is well understood that in the theory of pairs  $(X, \Delta)$ one can not separate the underlying variety $X$ from the divisorial part $\Delta$.  For example, if $X$ is a surface with quotient singularities only, then the pair $(X,D)$  is plt for some smooth curves  $D\subset X$ but not even lc for some other cases. It really matters how exactly $D$ sits inside $X$. 

Thus it is unexpected that K-flatness depends only on the divisor $D$, not on the ambient variety $X$, though maybe this is less  surprising if one thinks of K-flatness as a variant of flatness. 

On the other hand, not all K-flat deformations of $D$ are realized in a given $X$. For example, if we start with  $\bigl(\a^2, (xy=0)\bigr)$ then  every K-flat deformation of the pair induces a  flat deformation of $D_1=(xy=0)\subset \a^2$.  If we have 
$\bigl((xy=z^2), (z=0)\bigr)$ then there are K-flat deformations of the pair that induce a  non-flat deformation of 
$D_2=(xy-z^2=z=0)\cong D_1$. 
\end{say}

\begin{say}[Push-forward] \label{K.flat.push.fwd}
Let $f:X\to S$ and $g:Y\to S$ be  flat, projective morphisms with $S_2$ fibers of pure dimension $n$ and  $\tau:X\to Y$ a finite morphism.
Let $D\subset X$ be a  K-flat relative Mumford divisor   such that $\tau_*D$ is also a 
relative Mumford divisor.  (That is, none of the irreducible components of $D_s$ is mapped to $\sing (Y_s)$.)  Then $\tau_*D$ is also K-flat.
(See (\ref{K-flat.defn.1}.7) and Section~\ref{divisorial.support.sect}
for the correct definition of $\tau_*D$.)

\end{say}

\subsection*{Application to moduli spaces of pairs}{\ }

 The construction of the
moduli space of stable  pairs given in \cite[Sec.4.9]{k-modbook} relies on
a suitable moduli theory of divisors. In characteristic 0, all
restrictions on bases schemes come from the moduli theory of divisors. In
\cite[Sec.4.7]{k-modbook} we used Cayley-Chow theory, which was at that time
worked out for seminormal base schemes, though later it was extended to
reduced base schemes.

Using (\ref{KDIV.exits.thm}) and the methods of  \cite[Chap.4]{k-modbook}, we get a  moduli theory of stable
pairs over arbitrary base schemes in characteristic 0. 

\begin{defn}\label{stable.pair.defn} A {\it family of stable pairs}
is a morphism  $f:(X, cD)\to S$ where
\begin{enumerate}
\item $f:X\to S$ is flat and projective,
\item   $D$ is a K-flat family of divisors on $X$,
\item $K_X+cD$ is $\q$-Cartier, relatively ample and
\item the fibers  $(X_s, cD_s)$ are semi-log-canonical.
\end{enumerate}

The basic invariants of the fibers are the dimension $n=\dim X_s$ and
the {\it volume}  $(K_{X_s}+cD_s)^n$. 

The role of the coefficient $c$ is murkier. It is well understood that, in order to get finite dimensional moduli spaces, 
one needs to control the coefficients of the divisorial part of  stable pairs $(X,\Delta)$; see for example \cite[4.68]{k-modbook}.
Fixing $c$ is one of the easiest way to ensure this control.

If we work with
$\q$-divisors  $\sum a_iD_i$, then a convenient choice of $c$ is
the reciprocal of the least common denominator of the $a_i$. 
Thus  $D:=c^{-1}(\sum a_iD_i)$ is a $\z$-divisor.

Fixing $c$ leads to the largest moduli spaces. In practice one may want to
impose additional restrictions,  handle the different divisors $D_i$ differently and allow real coefficients as well. One way to achieve these is the notion of {\it marked pairs}  \cite[Sec.4.7]{k-modbook}. 
\end{defn}

For now we focus on  the most general form of the basic existence theorem.

\begin{thm}\label{stable.pair.thm}  Fix constants $n,c,v$ and work with schemes over $\q$. 
Then the functor of  families of stable pairs $f:(X, cD)\to S$ of dimension $n$ and volume $v$  has a coarse moduli space
${\rm SP}(n,c,v)$
that is projective over $\q$.
\end{thm}

This theorem represents the culmination of the work of several decades; some of the main contributions are
\cite{ksb,  k-proj, al94, ale-pairs, kk-singbook, k-source, MR3671934, MR3779955, MR3779687}. Many parts of the proof work in positive characteristic, and even over $\z$, but there are several fundamental unsolved questions.

I hope---but do not claim---that K-flatness gives the `optimal' moduli theory
for divisors. By `optimal' I mean that
\begin{itemize}
\item it is defined over arbitrary schemes,
\item it agrees with the notion of Mumford divisors over reduced schemes,
\item it leads to  moduli spaces that include all families that one would wish to consider.
\end{itemize}
K-flatness satisfies the first 2 and it  has surprisingly nice additional features. Once its basic properties are
established, it is quite easy to work with, since we can mostly ignore
singularities that occur in codimension $\geq 3$. On the other hand, it is possible that ignoring codimension $\geq 3$ means that we allow too many infinitesimal deformations and there is a better,  more restrictive theory.

\subsection*{Problems and questions about K-flatness}{\ }

There are also some difficulties with K-flatness. I believe that they do not effect the  general moduli theory of stable pairs, but they may make explicit computations lengthy.

\begin{say}[Hard to compute]
The definition of K-flatness is quite hard to check, since for $X\subset \p^N$ we need to check not just linear projections  $\p^N_S\map \p^n_S$ (\ref{lin.projection.defn}) but 
all morphisms  $X\to \p^n_S$ involving  all linear systems on $X$.

On the other hand, at least in the examples in Sections~\ref{hyp.surf.K.flat.sec}--\ref{s.n.curve.K.flat.sec},
the computation of the restrictions imposed by linear projections is the hard part, the general cases then follow easier. 
It would be good to work out more space
 curves  $C\subset \a^3$.  Hopefully, once we first check linear projections 
$\a^3\to\a^2$,  we are  left with a very short list of possibilities and then the general projections work out. 

\end{say}

\begin{say}[Tangent space and obstruction theory] I do nor know how to write down the tangent space of $\kdiv(X/S)$. A handful of examples are computed in Sections~\ref{hyp.surf.K.flat.sec}--\ref{s.n.curve.K.flat.sec}, but they do not seem to suggest any general pattern.
The obstruction theory of K-flatness is completely open.
\end{say}

\begin{say}[The definition is not formal-local] One expects 
K-flatness to be a formal-local property on $X$,   but there are some (hopefully only technical) problems with this.  See (\ref{fKf.=.Kf.conj}) and (\ref{d=1.scf=kf=ff.prop}) for partial results. 
\end{say}

Over a DVR, every  K-flat deformation of a variety $X$ is a flat deformation of some scheme $X'$ such  that $\red X'=X$.  By (\ref{comp.w.flat.i}.1),
the torsion subsheaf of $\o_{X'}$ is supported on $\sing X$. It would be good to get a good a priori bound on the size of $\tors \o_{X'}$.

\begin{ques}[Bounding the torsion]\label{bd.torrs.ques} Let $(A,m,k)$ be an Artin scheme and $C_A\to \spec A$  a K-flat deformation of a pointed curve $(c,C)$ that is flat on $C\setminus\{c\}$. 
Let  $C_k$ be the central fiber and  $I=\ker[\o_{C_A}\to \o_C]$.
Thus  $\tors_cC_k=I/m\o_{C_A}$  (and  $C_A\to \spec A$ is flat iff $\tors_cC_k=0$ by (\ref{flat.over.art.lem})). 

 What is the best bound for $\tors_cC_k$, depending only on $C$?

\end{ques}

We are always interested in divisors that lie on a particular family of varieties $X\to S$, but, in view of (\ref{K.flat.indep.X}), the following seems also natural.

\begin{ques}[Universal deformation spaces] Let $D$ be a reduced, projective scheme over a field $k$. Is there a universal deformation space for its K-flat deformations?
\end{ques}

\subsection*{Examples}{\ }

The first example  shows that the space of first order deformations of the smooth divisor $(x=0)\subset \a^2$, that are Cartier away from the origin, is infinite dimensional. Thus working with generically flat divisors does not give a sensible moduli space.

\begin{exmp} For $g\in k[x,y]$ consider the ideal
$$
I_g=(x^n, xy+\epsilon g, \epsilon x)\subset k[x,y,\epsilon]\qtq{and set}  D_g=\spec k[x,y,\epsilon]/I_g.
$$
It is easy to check directly that
\begin{enumerate}
\item $D_g$ is Cartier away from the origin,
\item  $(I_g, \epsilon)/(\epsilon)=(x)\cap (x^n, y)$,
\item $D_g$ has no embedded points iff $g\notin (x^n, xy)$ and
\item $D_{g_1}=D_{g_2}$ iff $g_1-g_2\in (x^n, xy)$.
\end{enumerate}
More general computations are done in   (\ref{divs.sheaves.dim.2.say.4}). 
\end{exmp}

\begin{exmp}\label{no.morph.to.prod.exmp} Let $C$ be a smooth projective curve and $E$ a stable vector bundle over $E$ of rank $n+1\geq 2$ and of degree 0. Then there is no finite morphism  $\p_C(E)\to \p^n\times C$.
\end{exmp}

\begin{say}[Description of the sections]
We start by reviewing the divisor theory over Artin schemes in Section~\ref{inf.stud.dic.sch.sec}. The key notion of divisorial support is introduced and studied in
Section~\ref{divisorial.support.sect}.

Several versions of K-flatness are investigated in Section~\ref{K.flat.sect}. For our treatment, technically the most important is C-flatness, which is treated in detail in 
Section~\ref{CC.flat.sect}.  The ideal of Chow equations is introduced in
Section~\ref{chow.eq.ideal.sect}  and the main results are proved in 
Section~\ref{C.K.rep.thms.sect}.

 Sections~\ref{hyp.surf.K.flat.sec}--\ref{s.n.curve.K.flat.sec} are devoted to examples;
 we  describe K-flat deformations of plane curves and of seminormal curves  over $k[\epsilon]$. While the computations are somewhat lengthy, the answers are quite nice in both cases.
\end{say}

\begin{ack} 
Partial  financial support    was provided  by  the NSF under grant number
DMS-1901855.
\end{ack}

\section{Infinitesimal study of Mumford divisors}\label{inf.stud.dic.sch.sec}

\begin{say}\label{loc.glop.prob.claim}
The infinitesimal method to study families of objects  in algebraic geometry posits that we should proceed in 3 broad steps.
\begin{itemize}
\item Study families over Artin schemes.
\item Inverse limits  then give families over complete local schemes.
\item For arbitrary local schemes,  descend properties from the completion.
\end{itemize}  
This approach has been very successful for proper varieties and coherent sheaves on them. 
One of the problems we have with general (possibly non-flat)  families of divisors is that the global and the infinitesimal computations do not match up; in fact they say the opposite in some cases.
We discuss 2 instances of this:
\begin{itemize}
\item Relative Cartier divisors on non-proper varieties.
\item Generically flat families of divisors on surfaces.
\end{itemize} 
The surprising feature is that the two behave quite differently.
We state 2 special cases of the results where the contrasts between 
Artin and DVR bases are especially striking.

\medskip
{\it Claim \ref{loc.glop.prob.claim}.1.}
Let $\pi:X\to (s,S)$ be a smooth, affine morphism to a 
local scheme.
\begin{enumerate}
\item[(a)] If $S$ is Artin then the restriction map $\pic(X)\to \pic(X_s)$ is an isomorphism.
\item[(b)] If $S=\spec k[[t]]$  then $\pic(X)$ is frequently infinite dimensional.
\end{enumerate}
\medskip

Thus there can be many nontrivial line bundles on $X$ over $\spec k[[t]]$, but we do not see them when working over $\spec k[[t]]/(t^m)$; see (\ref{pic.over.artin.say}.3) and  (\ref{pic.const.ell.say}) for details.

\medskip
{\it Claim \ref{loc.glop.prob.claim}.2.}
Let $\pi:X\to (s,S)$ be a smooth morphism of relative dimension 2 to a  local scheme $S$.
\begin{enumerate}
\item[(a)] If $S$ is Artin and non-reduced, then relative class group $\cl(X/S)$ 
(\ref{mcl.defn}) is infinite dimensional.
\item[(b)] If $S=\spec k[[t]]$  then every relative Mumford  divisor $D\subset X$ is Cartier. 
\end{enumerate}
\medskip

As an  example,  one easily computes that
$$
\cl\bigl(\p^1_{k[[t]]}\bigr)\cong \z \qtq{but}
  \cl\bigl(\p^1_{k[[t]]/(t^m)}\bigr)\cong \z+k^{\infty} \qtq{for} m\geq 2.
$$
Note that if $D\subset X$ is a relative Mumford  divisor over $S=\spec k[[t]]/(t^m)$  then it  is Cartier on an open set $X^\circ\subset X$ whose complement is  finite.
We see that the study of  $\cl(X/S)$ is pretty much equivalent to
the study of  $\cdiv(X^\circ/S)$ for every $X^\circ$. 
Here $X^\circ$ is not affine, but it is the next simplest scheme, as far as cohomological dimension is concerned.
Indeed, 
\begin{itemize}\setcounter{enumi}{2}
\item  If $X$ is affine then $H^i(X, F)=0$, for every $i>0$ and for every coherent sheaf $F$ on $X$. 
\item If $X$ is an affine surface and $X^\circ\subset X$ is open  then 
$H^i(X^\circ, F^\circ)=0$, for every $i>1$ and for every coherent sheaf $F^\circ$ on $X^\circ$. 
\end{itemize}
\end{say}

\subsection*{Relative Picard group, examples}{\ }

\begin{say}[Picard group over Artin schemes]\label{pic.over.artin.say}
 Let $(A,m,k)$ be a local Artin ring and $X_A\to \spec A$ a flat morphism. Let $(\epsilon)\subset A$ be an ideal such that $I\cong k$ and set $B=A/(\epsilon)$.
We have an exact sequence
$$
0\to \o_{X_k}\stackrel{e}{\longrightarrow} \o^*_{X_A}\to \o^*_{X_B}\to 1,
\eqno{(\ref{pic.over.artin.say}.1)}
$$
where $e(h)=1+h\epsilon$ is the exponential map. 
We use its  long exact cohomology sequence and induction on $\len A$ to compute $\pic(X_A)$. There are 3 cases that are especially interesting for us.

\medskip
{\it Claim \ref{pic.over.artin.say}.2.} Let $X_A\to \spec A$ be a flat, affine  morphism. Then  the restriction map $\pic(X_A)\to \pic(X_k)$ is an isomorphism.
\medskip

Proof. We use the exact sequence
$$
H^1(X_k, \o_{X_k})\to \pic(X_A)\to \pic(X_B)\to H^2(X_k, \o_{X_k}).
\eqno{(\ref{pic.over.artin.say}.3)}
$$
Since  $X$ is affine, the two groups at the ends vanish, hence we get an isomorphism in the middle. Induction completes the proof. \qed
\medskip

\medskip
{\it Claim \ref{pic.over.artin.say}.4.} Let $X_A\to \spec A$ be a flat, proper  morphism. Assume that $H^0(X_k, \o_{X_k})=k$.
Then  the kernel of the restriction map $\pic(X_A)\to \pic(X_k)$ is a unipotent group scheme of dimension $\leq h^1(X_k, \o_{X_k})\cdot (\len A-1)$.
(If $\chr k=0$ then the kernel is $k$-vector space.)
 \medskip

Proof. By \cite[III.12.11]{hartsh},
$H^0(X_A, \o_{X_A})\to H^0(X_B, \o_{X_B})$ is surjective,  and so is
$H^0(X_A, \o^*_{X_A})\to H^0(X_B, \o^*_{X_B})$. 
Thus we get the exactness of
$$
0\to H^1(X_k, \o_{X_k})\to \pic(X_A)\to \pic(X_B)\to H^2(X_k, \o_{X_k}). \qed
\eqno{(\ref{pic.over.artin.say}.5)}
$$

\medskip
{\it Claim \ref{pic.over.artin.say}.6.} Let $X_A\to \spec A$ be a flat, affine morphism and  $Z\subset X_A$ a closed subset of codimension $\geq 2$. Set
$X_A^\circ:=X_A\setminus Z$. 
 Assume that $X_k$ is $S_2$.
Then  the kernel of the restriction map $\pic(X^\circ_A)\to \pic(X^\circ_k)$ is a unipotent group scheme of dimension $\leq h^1(X^\circ_k, \o_{X^\circ_k})\cdot (\len A-1)$.

 \medskip

Proof. Since $X_k$ is $S_2$, $H^0(X^\circ_k, \o_{X^\circ_k})\cong H^0(X_k, \o_{X_k})$
and similarly for $X_A$. Thus
$H^0(X^\circ_A, \o^*_{X^\circ_A})\to H^0(X^\circ_B, \o^*_{X^\circ_B})$ is surjective,  and the rest of the argument works as in (\ref{pic.over.artin.say}.5). \qed

\medskip
{\it Remark \ref{pic.over.artin.say}.7.}  Although (\ref{pic.over.artin.say}.6) is very similar to (\ref{pic.over.artin.say}.4), a key difference is that in (\ref{pic.over.artin.say}.6) the group
$ H^1(X^\circ_k, \o_{X^\circ_k})$ can be infinite dimensional. Indeed, 
$H^1(X^\circ_k, \o_{X^\circ_k})\cong H^2_Z(X_k, \o_{X_k})$ and it is
\begin{enumerate}
\item[(a)] infinite dimensional if $\dim X_k=2$, 
\item[(b)] finite dimensional if $X_k$ is $S_2$ and $\codim_{X_k}Z\geq 3$,  and
\item[(c)] 0 if $X_k$ is $S_3$ and $\codim_{X_k}Z\geq 3$.
\end{enumerate}
See, for example, \cite[Sec.10.2]{k-modbook} for these claims. 

\medskip
{\it Remark \ref{pic.over.artin.say}.8.}  If $H^2(X_k, \o_{X_k})=0$ then the
$\leq$ in (\ref{pic.over.artin.say}.4) and (\ref{pic.over.artin.say}.6)
are equalities. If the characteristic is 0 then equality holds even if
$H^2(X_k, \o_{X_k})\neq 0$; see \cite{blr}. 

\end{say}

The following immediate consequence of (\ref{pic.over.artin.say}.7.c) is especially useful for us; see also  \cite[4.36]{k-modbook}.

\begin{cor} \label{pic.codim3.cor}
Let $X\to S$ be a smooth morphism,  $D\subset X$ a closed subscheme and $Z\subset X$ a closed subset. Assume that
\begin{enumerate}
\item $D$ is a relative Cartier divisor on $X\setminus Z$,
\item  $D$ has no embedded points in $Z$ and
\item $\codim_{X_s}Z_s\geq 3$ for every $s\in S$.
\end{enumerate}
Then $D$ is a  relative Cartier divisor.\qed
\end{cor}

The following is essentially in \cite[XIII]{sga2}, 
see also  \cite[2.93]{k-modbook}.

\begin{thm}\label{cart.infinit.cond.thm} Let $X\to S$ be a flat morphism with $S_2$ fibers and $D$ a divisorial subscheme. Let $U\subset X$ be an open subscheme such that $D|_U$ is relatively Cartier and
$\codim_{X_s}(X_s\setminus U_s)\geq 2$ for every $s\in S$.

Then $D$ is relatively Cartier iff  the divisorial pull-back  $\tau^{[*]}D$ 
(\ref{fam.divs.ch3}.5)  is  relatively Cartier
for every Artin subscheme $\tau:A\into S$. \qed
\end{thm}

Over Artin rings, we have the following flatness criterion.
For a coherent sheaf $F$ let $\emb F$ denote the largest subsheaf whose support is the union of the (closures of the) embedded points of $F$.  

\begin{lem}\label{flat.over.art.lem}
 Let $(A,m,k)$ be an local  Artin ring,
$g:X\to \spec A$ a morphism and $F$ a coherent sheaf on $X$.
Assume that $F$ is generically flat over $A$ and $\emb F=0$. 
Then $F$ is flat over $A$ iff $\emb (F_k)=0$. 
\end{lem}

Proof. Choose  $\epsilon\in m$ such that $m\epsilon=0$.
If $F$ is flat over $A$ then $\epsilon F\cong F_k$, thus we get an injection
$\epsilon:\emb (F_k)\into \emb F$. Thus if $\emb F=0$ then so is
$\emb (F_k)$. 

Conversely, assume that $\emb (F_k)=0$. We may assume that $X$ is affine. 
By induction on $\len A$ we may assume that  $(F/\epsilon F)/\emb(F/\epsilon F)$ is flat over $A/(\epsilon)$. 
We claim that  $\emb(F/\epsilon F)\subset \emb F$.

By assumption $(F/\epsilon F)/\emb(F/\epsilon F)$  is a free $A/(\epsilon)$ module,
choose  basis elements $f_{\lambda}$ and lift them back to $\tilde f_{\lambda}\in H^0(X, F)$.

Let $(\epsilon F)^{(1)}\subset F$ be the preimage of $\emb(F/\epsilon F)$. 
Pick now $h\in (\epsilon F)^{(1)}$. 
The image of $h$ in $(F/mF)$ is 0, so 
$h=\sum a_ig_i$ for some $a_i\in m, g_i\in F$. Write each 
$g_i $ in the $\tilde f_{\lambda}$ basis. Thus we have
$$
g_i\equiv\tsum_{\lambda}c_{i\lambda} \tilde f_{\lambda}\mod (\epsilon F)^{(1)}.
$$
Since  $m(\epsilon F)^{(1)}=0$, we get that
$$
h=\tsum_{\lambda}\bigl(\tsum_i a_ic_{i\lambda}\bigr) \tilde f_{\lambda}.
$$
 This is zero modulo $(\epsilon F)^{(1)}$, so $ \tsum_i a_ic_{i\lambda}\in (\epsilon)$ for every $\lambda$. Thus $h\in \epsilon F$. 

Thus $\epsilon F\cong F/mF$ and $\emb(F/\epsilon F)=0$,
so $F$  is flat over $A$. 
\qed
\medskip

Relative Cartier divisors also have some unexpected properties over non-reduced base schemes. These do not cause theoretical problems, but it is good to keep them in mind. 

\begin{exmp}[Cartier divisors over {$k[\epsilon]$}]
Let $R$ be an integral domain over a field $k$. Relative principal ideals
in $R[\epsilon]$ over $k[\epsilon]$ are given as
$(f+g\epsilon)$ where  $f, g\in R$ and $f\neq 0$. 
We list some properties of such principal ideals that hold for any integral domain $R$. 
\begin{enumerate}
\item $(f+g_1\epsilon)=(f+g_2\epsilon)$ iff $g_1- g_2\in (f)$,
\item If $u\in R$ is a unit then so is $u+g\epsilon$ since $(u+g\epsilon)(u^{-1}-u^{-2}g\epsilon)=1$,
\item If $f$ is irreducible then so is $f+g\epsilon$ for every $g$,
\item  $(f+g\epsilon)(f-g\epsilon)=f^2$ shows that there is no unique factorization.
\item If the $f_i$ are pairwise relatively prime then 
$$
\tprod_i (f_i+g_i\epsilon)=\tprod_i (f_i+g'_i\epsilon)
\qtq{iff} (f_i+g_i\epsilon)= (f_i+g'_i\epsilon) \quad \forall i.
$$
\end{enumerate}
\end{exmp}

The following concrete example illustrates several of the above features.

\begin{exmp}[Picard group of a constant elliptic curve]\label{pic.const.ell.say}
Let $(0,E)$ be a smooth, projective elliptic curve. Over any base $S$ we have the constant family $\pi: E\times S\to S$ with the constant section
$s_0:S\cong \{0\}\times S$.  Let $L$ be a line bundle on $E\times S$. 
Then  $L\otimes \pi^*s_0^*L^{-1}$ has a canonical trivialization along
$\{0\}\times S$, hence it defines a morphism $S\to \pic(E)$. 
So the relative Picard group is computed by the formula
$$
\pic(E\times S/S)\cong \mor\bigl(S, \pic(E)\bigr).
\eqno{(\ref{pic.const.ell.say}.1)}
$$ 
Two consequences are worth mentioning.

\medskip
{\it Claim \ref{pic.const.ell.say}.2.} 
Let $(R, m)$ be a complete local ring.  Set $S=\spec R$  and $S_n=\spec R/m^n$.  Then
$$
\pic(E\times S/S)=\varprojlim \pic(E\times S_n/S_n).\qed
$$

\medskip
{\it Claim \ref{pic.const.ell.say}.3.} 
Let $S=\spec k[t]_{(t)}$ be the local ring of the affine line at the origin and
$\hat S=\spec k[[t]]$ its completion. Then 
$$
\pic(E\times S/S)\cong \pic(E) \qtq{but} \pic(E\times \hat S/\hat S)
\qtq{is infinite dimensional.} \qed
$$

Next consider the affine elliptic curve $E^\circ=E\setminus\{0\}$ and the
constant affine family $E^\circ\times S\to S$. Note that $\pic(E^\circ)\cong\pic^\circ(E)$. 

If $S$ is smooth and $D^\circ$ is a Cartier divisor on $E^\circ\times S$ then 
its closure $D\subset E\times S$ is also Cartier. More generally, this also holds if $S$ is normal, using \cite[4.21]{k-modbook}.  Thus (\ref{pic.const.ell.say}.1)
gives the following.

\medskip
{\it Claim \ref{pic.const.ell.say}.3.}  If $S$ is normal then 
$$
\pic(E^\circ\times S/S)\cong \mor\bigl(S, \pic^\circ(E)\bigr). \qed
$$

By contrast, (\ref{pic.over.artin.say}.3) gives the following.

\medskip
{\it Claim \ref{pic.const.ell.say}.4.}  If $S=\spec A$ is Artin  then 
$$\pic(E^\circ\times S/S)\cong \pic^\circ(E). \qed
$$

So $\pic(E^\circ\times S/S)$
has dimension 1 but
$\dim_k  \mor\bigl(S, \pic^\circ(E)\bigr)=\len A$.

The following is a good illustration.

\medskip
{\it Concrete Example \ref{pic.const.ell.say}.5.}
Start with the plane cubic  with equation $Y^2Z=X^3-Z^3$. In the affine plane $Z=1$ we get  $y^2=x^3-1$   (where  $x=X/Z, y=Y/Z$) and in the 
$Y=1$ plane we get  $v=u^3-v^3$  (where  $u=X/Y, v=Z/Y$). The diagonal in
$(y^2=x^3-1)\times (v=u^3-v^3)$ is a Cartier divisor which is defined by 2 equations   $yv=1$ and $ yu=x$.

At $(u=v=0)$ the local coordinate is $u$. Note that $u$ also vanishes at the points where $v^2+1=0$. If we invert it, then  we get that
$$
(u^{3r})=(v^r)\subset k\bigl[u,v,(v^2+1)^{-1}\bigr]/( u^3-v^3-v).
$$
What is the ideal
$$
(yv-1, yu-x, u^r)\subset k\bigl[x,y,u,v,(v^2+1)^{-1}\bigr]/(y^2-x^3+1,  u^3-v^3-v).
$$
Note that it contains
$$
(yv-1)(y^{r-1}v^{r-1}+\cdots+yv+1)=y^rv^r-1=y^r(v^2+1)^{-r}u^{3r}-1.
$$
Thus $1\in (yv-1, yu-x, u^r)$ and the ideal is the whole ring.
\end{exmp}

\subsection*{Relative Mumford divisors}{\ }

\begin{defn} \label{mcl.defn}
Let $f:X\to S$ be a flat morphism. Two relative Mumford divisors
$D_1, D_2\subset X$ are {\it linearly equivalent} if
$\o_X(-D_1)\cong \o_X(-D_2)$, and {\it linearly equivalent over $S$} if
$\o_X(-D_1)\cong \o_X(-D_2)\otimes f^*L$ for some line bundle $L$ on $S$. The linear equivalence classes over $S$ of relative Mumford divisors generate the {\it relative Mumford class group}  $\mcl(X/S)$. 

By definition, if $D$ is a Mumford divisor then there is a 
closed subset $Z\subset X$ such that  $\o_X(-D)|_{X\setminus Z}$ is locally free
and $\codim_{X_s}Z_s\geq 2$ for every $s\in S$.
This gives a natural identification 
$$
\mcl(X/S)=\textstyle{\lim}_Z \pic(X\setminus Z/S),
\eqno{(\ref{mcl.defn}.1)}
$$
where the limit is over all closed subsets $Z\subset X$ such that   $\codim_{X_s}Z_s\geq 2$ for every $s\in S$.

On a normal variety, a Mumford divisor is the same as a Weil divisor and the 
Mumford class group is the same as the  class group. 
If $f$ has normal fibers, then we get the 
{\it relative  class group}  $\cl(X/S):=\mcl(X/S)$. 

As with the Picard group, this may not be the optimal definition when $S$ is projective, but we will use this notion mostly when $S$ is local, and then this seems the right definition.
\end{defn}

\begin{prop} 
\label{divs.sheaves.dim.2.say.1}  Let   $(A, k)$ be a local Artin ring,
$k\cong (\epsilon)\subset A$ an ideal and $B=A/(\epsilon)$.
Let 
$(R_A, m)$ be a flat, local, $S_2$,   $A$-algebra  of dimension 2 and
set $X_A:=\spec_A R_A$. 
Let $f_B\in R_B$ be a non-zerodivisor and set $C_B:=(f_B=0)\subset X_B$.

Then the set of all relative Mumford divisors
 $D_A\subset X_A$  such that $\pure\bigl((D_A)|_B\bigr)=C_B$
is a torsor under  the infinite dimensional $k$-vector space $H^1_m(C_k, \o_{C_k})$. 
\end{prop}

Proof. We can lift $f_B$ to $f_A\in R_A$. 
Choose $y\in m$ that is not a zerodivisor on $C_B$ and such that
$D_A$ is a principal divisor on $X_A\setminus (y=0)$. After inverting $y$, we can write the ideal of $D$ as
$$
I_{(y)}=(f_A+\epsilon y^{-r}g_k)\qtq{where} g\in R_k, r\in\n.
\eqno{(\ref{divs.sheaves.dim.2.say.1}.1)}
$$
We can multiply $f_A+\epsilon y^{-r}g_k $ by  $u+\epsilon y^{-s}v$ where $u$ is a unit in $R$.
This changes $g_k$ to $ug_k+vy^{r-s}f_A$. Thus the relevant information is carried by
the residue class  
$$
\overline{y^{-r}g_k}\in  H^0(C_k^\circ, \o_{C_k^\circ}),
\eqno{(\ref{divs.sheaves.dim.2.say.1}.2)}
$$
where $C_k^\circ\subset C_k$ denotes the complement of the closed point.

If the residue class is in $H^0(C_k, \o_{C_k})$ then we get a Cartier divisor.
Thus the non-Cartier divisors are parametrized by
$$
H^0(C_k^\circ, \o_{C_k^\circ})/H^0(C_k, \o_{C_k})\cong H^1_m(C_k, \o_{C_k}).
\eqno{(\ref{divs.sheaves.dim.2.say.1}.3)}
$$
We compute in (\ref{divs.to.lb.loc.say}.2) that different elements of
$H^1_m(C_k, \o_{C_k}) $ give non-isomorphic divisors. \qed
\medskip

\begin{cor} \label{divs.sheaves.dim.2.cor} Let   $(A, k)$ be a local Artin ring,
$k\cong (\epsilon)\subset A$ an ideal and $B=A/(\epsilon)$.
Let 
$(R_A, m)$ be a flat, local, $S_2$,   $A$-algebra  of dimension 2. 
Let $f_A\in R_A$ and $g_k\in R_k$  be a non-zerodivisors, and
$y$ a non-zerodivisor modulo both $f_A$ and $g_k$. 

 For the divisorial ideal $I:=R_A\cap (f_A+\epsilon y^{-r}g_k)R_A[y^{-1}]$  the following are equivalent.
\begin{enumerate}
\item  $I$ is a principal ideal.
\item The  residue class $\overline{y^{-r}g_k} $ lies in  $R_k/(f_k)$.
\item  $g_k\in (f_k, y^r)$.
\end{enumerate}
\end{cor}

Proof. $I$ is a principal ideal iff it has a generator of the form
$f_A+\epsilon h_k$ where $h_k\in R_k$. This holds iff
$$
f_A+\epsilon y^{-r}g_k=(1+\epsilon y^{-s}b_k)(f_A+\epsilon h_k)\qtq{for some} b_k\in R_A.
$$
Equivalently, iff  $y^{-r}g_k=h_k+y^{-s}b_kf_k$.
If $r>s$ then $g_k=y^rh_k+y^{r-s}b_kf_k$ which is impossible since $y$ is not a zerodivisor modulo $g_k$. If $r<s$ then $y^{s-r}g_k=y^sh_k+b_kf_k$ which is impossible since $y$ is not a zerodivisor modulo $f_k$. Thus
$r=s$ and then $g_k=y^rh_k+b_kf_k$ is equivalent to $g_k\in (f_k, y^r)$.\qed

\begin{cor} \label{divs.sheaves.dim.2.cor.2} 
Using the notation of (\ref{divs.sheaves.dim.2.cor}), assume that $f_A-f'_A\in \epsilon m_R^N$ and $ g_k-g'_k\in m_R^N$
 some $N\gg 1$  (depending on $f_k, g_k$ and $r$). Then
$(f_A+\epsilon y^{-r}g_k)$ defines a relative Cartier divisor
iff $(f'_A+\epsilon y^{-r}g'_k)$ does.
\end{cor}

Proof. Choose $N$ such that $m^N\subset (f_k, y^r)$. 
Then $(f'_k, y^r)=(f_k, y^r)$ and $g_k-g'_k\in (f_k, y^r)$. \qed
\medskip

{\it Remark  \ref{divs.sheaves.dim.2.cor.2}.1.} If $R_k$ is regular then
we can choose $y$ to be a general element of $m\setminus m^2$. Then
 $\dim R_k/( f_k, y^r)=r\cdot \mult f_k$, so
$N=r\cdot \mult f_k$ works. 

\medskip

The connection between (\ref{divs.sheaves.dim.2.say.1}) and (\ref{pic.over.artin.say}) is given by the following.

\begin{say}\label{divs.to.lb.loc.say}
Let $X$ be an affine, $S_2$  scheme and
$D:=(s=0)\subset X$ a Cartier divisor. Let $Z\subset D$ be a closed subset that has codimension $\geq 2$ in $X$. Set
$X^\circ:=X\setminus Z$ and $D^\circ:=D\setminus Z$.
Restricting the exact sequence
$$
0\to \o_{X}\stackrel{s}{\to} \o_{X} \to \o_{D}\to 0
$$
to $X^\circ $ and taking cohomologies we get
$$
0\to H^0(X^\circ, \o_{X^\circ})\stackrel{s}{\to} H^0(X^\circ, \o_{X^\circ})\to H^0(D^\circ, \o_{D^\circ})\stackrel{\partial}{\to} H^1({X^\circ},\o_{X^\circ}).
$$
Note that $H^0(X^\circ, \o_{X^\circ})=H^0(X, \o_X)$ since $X$ is $S_2$ and its image in $H^0(D^\circ, \o_{D^\circ})$ is $ H^0(D, \o_{D}) $. Thus $\partial$ becomes the injection
$$
\partial\colon  H^1_Z(D, \o_D)\cong H^0(D^\circ, \o_{D^\circ})/ H^0(D, \o_{D})
{\into} H^2_Z(X, \o_X).
\eqno{(\ref{divs.to.lb.loc.say}.1)}
$$
We are especially interested in the case when $(x,X)$ is local, 2-dimensional and $Z=\{x\}$. In this case (\ref{divs.to.lb.loc.say}.1) becomes
$$
\partial\colon  H^1_x(D, \o_D)
{\into} H^2_x(X, \o_X).
\eqno{(\ref{divs.to.lb.loc.say}.2)}
$$\end{say}

We can be especially explicit in the smooth case.
(Note that,  by the Weierstrass preparation theorem, almost every curve in $\hat{\a}_{uv}$ is defined by a  monic polynomial in $v$.)

\begin{lem} \label{ci.bais.lem.5}  Let $f\in k[[u]][v]$ be a  monic polynomial in $v$ of degree $n$ defining a curve $C_k\subset \hat\a^2_{uv}$.  
Let  $D\subset \hat\a^2_{k[\epsilon]}$  be a relative Mumford divisor
 such that $\pure(D_k)=C_k$. 
Then the restriction of $D$ to the complement of $(u=0)$ 
can be uniquely written as
$$
f+\epsilon\tsum_{i=0}^{n-1} v^i\phi_i(u)=0\qtq{where}  \phi_i(u)\in u^{-1}k[u^{-1}].
$$
Thus the set of all such $D$
is naturally isomorphic to  the infinite dimensional $k$-vector space $H^1_m(C_k, \o_{C_k})\cong \oplus_{i=0}^{n-1} u^{-1}k[[u^{-1}]]$. 
\end{lem}

Proof.  Note that   $k[[u,v]]/(f)\cong \oplus_{i=0}^{n-1} v^ik[[u]]$ as a $k[[u]]$-module, 
giving isomorphism
$$
H^0(C_k, \o_{C_k})\cong \oplus_{i=0}^{n-1} v^ik[[u]]\qtq{and}
H^0(C_k^\circ, \o_{C_k^\circ})\cong \oplus_{i=0}^{n-1} v^ik((u)).
\eqno{(\ref{ci.bais.lem.5}.1)}
$$
That is,  if $g\in k((u))[v]$ is a
 polynomial of degree $<n$ in $v$ then $g|_{C^\circ}$ extends to  a regular function on $C$  iff $ g\in k[[u]][v]$. \qed

We can also restate (\ref{ci.bais.lem.5}.1)
 as
$$
H^1_m(C_k, \o_{C_k})\cong \oplus_{i=0}^{n-1} v^ik((u))/k[[[u]]
\cong \oplus_{i=0}^{n-1} v^iu^{-1}k[u^{-1}].
\eqno{(\ref{ci.bais.lem.5}.2)}
$$

\begin{exmp}  \label{divs.sheaves.dim.2.say.4} 
Consider next the special case of (\ref{ci.bais.lem.5}) when $f=v$.
We can then write the restriction of $D$ as
 $(v+\phi(u)\epsilon=0)$ where $\phi\in u^{-1}k[u^{-1}]$.
Let $r$ denote the pole-order of $\phi$ and set $q(u):=u^r\phi(u)$. 
\medskip

{\it Claim \ref{divs.sheaves.dim.2.say.4}.1.} The ideal of $D$ is
$$
I_D=\bigl(v^2, vu^r+q(u)\epsilon, v\epsilon\bigr).
$$
Thus the fiber over the closed point is
$k[[u,v]]/(v^2, vu^r)$. Its torsion submodule is isomorphic to  $k[[u,v]]/(v, u^r)\cong k[u]/(u^r)$. 
\medskip

Proof. To see this note first that
$v^2=(v+\phi(u)\epsilon)(v-\phi(u)\epsilon)$, 
$vu^r+q(u)=(v+\phi(u)\epsilon)u^r$ and 
$v\epsilon=(v+\phi(u)\epsilon)\epsilon$ are  elements of  $I_D$.
Next note that  $q(u)$ is a polynomial with nonzero constant term, hence invertible in $k[[u,v]] $. Therefore
$$
k[[u,v]][\epsilon]/\bigl(v^2, vu^r+q(u)\epsilon, v\epsilon\bigr)\cong
k[[u,v]]/\bigl(v^2,  v^2u^rq(u)^{-1}\bigr)=  k[[u,v]]/(v^2)
$$
has no embedded points. \qed
\medskip

The ideals of relative Mumford divisors in $k[[u,v]][\epsilon]$ are likely to be more complicated in general. At least the direct generalization of
(\ref{divs.sheaves.dim.2.say.4}.1) does not always give the correct generators.

For example,  let  $f=v^2-u^3$ and consider the ideal $I\subset k[[u,v]][\epsilon]$ extended from
 $\bigl((v^2-u^3)+u^{-3}v\epsilon\bigr)$.   The above procedure gives the elements
$$
 (v^2-u^3)^2, \   u^3(v^2-u^3)+v\epsilon, \ (v^2-u^3)\epsilon \in I.
$$
However, $u^3(v^2-u^3)+v\epsilon=v^2(v^2-u^3)+v\epsilon$ and we can cancel the $v$ to get that 
$$
I=\bigl( (v^2-u^3)^2,    v(v^2-u^3)+\epsilon, (v^2-u^3)\epsilon\bigr).
\eqno{(\ref{divs.sheaves.dim.2.say.4}.2)}
$$

\medskip
Using  the isomorphism
$R[\epsilon]/(f^2, fg+\epsilon, f\epsilon)\cong R/(f^2, -f^2g)\cong R/(f^2)$, the above examples can be generalized to the non-smooth case as follows.
\medskip

{\it Claim \ref{divs.sheaves.dim.2.say.4}.3.}   Let
$(R, m)$ be a local, $S_2$,   $k$-algebra  of dimension 2 and 
 $f,g\in m$  a system of parameters.  Then
$J_{f,g}=(f^2, fg+\epsilon, f\epsilon)$ is (the ideal of) a relative Mumford divisor in $R[\epsilon]$
whose central fiber is $R/(f^2, fg)$, with embedded subsheaf  isomorphic to $R/(f, g)$. \qed
\medskip

\end{exmp}

\section{Divisorial support}\label{divisorial.support.sect}

There are at least 3 ways to associate a divisor to a sheaf
(\ref{div.supp.defn.1}) but only one of them---the divisorial support---behaves well in flat families. In this Section we develop this notion and a method to compute it. The latter  is especially important for the applications.

\begin{defn}[Divisorial support of a sheaf]\label{div.supp.defn.1}
 Let $X$ be a scheme and $F$ a coherent sheaf on $X$.
One usually defines its {\it support} $\supp F$ and its
{\it scheme-theoretic support} $\ssupp F:=\spec_X (\o_X/\ann F)$.

Assume next that $X$ is regular at every generic point $x_i\in\supp F$ that has codimension 1 in $X$. Then there is a unique divisorial sheaf \cite[3.50]{k-modbook} associated to the Weil divisor  
$\tsum \len( F_{x_i})\cdot  [\bar x_i]$.
We call it the {\it divisorial support} of $F$ and denote it by $\dsupp F$. 

If every associated point of $F$ has codimension 1 in $X$ then we have
inclusions of subschemes
$$
\supp F\subset \ssupp F \subset \dsupp F.
\eqno{(\ref{div.supp.defn.1}.1)}
$$
In general all 3 subschemes are different, though with the same support.

Our aim is to develop a relative version of this notion and some ways of computing it in families. Let $X\to S$ be a morphism and $F$ a coherent sheaf on $X$.
Informally, we would like the relative divisorial support of $F$, denoted by $\dsupp_SF$,  to be a scheme over $S$ whose fibers are
$\dsupp (F_s)$ for all $s\in S$. If $S$ is reduced, this requirement uniquely determines $\dsupp_SF$ but in general there are 2 problems.
\begin{itemize}
\item Even in nice situations, this requirement may be impossible to meet.
\item  For non-reduced base schemes $S$, the fibers alone do not determine $\dsupp_SF$.
\end{itemize}
The right concept is developed through a series of Definition-Lemmas.
Each one is a definition, where we need to check that it is independent of the choices involved, and that it coincides with our naive definition over reduced schemes.
\end{defn}

We start with a very elementary case which, however,
turns out to be crucial.

\begin{defn-lem}[Divisorial support I]\label{char.pol.=.ds.1d}
Let  $C$ be a smooth curve and $M$  a torsion sheaf on $C$.
Thus it can be written as    $M\cong \oplus_j \o_C/\o_C(-n_jP_j)$
where $P_j\in C$ and $n_j\in \n$  (repetitions allowed). Then
$$
\dsupp(M)= \spec_C\bigl( \o_C/\o_C(-\tsum_jn_jP_j)\bigr).
\eqno{(\ref{char.pol.=.ds.1d}.1)}
$$
Let $\pi:C_1\to C_2$ be an \'etale  morphism of smooth curves and $M$  a torsion sheaf on $C_2$. Then we get that 
$$
\dsupp (\pi^*M)=\pi^*\dsupp (M).
\eqno{(\ref{char.pol.=.ds.1d}.2)}
$$
Thus the computation of   $\dsupp(M)$ is an \'etale-local question.
In order to develop another formula, we can work on $\a^1$. 
 Let $g_j\in k[x]$ be monic polynomials and set $M=\oplus_j k[x]/(g_j)$.
 Multiplication by $x$ is an endomorphism $\mu_x$ of the $k$-vector space $M$. We claim that its
characteristic polynomial is 
$$
\chi(\mu_x)(t)=(-1)^d\tprod_j g_j(t)\qtq{where} d=\deg \tprod_j g_j.
\eqno{(\ref{char.pol.=.ds.1d}.3)}
$$
Thus
$$
\dsupp M= \bigl(\tprod_jg_j=0\bigr)=\bigl(\chi(\mu_x)(x)=0\bigr).
\eqno{(\ref{char.pol.=.ds.1d}.4)}
$$
In particular, we could use any \'etale coordinate instead of $x$ in (\ref{char.pol.=.ds.1d}.4).

It is  enough to do check these  for 1 polynomial. 
Thus let $g=x^n+a_{n-1}x^{n-1}+\cdots+a_0$ be a monic polynomial.
In the module  $k[x]/(g)$ choose a $k$-basis
$e_i=x^i$ for $i=0,\dots, n-1$. Thus multiplication by $\mu_x$ is 
$$
\mu_x(e_i)= e_{i+1}\qtq{for $i<n-1$ and} \mu_x(e_{n-1})= -\tsum_i a_ie_i.
$$ 
We compute that its characteristic polynomial is
$(-1)^ng(x)$. 
For example, if $n=4$ then the matrix of  $\mu_x$ is  
$$
\left(\begin{array}{cccc}
0 & 0 & 0 & -a_0\\
1 & 0 & 0 & -a_1\\
0 & 1 & 0 & -a_2\\
0 & 0 & 1 & -a_3
\end{array}
\right).
$$
Expanding by the last column gives the characteristic polynomial
$$
\det \left(\begin{array}{cccc}
-x & 0 & 0 & -a_0\\
1 & -x & 0 & -a_1\\
0 & 1 & -x & -a_2\\
0 & 0 & 1 & -x-a_3
\end{array}
\right)
=(-1)^4\bigl(x^4+a_3x^2+\cdots+a_0\bigr). \qed
$$
\end{defn-lem}

\begin{defn-lem}[Divisorial support II]\label{dsupp.car.1d.dl}
Let $X\to \a^1_S$ be an \'etale morphism.
Let $F$ be a coherent sheaf on $X$ that is finite and flat (hence locally free) over $S$. Let $t$ denote a coordinate on $\a^1_S$.
Multiplication by $t$ is an endomorphism $\mu_t$ of the locally free $\o_S$-module $F$;
let   $\chi(\mu_t)(*)$ be  its characteristic polynomial.
Then the divisorial support  of $F$ over $S$ is 
$$
\dsupp_S(F)=  \bigl(\chi(\mu_t)(t)=0\bigr)\subset X.
$$
It is a relative Cartier divisor.
\end{defn-lem}

Proof of consistence.  We need to show that this is independent of the choice of $t$.
This is an \'etale-local question, 
We may thus assume that $S$ is the spectrum of a Henselian local ring  $(R, m, k)$, $X$ is the spectrum of
 $R\langle t\rangle$ (the Henselization of $R[t]$)  and $F$ is the sheafification of the free $R$-module $M$.

Multiplication by $t$ is an endomorphism of  $M$, its characteristic polynomial is the defining equation of $\dsupp(M)$.
However, the choice of $t$ is not unique; we could have used any other local parameter
$t'=a_1t+a_2t^2+\cdots$ where $a_1\notin m$. 

Assume first that $S$ is reduced. Then the independence of $t$ can be checked over the generic fibers, where we recover the computation of 
(\ref{char.pol.=.ds.1d}). 

If $S$ is arbitrary, then we use that 
 the pair  $\bigl(M, \mu_t:M\to M\bigr)$ is   induced from the universal endomorphism $\mu^u$ of $M_S:=\oplus_i R^ue_i$ 
over the Henselian ring $R^u:=k\langle t_{ij}:1\leq i,j\leq r\rangle$ given by 
$$
(e_1,\dots, e_r)^{\rm tr}\mapsto  \bigl(t_{ij}\bigr)\cdot (e_1,\dots, e_r)^{\rm tr}.
$$
We already noted that independence of $t$ holds over $\spec R^u$,
hence it also holds  after pulling back to $S$. \qed

\begin{defn-lem}[Divisorial support III] \label{dsupp.car.nd.dl}
Let $X\to S$ be a smooth morphism of pure relative dimension $n$.
Let $F$ be a coherent sheaf on $X$ that is flat over $S$  with CM
 fibers    of pure dimension $n-1$.
 Then its divisorial support   $\dsupp_S(F)$ is defined and it is relatively Cartier over $S$.
\end{defn-lem}

Proof  of consistence. 
Being relatively Cartier can be checked \'etale-locally.
Thus we may assume that $S$ is the spectrum of  a Henselian local ring  $(R, m, k)$ with infinite residue field and   $M$ is a finite  $R\langle x_1,\dots, x_n\rangle $-module that is flat over $R$ and
such that $M_k$ is CM and  of dimension $n-1$.

After a linear coordinate change, we may also assume that
 $M_k$ is finite over $k\langle  x_1,\dots, x_{n-1}\rangle $. Since $M_k$ is CM, it is also free, thus $M$ is also free over $R\langle x_1,\dots, x_{n-1}\rangle$.

Multiplication by $x_n$ is an endomorphism of the free $R\langle x_1,\dots, x_{n-1}\rangle$-module $M$, its characteristic polynomial is the defining equation of $\dsupp(M)$.

It remains to show that this equation is 
independent of the choices that we made, up to a unit. On a scheme two 
locally principal  divisors agree iff they agree at their generic points. 

Thus after further localization at a generic point of $\supp M$ and flat base change  \cite[10.47]{k-modbook}, 
we may assume that the  relative dimension is 1. 
This was already treated in (\ref{dsupp.car.1d.dl}). \qed
\medskip

The following  properties are especially important.

\begin{cor} \label{dsupp.rets.lem.3-2} 
Continuing with the  notation and assumptions of (\ref{dsupp.car.nd.dl}), let $h:S'\to S$ be a  morphism. By base change we get
$g':X'\to S'$ and $h_X:X'\to X$.  Then
$$
h_X^*(\dsupp F)=\dsupp (h_X^*F).
$$
\end{cor}

Proof.  This is an \'etale-local question on $X$, thus, as in the proof of
(\ref{dsupp.car.nd.dl}) we may assume that  $M$ is  free over $R\langle x_1,\dots, x_{n-1}\rangle$.  The base change  $S'\to S$ corresponds to a ring extension
$R\to R'$, and the characteristic polynomial commutes with ring extensions.
\qed

\begin{cor} \label{dsupp.rets.lem.3-1} 
Let $X\to S$ be a smooth morphism of pure relative dimension $n$.
Let $F$ be a coherent sheaf on $X$ that is flat over $S$  
 fibers    of pure dimension $n-1$.
 Then $\dsupp_SF$ is a relative Cartier divisor.
\end{cor}

Proof. Unlike in (\ref{dsupp.car.nd.dl}), we do not assume that the fibers are CM. However, if $x\in X_s$ is a point of codimension $\leq 2$, then $F_s$ is CM at $x$, hence $\dsupp_SF$ is a relative Cartier divisor at $x$ by (\ref{dsupp.car.nd.dl}). Since $X\to S$ is smooth, $\dsupp_SF$ is a relative Cartier divisor everywhere by (\ref{pic.codim3.cor}). \qed

\begin{cor} \label{dsupp.rets.lem.3} 
Continuing with the  notation and assumptions of (\ref{dsupp.car.nd.dl}), let $D\subset X$ be a relative Cartier divisor that is also  smooth over $S$. Assume that $D$ does not contain any
 generic point of   $\supp F_s$ for any $s\in S$. Then
$$
\dsupp (F|_D)=(\dsupp F)|_D.
$$
\end{cor}

Proof. As in the  proof of  (\ref{dsupp.car.nd.dl}),
we can choose  local coordinates such that $D=(x_1=0)$. 
Then $\dsupp F $ is computed from the characteristic polynomial of
$M$ over $R\langle x_1,\dots, x_{n-1}\rangle$ and 
$\dsupp (F|_D) $ is computed from the characteristic polynomial of
$M/x_1M$ over $R\langle x_2,\dots, x_{n-1}\rangle$. \qed

\medskip

Now we are ready to define the sheaves for which the
relative divisorial support makes sense, but first we have to distinguish associated points that come from the base from the other ones.

\begin{defn}\label{vpure.defn}
 Let $g:X\to S$ be  a morphism and $F$ a coherent sheaf on $X$ such that 
$\supp F\to S$ has pure relative dimension $d$.
An associated point $x\in \ass(F)$ is called {\it vertical}
if $x$ is not a generic point of  $\supp (F_{g(x)})$. 

We say that $F$ is {\it vertically pure} if it has no vertical
associated points. 

If $F$ is  generically flat over $S$ (\ref{pure.over.defn}), then there is a unique largest subsheaf  $\vtors_S(F)\subset F$---called the {\it  vertical torsion} of $F$---such that every fiber of the structure map $\supp(\vtors_S(F))\to S$ has dimension $<d$. 

Then  $\vpure (F):=F/\vtors_S(F)$ has no vertical associated primes.

All these notions make sense for subschemes of $X$ as well.
\end{defn}

\begin{defn}\label{pure.over.defn}
 Let $X\to S$ be  a morphism and $F$ a coherent sheaf on $X$.
We say that $F$ is a {\it generically flat family of pure sheaves}  of dimension $d$ over $S$ if the following hold.
\begin{enumerate}
\item $F$ is flat at every generic point of $F_s$ for every $s\in S$ and 
\item $\supp F\to S$ has pure relative dimension $d$.
\end{enumerate}
We usually do not care about  vertical
associated points on $F$, thus we frequently replace  $F$ by $\vpure (F)=F/\vtors_S(F)$   and then the following condition is also satisfied.
\begin{enumerate}\setcounter{enumi}{2}
\item  $F$ is vertically pure. 
\end{enumerate}
We say that $Z\subset X$ is a {\it generically flat family of pure subschemes}
if its structure sheaf $\o_Z$ has this property.

The following properties are clear from the definition.
\begin{enumerate}\setcounter{enumi}{3}
\item  Conditions (1--2) are preserved by any base change $S'\to S$ and
(3) is preserved by flat base change.
\item If (3) holds  then  the generic fibers $F_g$ are pure of dimension $d$, but special fibers
may have embedded points outside $\fcm(F)$ (\ref{flat.locus.defn}). 
\end{enumerate}
\end{defn}

\begin{defn}\label{flat.locus.defn}
 Let $X\to S$ be a morphism and  $F$ a coherent sheaf on $X$.
The {\it flat locus} of $F$ is the largest open subset  $U\subset \supp F$ such that $F|_U$ is flat over $S$. We denote it by
$\floc (F)$.

It is sometimes more convenient to work with the 
 {\it flat-CM locus} of $F$. It is the largest open subset  $U\subset \supp F$ such that $F|_U$ is flat with CM fibers over $S$. We denote it by
$\fcm (F)$.

These properties are unchanged if we replace $X$ by
$\ssupp F$. Thus we may assume that $\supp F=X$, or even that $\ann(F)=0$, whenever it is convenient.
\end{defn}

\begin{defn-lem}[Divisorial support IV] \label{dsupp.car.nd.dl.4}
Let $g:X\to S$ be a flat morphism of pure relative dimension $n$
and  $g^\circ: X^{\circ}\to S$ the smooth locus of $g$.

Let $F$ be a coherent sheaf on $X$ that is   generically flat and pure over $S$ of dimension $n-1$.
Assume that for every $s\in S$, every generic point of $F_s$ is contained in 
$X^{\circ}$.

Set $Z:=\supp F\setminus \bigl(\fcm(F)\cap X^{\circ}\bigr)$,  $U:=X\setminus Z$ and $j:U\into X$ the natural injection.  We define the
{\it divisorial support} of $F$ over $S$ as 
$$
\dsupp_S(F):=\overline{\dsupp_S(F|_U)},
\eqno{(\ref{dsupp.car.nd.dl.4}.1)}
$$
the scheme-theoretic closure of $\dsupp_S(F|_U) $. This makes sense since
 the latter is already defined by (\ref{dsupp.car.nd.dl}).

Note that $\dsupp_S(F)$ is a generically flat family of pure subschemes
 of dimension $n-1$  over $S$ and it is 
 relatively Cartier on $U$.

It  is enough to check the following equalities  at codimension 1 points, which  follow from (\ref{char.pol.=.ds.1d}).

\medskip
{\it Claim \ref{dsupp.car.nd.dl.4}.2.}  
 Let $X_i\to S$ be   flat morphisms of pure relative dimension $n$ and
 $\pi:X_1\to X_2$ a finite morphism.
Let $D\subset X_1$ be a relative Mumford divisor. 
Assume  that
 $\red D_s\to \red(\pi(D_s))$ is birational and 
 $\pi$ is \'etale at every generic point of $D_s$. Then
$$
\dsupp(\pi_*\o_D)=\pi(D), \qtq{the scheme-theoretic image of $D$.} \qed
$$

\medskip
{\it Claim \ref{dsupp.car.nd.dl.4}.3.}  Let $X_i\to S$ be   flat morphisms of pure relative dimension $n$ and
 $\pi:X_1\to X_2$ a finite morphism. Let
$F_1$ be a coherent sheaf on $X_1$ that is   generically flat and pure over $S$ of dimension $n-1$. Set $F_2:=\pi_*F_1$. 
Assume that $g_i$ is smooth at  every generic point of $(F_i)_s$ for every $s\in S$. Then
$$
\dsupp_S(\pi_*F_1)=\dsupp_S\bigl(\pi_*\dsupp_S(F_1)\bigr). \qed
$$

\end{defn-lem}

The next claim directly  follows from \cite[4.36]{k-modbook}.

\begin{lem} \label{dsupp.red.sm.lem}
Let $S$ be a reduced scheme and $g:X\to S$  a smooth  morphism of pure relative dimension $n$
Let $F$ be a coherent sheaf on $X$ that is   generically flat and pure over $S$ of dimension $n-1$. Then $\dsupp_S(F)$ is a relative Cartier divisor. \qed
\end{lem}

Divisorial support commutes with restriction to a divisor, whenever
everything makes sense.  We just need to make enough assumptions that guarantee that 
(\ref{dsupp.rets.lem.3}) applies on a dense set of every fiber.

\begin{cor} \label{dsupp.rets.lem.4}
Continuing with the  notation and assumptions of (\ref{dsupp.car.nd.dl.4}), let  $D\subset X$ be a relative Cartier divisor. Assume that there is an open set $D^\circ\subset D$ such that
\begin{enumerate}
\item $g|_D$ is smooth on $D^\circ$,
\item $D^\circ_s$ is dense in $ D_s$ for every $s\in S$,
\item $D$ does not contain any generic point of   $\supp F_s$ for any $s\in S$, and  
\item   $D^\circ\subset \fcm (F)$. 
\end{enumerate}
 Then
$$
\dsupp (F|_D)=\vpure\bigl((\dsupp F)|_D\bigr).\qed
$$
\end{cor}

Various Bertini-type theorems show that the  above assumptions are  quite easy to satisfy, at least locally. 

\begin{cor} \label{dsupp.rets.lem.lem4.c}  Continuing with the above notation,
let $|D|$ be a linear system  on $X$   that is base point free in characteristic 0 and very ample in general. 
Fix  $s\in S$ and let $D\in |D|$ be a general member. 
Then there is an open neighborhood  $s\in S^\circ\subset S$ such that 
 $$
\dsupp (F|_D)=(\dsupp F)|_D  \qtq{holds over $S^\circ$.}
\eqno{(\ref{dsupp.rets.lem.lem4.c}.1)}
$$
\end{cor} 

Proof. We apply the usual Bertini theorems  to  $X_s$. We get that
$D_s$ satisfies conditions (\ref{dsupp.rets.lem.4}.1--4), and then they also hold over some open neighborhood  $s\in S^\circ\subset S$. 

This gives (\ref{dsupp.rets.lem.lem4.c}.1),  modulo 
vertical torsion. Finally note that  there is no such torsion for general $D$ by \cite[10.9]{k-modbook}. \qed

\begin{lem} \label{div.supp.flat.bc.lem} Divisorial support commutes with  base change. That is,
let $g:X\to S$ be a flat morphism of pure relative dimension $n$ and
$F$  a generically flat family of pure sheaves  of dimension $n-1$ over $S$. 
Assume that for every $s\in S$, every generic point of $\supp F_s$ is contained in 
the smooth locus of $g$.

Let $h:S'\to S$ be a morphism. By base change we get
$g':X'\to S'$ and $h_X:X'\to X$.  Then
$$
h_X^{[*]}(\dsupp F)=\dsupp (h_X^*F).
$$
\end{lem}

Proof.  Set $U:=\fcm (F) \subset X$ with injection $j:U\into X$. 
Set $U':=h_X^{-1}(U)$  and $h_U:U'\to U$ the restriction of $h_X$. Then 
$h_U^*(\dsupp F|_U)=\dsupp \bigl(h_U^*(F|_U)\bigr)$ by (\ref{dsupp.rets.lem.3-2}).  

By (\ref{pure.over.defn}.4)  $h_X^{[*]}(\dsupp F) $ is a generically flat family of pure divisors and it agrees with $ \dsupp (h_X^*F)$ over $U'$. Thus the 2 are equal. \qed

\begin{defn}[Divisorial support of cycles]\label{div.supp.cycles.defn}
Let $S$ be a seminormal scheme  and  $Z$ a well defined family of $d$-cycles on $\p^n_S$ as in \cite[I.3.10]{rc-book}. 

Let $\rho:\supp Z\to \p^{d+1}_S$ be a finite morphism. 
Then $\rho_*Z$ is a  well defined family of $d$-cycles on $\p^{d+1}_S$.

If all the residue characteristics are 0, or if $Z$ satisfies the field of definition condition  \cite[I.4.7]{rc-book}, then
there is a unique relative Cartier divisor  $D\subset \p^{d+1}_S$
whose associated cycle is $\rho_*Z$; see \cite[I.3.23.2]{rc-book}. We denote it by   $\dsupp(\rho_*Z)$.  As a practical matter, we usually think of
 $\rho_*Z$ and $\dsupp(\rho_*Z)$ as the same object.

Let  $F$ be a coherent sheaf on $\p^n_S$ that is   generically flat and pure over $S$ of dimension $d$. One can   associate to it a  cycle $Z(F)$ that is a well defined family of $d$-cycles over $S$  (cf.\ \cite[I.3.15]{rc-book}). Let $\rho:\supp F\to \p^{d+1}_S$ be a finite morphism.
As in  (\ref{dsupp.car.nd.dl.4}.3) we get that
$$
\dsupp_S(\rho_*F)=\dsupp_S\bigl(\rho_*Z(F)\bigr).
\eqno{(\ref{div.supp.cycles.defn}.1)}
$$
Thus $\dsupp_S(\rho_*F) $ can be defined entirely in terms of cycles.
Note, however, that here the right hand side is defined only for seminormal schemes. 

One of the main aims defining the divisorial support for sheaves is to be able to work over arbitrary schemes.

\end{defn}

\section{Variants of K-flatness}\label{K.flat.sect}

We introduce 5 versions of K-flatness, which may well be equivalent to each other. From the technical point of view Cayley-Chow-flatness (or C-flatness)
is the easiest to use, but a priori it depends on the choice of a projective embedding. Then most of the work in the next 2 sections  goes to proving that a modified version  (stable C-flatness) is equivalent to K-flatness, hence independent of the projective embedding.

\begin{say}[Projections of $\p^n$]\label{proj.formulas.exmp}
Projecting $\p^{n}_{\mathbf x}$ from the point $(a_0 : \cdots : a_{n})$
to the $(x_{n}=0)$ hyperplane is given by
$$
\pi: (x_0 : \cdots : x_{n}) \to (a_{n}x_0-a_0x_{n} : \cdots : 
a_{n}x_{n-1}-a_{n-1}x_{n}).
\eqno{(\ref{proj.formulas.exmp}.1)}
$$
It is convenient to normalize $a_{n}=1$ and then we get
$$
\pi: (x_0 : \cdots : x_{n}) \to (x_0-a_0x_{n} : \cdots : 
x_{n-1}-a_{n-1}x_{n}).
\eqno{(\ref{proj.formulas.exmp}.2)}
$$
Similarly, a Zariski open set of projections of $\p^{n}_{\mathbf x}$ 
to $L^r=(x_n=\cdots=x_{r+1}=0)$ is given by 
$$
\pi: (x_0 : \cdots : x_{n}) \to \bigl(x_0-\ell_0(x_{r+1},\dots,x_{n}) : \cdots : 
x_{r}-\ell_r(x_{r+1},\dots,x_{n})\bigr),
\eqno{(\ref{proj.formulas.exmp}.3)}
$$
where the $\ell_i$ are linear forms.

Note that in affine coordinates, when we set $x_0=1$, the projections become
$$
\pi: (x_1,\dots,x_{n}) \to \Bigl(\frac{x_{1}-\ell_1}{1-\ell_0},\dots,
\frac{x_{r}-\ell_r}{1-\ell_0}\Bigr),
\eqno{(\ref{proj.formulas.exmp}.4)}
$$
where again the  $\ell_i$ are (homogeneous) linear forms in the $x_{r+1},\dots,x_{n}$. The coordinate functions  have a non-linear expansion
$$
\frac{x_{i}-\ell_i}{1-\ell_0}=(x_{i}-\ell_i)(1+\ell_0+\ell_0^2+\cdots ).
\eqno{(\ref{proj.formulas.exmp}.5)}
$$
Finally, non-linear projections are given as
$$
\pi: (x_1,\dots,x_{n}) \to \bigl(x_{1}-\phi_1(x_{1},\dots,x_{n}),\dots,
x_{r}-\phi_r(x_{1},\dots,x_{n})\bigr),
\eqno{(\ref{proj.formulas.exmp}.6)}
$$
where $\phi_i(x_1,\dots,x_{r}, 0,\dots, 0)\equiv 0$ for every $i$.
\end{say}

\begin{defn}\label{lin.projection.defn}  Let $E$ be a vector bundle over a scheme $S$ and $F\subset E$ a vector subbundle. This induces a natural 
{\it linear projection} map
$\pi:\p_S(E)\map \p_S(F)$. If $S$ is local then $E, F$ are free. After choosing bases, $\pi$ is given by a matrix of constant rank with entries in  $\o_S$. 
We call these {\it $\o_S$-projections} if we want to emphasize this.
If $S$ is over a field $k$, we can also consider {\it $k$-projections,} 
 given by a matrix with entries in $k$. These, however, only make good sense if we have a canonical trivialization of $E$; this rarely happens for us.

\end{defn}

We can now formulate  various versions of K-flatness and their basic relationships.

\begin{defn} \label{CK.flat.verss.defns}
Let $(s,S)$ first be a local scheme with infinite residue field and $F$ a
 generically flat family of  pure, coherent sheaves of relative dimension $d$  on $\p^n_S$   (\ref{pure.over.defn}), with scheme-theoretic support  $Y:=\ssupp F$.  
\begin{enumerate}
\item  $F$ is {\it C-flat} over $S$ iff $\dsupp (\pi_*F)$ is Cartier over $S$ for every $\o_S$-projection  $\pi: \p^n_S\map \p^{d+1}_S$ (\ref{lin.projection.defn})  that is finite on $Y$.
\item $F$ is {\it stably C-flat} iff $(v_m)_*F$ is C-flat for every
Veronese embedding  $v_m: \p^n_S\into \p^N_S$ (where $N=\binom{n+m}{n}-1$).
\item  $F$ is {\it K-flat} over $S$ iff $\dsupp (\rho_*F)$ is Cartier over $S$ for every finite morphism  $\rho: Y\to \p^{d+1}_S$.
\item  $F$ is {\it locally K-flat} over $S$ at $y\in Y$ iff $\dsupp (\rho_*F)$ is Cartier over $S$ at $\rho(y)$ for every finite morphism  $\rho: Y\to \p^{d+1}_S$ 
for which  $\{y\}=\rho^{-1}(\rho(y))$.
\item $F$ is {\it formally K-flat} over $S$ at a closed point  $y\in Y$  iff $\dsupp (\rho_*\hat F)$ is Cartier over $\hat S$  for every finite morphism  $\rho: \hat Y\to \hat{\a}^{d+1}_{\hat S}$,  where $\hat S $  (resp.\  $\hat Y $) denotes the completion of   $S$ at $s$  (resp.\ $Y$ at $y$).
\end{enumerate}

{\it Base change properties \ref{CK.flat.verss.defns}.6.} 
We see in (\ref{CC.flat.4.vers.thm.cor.1}) that being C-flat is
 preserved by arbitrary base changes
and  the property descends from faithfully flat base changes. This then implies the same for stable C-flatness. Once we prove that the latter is equivalent to 
K-flatness, the latter also has the same base change properties. Most likely the same holds for formal K-flatness.
\medskip

{\it General base schemes \ref{CK.flat.verss.defns}.7.} 
We say that any of  the above 
 notions (1--3) holds for a local base scheme $(s,S)$ (with finite residue field) 
if it holds after some  faithfully flat base change
$(s',S')\to (s, S)$, where $k(s')$ is infinite. Property (6) assures 
that this is independent of the choice of $S'$.

Finally we say that any of  the above 
 notions (1--3) holds for an  arbitrary base scheme $S$
if it holds for all its localizations. 

\medskip

{\it Comment on the notation \ref{CK.flat.verss.defns}.8.} 
Here C stands for the initial of either Cayley or Chow and, as before,  K stands for the first syllable of Cayley. 
\end{defn}

\begin{variants} \label{CK.flat.verss.vars}  These definitions each have other versions and relatives. I believe that each of the above 5 are natural and maybe even optimal, though they may not be stated in the cleanest form.
Here are some other possibilities and equivalent versions.
\begin{enumerate}
\item   It could have been better to define C-flatness using the Cayley-Chow form; the equivalence is proved in (\ref{CC.flat.4.vers.thm}). 
The Cayley-Chow form version matches better with the study of Chow varieties;  the definition in (\ref{CK.flat.verss.defns}.1)  emphasizes the similarity with the other 4.
\item  In (\ref{CK.flat.verss.defns}.3) we get an equivalent notion if we allow all  finite morphisms  $\rho: Y\to W$,   where $W\to S$ is any smooth, projective morphism of pure relative dimension $d+1$ over $S$.
Indeed, let $\pi:W\to \p^{d+1}_S$ be a finite morphism.  If $F$ is K-flat then
$\dsupp \bigl((\pi\circ \rho)_*F\bigr)$ is a relative Cartier divisor, hence
$\dsupp ( \rho_*F)$ is $F$ is K-flat by (\ref{dsupp.car.nd.dl.4}.3). Since $W\to S$ is smooth,
$\dsupp ( \rho_*F)$ is a relative Cartier divisor by (\ref{mumf.d.Cf.cor.2}). 
\item   It would be natural to consider an  affine version of $C$-flatness:
We start with a coherent sheaf $F$ on $\a^n_S$ and require that 
 $\dsupp (\pi_*F)$ be Cartier over $S$ for every projection  $\pi: \a^n_S\to \a^{d+1}_S$ that is finite on $Y$. 

The problem is that 
the relative affine version of  Noether's normalization theorem does not hold, thus there may not be any such projections; see  \cite[10.73.7]{k-modbook}.
This is why 
(\ref{CK.flat.verss.defns}.4) is stated for projective morphisms only.

Nonetheless, the notions (\ref{CK.flat.verss.defns}.1--4) are \'etale local on $X$ and most likely the following Henselian version of (\ref{CK.flat.verss.defns}.5) does work.
\item Assume that $f:(y,Y)\to (s,S)$ is a local morphism of pure relative dimension $d$ of Henselian local schemes such that  $k(y)/k(s)$ is finite. Let  $F$ be  a coherent sheaf on $X$ that is  pure of relative dimension $d$ over $S$.
Then  $F$ is {\it  K-flat} over $S$  iff $\dsupp (\rho_*F)$ is Cartier over $S$ for every  finite morphism  $\rho: Y\to \spec \o_S\langle  x_0,\dots, x_d\rangle$
(where $ R\langle  \mathbf{x} \rangle$ denotes the Henselization of 
$ R[\mathbf{x} ]$).

\end{enumerate}

\end{variants}

It is possible that in fact all 5 versions (\ref{CK.flat.verss.defns}.1--5) are equivalent to each other, but for now we  can prove only 8 of the 10 possible implications.
Four of them  are  easy to see.

\begin{prop} \label{CK.flat.verss.imps} Let $F$  be a  generically flat family of  pure, coherent sheaves of relative dimension $d$  on $\p^n_S$. 
Then
$$\mbox{formally K-flat} \Rightarrow \mbox{K-flat}\Rightarrow \mbox{locally K-flat} \Rightarrow \mbox{stably C-flat}\Rightarrow \mbox{C-flat}.
$$
\end{prop}

Proof. A divisor $D$ on a scheme $X$ is Cartier iff its completion $\hat D$ is Cartier on $\hat X$ for every $x\in X$ by (\ref{cart.infinit.cond.thm}). Thus 
$\mbox{formally K-flat} \Rightarrow \mbox{K-flat}$. 

$\mbox{K-flat}\Rightarrow \mbox{locally K-flat}$ is clear and
$\mbox{locally K-flat} \Rightarrow  \mbox{stably C-flat}$  follows from
(\ref{loc.C.flat.C.flat.lem}).  Finally $\mbox{stably C-flat}\Rightarrow \mbox{C-flat}$ is proved in (\ref{bert.sec.lem.1}). \qed
\medskip

A key technical  result of the paper is the following, proved in (\ref{K-flat.=.sC-flat.thm.pf}).

\begin{thm} \label{K-flat.=.sC-flat.thm} K-flatness is equivalent to stable C-flatness.
\end{thm}

It is quite likely that  our methods will prove the following.

\begin{conj} \label{fKf.=.Kf.conj}
Formal K-flatness is equivalent to K-flatness.
\end{conj}

We prove the special case of  relative dimension 1 in (\ref{d=1.scf=kf=ff.prop}); this is also a key step in the proof of (\ref{K-flat.=.sC-flat.thm}).

The remaining question is whether C-flat implies stably C-flat.
This holds in  the examples that I computed, but I have not been able to compute many  and I do not have any conceptual argument why 
these 2 notions should be equivalent. 
See also (\ref{id.ch.eq.defn}) for a related question about the ideal of Chow equations.

\begin{ques}  \label{CK.flat.verss.imps.4} 
Is C-flatness  equivalent to stable C-flatness?
\end{ques}

All of the above properties are automatic over reduced schemes and they can be checked on Artin subschemes.

\begin{prop}\label{over.red.K.autom.ptop}
Let $S$ be a  reduced scheme and $F$ a 
 generically flat family of  pure, coherent sheaves of relative dimension $d$  on $\p^n_S$.  Then $F$ is K-flat over $S$.
\end{prop}

Proof. This follows from (\ref{dsupp.red.sm.lem}).\qed

\begin{prop}\label{over.art.enough.ptop}
Let $S$ be a  scheme and $F$ a 
 generically flat family of  pure, coherent sheaves of relative dimension $d$  on $\p^n_S$.
Then $F$ satisfies one of the properties (\ref{CK.flat.verss.defns}.1--5) iff   $\tau^*F$ satisfies  the same property
for every Artin subscheme $\tau:A\into S$.
\end{prop}

Proof. Let $\pi:X\to \p^{d+1}_S$ be a finite morphism. 
By  (\ref{cart.infinit.cond.thm}) $\dsupp_S(\pi_*F)$ is Cartier iff
$\dsupp_A\bigl((\pi_A)_*\tau^*F\bigr)$ is Cartier for every Artin subscheme $\tau:A\into S$.  Thus the Artin versions imply the global ones in all cases.

To check the converse, we may localize at $\tau(A)$.  The claim is clear  if every finite morphism
$\pi_A:X_A\to \p^{d+1}_A$ can be extended to $\pi:X\to \p^{d+1}_S$.
This is obvious for C-flatness and stable C-flatness, but it need not hold
for K-flatness.

However, we see in (\ref{bert.sec.lem.1-1}) that it is enough to extend 
it after composition with a high enough
Veronese embedding. \qed

\section{Cayley-Chow flatness}\label{CC.flat.sect}

Let $Z\subset \p^n$ be a subvariety of dimension $d$. Cayley  \cite{cayley1860, cayley1862}\footnote{The titles of these articles are identical.} associates to it a hypersurface
$$
\chh(Z):=\{L\in \grass(n{-}d{-}1,\p^n)\colon Z\cap L\neq \emptyset\}\subset \grass(n{-}d{-}1,\p^n),
$$
called the  Cayley or Chow hypersurface; its equation is called the
Cayley or Chow form. 

We extend this definition to coherent sheaves on $\p^n_S$ over an arbitrary base scheme.  We use 2 variants, but the  proof of (\ref{CC.flat.4.vers.thm})  needs 2 other versions as well. All of these are defined in the same way, but
$\grass(n{-}d{-}1,\p^n) $ is replaced by other universal varieties.

\begin{defn}[Cayley-Chow hypersurfaces]\label{CC.defn.4.vers}
Let $S$ be a  scheme and $F$ a  generically flat family of pure, coherent sheaves  of dimension $d$ on $\p^n_S$  (\ref{pure.over.defn}). 
We define 4 versions of the Cayley-Chow hypersurface associated to $F$ as follows.

In all 4 versions the left hand side map $\sigma$ is a smooth fiber bundle.

\medskip
{\it Grassmannian version \ref{CC.defn.4.vers}.1.} Consider the diagram
$$
\begin{array}{lcl}
& \flag_S\bigl((\mbox{point}), n{-}d{-}1, \p^n\bigr)&\\
\sigma_g\swarrow && \searrow\pi_g\\[1ex]
\p^n_S && \grass_S(n{-}d{-}1, \p^n)
\end{array}
$$
where the flag variety parametrizes pairs  $(\mbox{point})\in L^{n-d-1}\subset \p^n$.  Set
$$
\chh_g(F):=\dsupp_S\bigl((\pi_g)_*\sigma_g^*F\bigr).
$$

\medskip
{\it Product version \ref{CC.defn.4.vers}.2.} Consider the diagram
$$
\begin{array}{lcl}
& \icorr_S\bigl((\mbox{point}), (\check{\p}^n)^{d+1}\bigr)&\\
\sigma_p\swarrow && \searrow\pi_p\\[1ex]
\p^n_S && (\check{\p}^n)^{d+1}_S
\end{array}
$$
where the incidence variety parametrizes   $(d+2)$-tuples
$\bigl((\mbox{point}), H_0,\dots, H_d\bigr)$ satisfying
$(\mbox{point})\in H_i$ for every $i$.  Set
$$
\chh_p(F):=\dsupp_S\bigl((\pi_p)_*\sigma_p^*F\bigr).
$$

\medskip
{\it Flag version \ref{CC.defn.4.vers}.3.} Consider the diagram
$$
\begin{array}{lcl}
& \operatorname{PFlag}_S(0, n{-}d{-}2, n{-}d{-}1,\p^n)&\\
\sigma_f\swarrow && \searrow\pi_f\\[1ex]
\p^n_S && \flag_S(n{-}d{-}2, n{-}d{-}1, \p^n)
\end{array}
$$
where the $\operatorname{PFlag}$  parametrizes triples  
$\bigl((\mbox{point}),  L^{n-d-2},  L^{n-d-1}\bigr)$  such that
$(\mbox{point})\in L^{n-d-1}$ and $L^{n-d-2}\subset  L^{n-d-1}$
(but the point need not lie on $L^{n-d-2}$).
 Set
$$
\chh_f(F):=\dsupp_S\bigl((\pi_f)_*\sigma_f^*F\bigr).
$$

\medskip
{\it Incidence version \ref{CC.defn.4.vers}.4.} Consider the diagram
$$
\begin{array}{lcl}
& \icorr_S\bigl((\mbox{point}), L^{n-d-1}, (\check{\p}^n)^{d+1}\bigr)&\\
\sigma_i\swarrow && \searrow\pi_i\\[1ex]
\p^n_S && \icorr_S\bigl(L^{n-d-1}, (\check{\p}^n)^{d+1}\bigr)
\end{array}
$$
where the incidence variety parametrizes   $(d+3)$-tuples
$\bigl((\mbox{point}), L^{n-d-1}, H_0,\dots, H_d\bigr)$ satisfying
$(\mbox{point})\in L^{n-d-1}\subset H_i$ for every $i$.  Set
$$
\chh_i(F):=\dsupp_S\bigl((\pi_i)_*\sigma_i^*F\bigr).
$$
\end{defn}

\begin{thm} \label{CC.flat.4.vers.thm}
Let $S$ be a scheme and $F$ a  generically flat family of pure, coherent sheaves  of dimension $d$ on $\p^n_S$. 
 The following are equivalent.
\begin{enumerate}
\item  $\chh_p(F)\subset (\check{\p}^n)^{d+1}_S$ is Cartier over $S$.
\item   $\chh_g(F)\subset \grass_S(n-d-1,\p^n)$ is Cartier over $S$.
\end{enumerate}
\noindent If $S$ is local with infinite residue field then these are also equivalent to
 \begin{enumerate}\setcounter{enumi}{2}
\item $\dsupp (\pi_*F)$ is Cartier over $S$ for every $\o_S$-projection  $\pi: \p^n_S\map \p^{d+1}_S$ (\ref{lin.projection.defn}) that is finite on $\supp F$.
\end{enumerate}
\end{thm}

Proof. 
The extreme cases $d=0$ and $d=n-1$ are somewhat exceptional, so we deal with them first. 

If $d=n-1$ then $\grass_S(n-d-1,\p^n_S)=\grass_S(0,\p^n_S)\cong \p^n_S$ and
the only projection is the identity. Furthermore 
$\chh_g(F)=\dsupp_S(F)$  by definition, so (2--4) are equivalent.
If these hold then $\chh_p(F)=\chh_p\bigl(\dsupp_S(F)\bigr)$
is also flat by (\ref{dsupp.car.nd.dl}). 
Conversely, for (1) $\Rightarrow$ (2) the argument in (\ref{CC.flat.4.vers.thm.pf.1}) works. 

If $d=0$ then $F$ is flat over $S$ and  (1--3) hold by (\ref{dsupp.rets.lem.3-1}). 

We may thus assume from now on that  $0<d<n-1$. These cases are discussed in
(\ref{CC.flat.4.vers.thm.pf.1}--\ref{CC.flat.4.vers.thm.pf.2}).

\begin{say}[Proof of (\ref{CC.flat.4.vers.thm}.1) $\Leftrightarrow$ (\ref{CC.flat.4.vers.thm}.2)]\label{CC.flat.4.vers.thm.pf.1}
To go between the product and the Grassmannian versions, the basic diagram is the following.
$$
\begin{array}{rcl}
& \icorr_S(L^{n-d-1}, \p^n_S)&\\
\swarrow && \searrow(\p^{d})^{d+1}{-}\mbox{bundle}\\[1ex]
(\check{\p}^n)^{d+1}_S && \grass_S(n{-}d{-}1, \p^n_S)
\end{array}
$$
The right hand side projection
$$
\pi_2: \icorr_S(L^{n-d-1}, \p^n_S)\to \grass_S(n{-}d{-}1, \p^n_S)
$$
is a $ (\p^{d})^{d+1}$-bundle. Therefore
$\chh_i(F)=\pi_2^*\chh_g(F)$. 
Thus $\chh_g(F)$ is Cartier over $S$ iff $\chh_i(F)$ is Cartier over $S$.
It remains to compare $\chh_i(F)$  and $\chh_p(F)$. 

The left hand side projection
$$
\pi_1: \icorr_S(L^{n-d-1}, \p^n_S)\to (\check{\p}^n)^{d+1}_S
$$
is birational. It is an isomorphism over
$(H_0,\dots, H_d)\in (\check{\p}^n)^{d+1}_S$ iff $\dim (H_0\cap\dots\cap H_d)=n{-}d{-}1$, the smallest possible. 
That is, when the rank of the matrix formed from the equations of the $H_i$ is $\leq d$.  Thus  $\pi_1^{-1}$ is an isomorphism outside a subset
of codimension $n+1-d$ in each fiber of $\pi_2$.

Therefore, if $\chh_i(F)$ is Cartier over $S$ then
$\chh_p(F)$ is Cartier over $S$, outside a subset
of codimension $n+1-d\geq 3$ on each fiber of $\pi_2$. 
Then
$\chh_p(F)$ is Cartier over $S$ everywhere by (\ref{pic.codim3.cor}).  

Conversely, if $\chh_p(F)$ is a relative Cartier divisor then so is
$\pi_1^*\chh_p(F)$, which is the union of $\chh_i(F)$ and of
the exceptional divisors. The latter are all relatively Cartier
  hence so is $\chh_i(F)$. 
\end{say}

\begin{say}[Proof of (\ref{CC.flat.4.vers.thm}.2) $\Leftrightarrow$ (\ref{CC.flat.4.vers.thm}.3--5)]\label{CC.flat.4.vers.thm.pf.2}
To go between the  Grassmannian version and the projection versions, the basic diagram is the following.
$$
\begin{array}{rcl}
& \flag_S(n{-}d{-}2, n{-}d{-}1, \p^n_S)&\\
\p^{n-d-1}{-}\mbox{bundle}\swarrow && \searrow\p^{d+1}{-}\mbox{bundle}\\[1ex]
\grass_S(n{-}d{-}1, \p^n) && \grass_S(n{-}d{-}2, \p^n)
\end{array}
$$
The left hand side projection
$$
\rho_1:\flag_S(n{-}d{-}2, n{-}d{-}1, \p^n_S)\to \grass_S(n{-}d{-}1, \p^n_S)
$$
is a $\p^{n-d-1}$-bundle and
$\chh_f(X)=\rho_1^* \chh_g(X)$.
Thus $\chh_g(F)$ is Cartier over $S$ iff $\chh_f(F)$ is Cartier over $S$.

The right hand side projection
$$
\rho_2:\flag_S(n{-}d{-}2, n{-}d{-}1, \p^n_S)\to \grass_S(n{-}d{-}2, \p^n_S)
$$
is a $\p^{d+1}$-bundle. 
Pick  $L\in \grass(n{-}d{-}2, \p^n_S)$. The fiber of $\rho_2$ over $[L]$
is the set of all $n{-}d{-}1$-planes that contain $L$; we can identify this with the target of the projection $\pi_L:\p^n\map L^\perp$.
So, if  $\chh_f(F)$ is Cartier over $S$ then 
$\dsupp\bigl((\pi_L)_*(F)\bigr)=\chh_f(F)|_{L^\perp}$
is also Cartier over $S$.

Conversely, assume that  $\dsupp\bigl((\pi_L)_*(F)\bigr)$
is  Cartier over $S$ for general $L$.

This gives an $L^\perp$ in $\grass_S(n{-}d{-}1, \p^n_S)$
where $\chh_g(F)$ is Cartier. Since $\dim_S L^\perp=d+1\geq 2$,
this implies that 
$\chh_g(F)$ is Cartier over $S$
outside a subset
of codimension $\geq 3$ on each fiber  of $\rho_2$. 
Then
$\chh_g(F)$ is Cartier over $S$ everywhere by (\ref{pic.codim3.cor}). 
\end{say}

\begin{cor} \label{CC.flat.4.vers.thm.cor.1}
Let $S$ be a   scheme  and  $F$ a  generically flat family of pure, coherent sheaves  of dimension $d$ on  $\p^n_S$.
Let $h:S'\to S$ be a  morphism. By base change we get
$g':X'\to S'$ and $F'=\vpure(h_X^*F)$ (\ref{vpure.defn}).  
\begin{enumerate}
\item If $F$ is C-flat, then so is $F'$.
\item If  $F$ is C-flat and $h$ is faithfully flat then $F$ is C-flat.
\end{enumerate}
\end{cor}

Proof.  We may assume that $S$ is local with infinite residue field.
Being C-flat is  exactly (\ref{CC.flat.4.vers.thm}.3) which is equivalent to
(\ref{CC.flat.4.vers.thm}.1).  $F\mapsto \chh_p(F)$ commutes with base change
by (\ref{div.supp.flat.bc.lem}) and, if  $h$ is faithfully flat, then a divisorial sheaf is Cartier iff its pull-back is (cf.\ \cite[4.22]{k-modbook}). \qed

\begin{defn}\label{loc.C.flat.defn}
Let $S$ be a local scheme  with infinite residue field and 
$F$ a generically flat family of pure, coherent  sheaves of dimension $d$ over $S$ (\ref{pure.over.defn}). 
 $F$ is {\it locally C-flat} over $S$ at $y\in Y:=\ssupp F$ iff $\dsupp (\pi_*F)$ is Cartier over $S$ at $\pi(y)$ for every  $\o_S$-projection  $\pi: \p^n_S\map \p^{d+1}_S$ that is finite on $Y$
for which  $\{y\}=\pi^{-1}(\pi(y))\cap Y$.
\end{defn}

\begin{lem}\label{loc.C.flat.C.flat.lem}
Let $S$ be a local scheme  with infinite residue field and 
$F$ a generically flat family of pure, coherent  sheaves of dimension $d$ on $\p^n_S$. Then $F$ is C-flat iff it is locally C-flat at every point.
\end{lem}

Proof. It is clear that C-flat implies locally C-flat.

Conversely, assume that $F$ is locally C-flat. Set
$Z_s:=\supp(F_s)\setminus \fcm (F)$ and pick  points
$\{y_i: i\in I\}$, one  in each irreducible component of $Z_s$.
If $\pi: \p^n_S\map \p^{d+1}_S$ is a general $\o_S$-projection, then 
  $\{y_i\}=\pi^{-1}(\pi(y_i))\cap Y$ for all $i\in I$.

Note that $\dsupp(\pi_*F)$ is a relative Cartier divisor along
$\p^{d+1}_s\setminus \pi(Z_s)$ by (\ref{dsupp.car.nd.dl}) and 
it is also relative Cartier at the points $\pi(y_i)$ for  $i\in I$ since
$F$ is locally C-flat.
Thus  $\dsupp(\pi_*F)$ is a relative Cartier divisor outside a
codimension $\geq 3$ subset of $\p^{d+1}_s$, hence a relative Cartier divisor everywhere by (\ref{pic.codim3.cor}). \qed

\begin{cor}\label{mumf.d.Cf.cor.2} Let $(s,S)$ be a local scheme and $X\subset \p^n_S$ a closed subscheme  that is flat over $S$ of pure relative dimension $d+1$. Let $D\subset X$ be a relative Mumford divisor. Let $x\in X_s$ be a  smooth  point. Then $\o_D$ is locally C-flat at $x$ iff
$D$ is a relative Cartier divisor at $x$.
\end{cor}

Proof.  We may assume that $S$ has infinite residue field.
 A general linear projection $\pi:X\to \p^{d+1}_S$ is \'etale at $x$.
Thus $D$ is a relative Cartier divisor at $x$ iff
$\pi(D)$ is a relative Cartier divisor at $\pi(x)$.
By (\ref{dsupp.car.nd.dl.4}.2)   
the latter holds iff 
 $\o_D$ is locally C-flat at $x$. \qed

\begin{cor}\label{mumf.d.Cf.cor.3}  Let $S$ be a scheme and 
$F$ a generically flat family of pure, coherent  sheaves of dimension $d$ over $S$ (\ref{pure.over.defn}). If $F$ is flat at $y\in Y:=\ssupp F$ then it is also locally C-flat at $y$.
\end{cor}

Proof. By \cite[10.8]{k-modbook}, $F_s$ is CM outside a subset $Z_s\subset Y_s$ of codimension
$\geq 2$. Let $W_s\subset Y_s$ be the set of points where $F$ is not flat.

Let $\pi:Y\to \p^{d+1}_S$ be a general linear projection.
By (\ref{dsupp.car.nd.dl})  $\dsupp(\pi_*F)$ is a relative Cartier divisor
outside  $\pi(Z_s\cup W_s)$. Moreover, we may assume that $\pi(y)\notin \pi(W_s)$. Thus, in a neighborhood of $\pi(y)$, 
$\dsupp(\pi_*F)$ is a relative Cartier divisor outside 
 the codimension $\geq 3$ subset $\pi(Z_s)$. 
Thus $\dsupp(\pi_*F)$ is a relative Cartier divisor at $y$ by 
(\ref{pic.codim3.cor}). \qed

\begin{lem}\label{bert.sec.lem.1-1} Let $S$ be a scheme and 
$F$ a generically flat family of pure, coherent  sheaves of dimension $d$ on $\p^n_S$.   Set $Y:=\ssupp F$ and let $\pi:Y\to \p^{d+1}_S$ be  a finite morphism.
Let $g_m:Y\into \p^N_S$ be an  embedding such that  $g_m^*\o_{\p^{N}_S}(1)\cong \pi^*\o_{\p^{d+1}_S}(m)$ for some $m\gg 1$.

If $(g_m)_*F$ is C-flat then  $\dsupp(\pi_*F)$ is a relative Cartier divisor.
\end{lem}

Proof.  We may assume that $S$ is local with infinite residue field.  Choosing $d+2$ general sections of $\o_{\p^{d+1}_S}(m)$ gives a
morphism $w_m: \p^{d+1}_S\to \p^{d+1}_S$ and there is a  linear projection
$\rho:\p^N_S\map \p^{d+1}_S$ such that
$w_m\circ \pi=\rho\circ g_m$. By assumption 
$\dsupp\bigl((\rho\circ g_m)_*F\bigr)$ is a relative Cartier divisor, hence so is 
$$
\dsupp\bigl((w_m\circ \pi)_*F\bigr)=\dsupp\bigl((w_m)_*\dsupp(\pi_*F)\bigr),
$$
where the equality follows from (\ref{dsupp.car.nd.dl.4}.3). 

Pick a point  $x\in \dsupp(\pi_*F)$. A general $w_m$
is \'etale at $x$ and $\{x\}=w_m^{-1}(w_m(x))\cap \dsupp(\pi_*F)$. 
Thus  $w_m: \dsupp(\pi_*F)\to \dsupp\bigl((w_m\circ \pi)_*F\bigr)$
is  \'etale at $x$. Thus 
$\dsupp(\pi_*F)$ is Cartier at $x$. 
 \qed

\begin{cor}\label{bert.sec.lem.1} Let $S$ be a scheme and 
$F$ a generically flat family of pure, coherent  sheaves of dimension $d$ on $\p^n_S$ (\ref{pure.over.defn}). 
Let $v_m:\p^n_S\into \p^N_S$ be the $m$th Veronese embedding.
If $(v_m)_*F$ is C-flat then so is $F$. \qed
\end{cor}

\subsection*{Bertini theorems for C-flatness}

\begin{lem} \label{bert.sec.lem.2} Let $(s,S)$ be a local scheme and 
$F$ a generically flat family of pure, coherent  sheaves of dimension $d\geq 1$ on $\p^n_S$ (\ref{pure.over.defn}). 
Set $Z_s:=\supp (F_s)\setminus \fcm(F)$. Let $H\subset \p^n_S$ be a hyperplane that does not contain any irreducible component of $Z_s$. 

If $F$ is C-flat then so is $F|_H$.
\end{lem}

Proof. We may assume that the residue field is infinite. 
Every projection $H\map \p^{d}_S$ is obtained as the restriction of a projection $\p^n_S\map \p^{d+1}_S$. The rest follows from (\ref{dsupp.rets.lem.4}). \qed

\begin{lem} \label{bert.sec.lem.3} Let $(s,S)$ be a local scheme and 
$F$ a generically flat family of pure, coherent  sheaves of dimension $d\geq 2$ on $\p^n_S$.
Then  $F$ is C-flat iff $F|_H$ is C-flat for an open, dense set of hyperplanes
$H$.
\end{lem}

Proof. One direction follows from (\ref{bert.sec.lem.2}). 
Conversely, if $F|_H$ is C-flat for an open, dense set of hyperplanes
$H$ then there is an open, dense set of projections
$\pi:\p^n_S\map \p^{d+1}_S$ such that for an open, dense set of hyperplanes
$L\subset \p^{d+1}_S$,  the restriction of $F$ to $\pi^{-1}(L)$ is C-flat. 
Thus  $\dsupp(\pi_* F)$ is a relative Cartier divisor 
in an open neighborhood of  such an $L$,
by (\ref{dsupp.rets.lem.lem4.c}).  Since $d\geq 2$, 
this implies that $\dsupp(\pi_* F)$ is a relative Cartier divisor  everywhere by (\ref{pic.codim3.cor}). 
Thus  $F$ is C-flat by (\ref{CC.flat.4.vers.thm}). \qed

\begin{lem} \label{bert.sec.lem.4} Let $(s,S)$ be a local scheme and 
$F$ a stably C-flat family of pure, coherent  sheaves of dimension $d\geq 1$ over $S$.
Set $Y:=\ssupp F$, $Z_s:=Y\setminus \fcm(F)$ and let  $D\subset Y$ be a relative Cartier divisor
 that does not contain any irreducible component of $Z_s$. 
Then $F|_D$ is also stably C-flat.
\end{lem}

Proof. We  may assume that the residue field is infinite. 
By (\ref{loc.C.flat.C.flat.lem}) it is sufficient to prove that  $F|_D$ is locally C-flat.
Pick a point $y\in D$ and let $H\supset D$ be a hypersurface section of $Y$ that does not contain any irreducible component of $Z_s$ and such that
$H$ equals $D$ in a neighborhood of $y$. 
After a Veronese embedding $H$ becomes a hyperplane section, and then 
 (\ref{bert.sec.lem.2}) implies that 
$F|_H$ is stably C-flat. Hence $F|_H$ is locally C-flat and so
$F|_D$ also locally C-flat at $y$.  \qed

\begin{defn}\label{stab.C.flat.H.defn}
Let $S$ be a local scheme and 
$F$ a generically flat family of pure, coherent  sheaves of dimension $d\geq 1$ over $S$. Set $Y:=\ssupp F$ and let  $L$ be a relatively ample line bundle on $Y$. 
We say that  $F$ is {\it stably C-flat} for $L$  over $S$ iff $\tau_*F$ is
C-flat for every embedding  $\tau:X\to \p^N_S$ 
such that $\tau^*\o_{\p^N_S}(1)\cong L^m$ for some $m\geq 1$.

By (\ref{bert.sec.lem.1-1}) this notion is unchanged if we replace $L$ by $L^r$ for some $r>0$.
\end{defn}

\begin{lem} \label{bert.sec.lem.5} Let $(s,S)$ be a local scheme and $F$ a generically flat family of pure, coherent  sheaves of dimension $d\geq 1$ over $S$.
 Let $L, M$ be relatively ample line bundles 
on $X$. Then $F$ is stably C-flat for $L$ iff it is stably C-flat for $M$.
\end{lem}

Proof. Assume that  $F$ is stably C-flat for $M$.
We may assume that $L$ is very ample. Repeatedly using 
 (\ref{bert.sec.lem.4}) we get that, for general $L_i\in |L_i|$, the  restriction of $F$ to the  complete intersection curve
$L_1\cap \cdots\cap L_{d-1}\cap Y$ is stably C-flat  for $M$.  
Thus the restriction of $F$ to 
$L_1\cap \cdots\cap L_{d-1}\cap Y$ is formally K-flat by (\ref{d=1.scf=kf=ff.prop}).
Using  (\ref{d=1.scf=kf=ff.prop}) in the other direction for $L$,
we get that the restriction of $F$ to 
$L_1\cap \cdots\cap L_{d-1}\cap Y$ is stably C-flat  for $L$. 
Now we can use (\ref{bert.sec.lem.3}) to conclude that
$F$ is stably C-flat for $L$. \qed

\begin{prop} \label{d=1.scf=kf=ff.prop}
Let $(s,S)$ be a local scheme and $F$ a generically flat family of pure, coherent  sheaves of dimension $1$ over $S$. Then $F$ is 
stably C-flat $\Leftrightarrow$  K-flat $\Leftrightarrow$ formally K-flat.
\end{prop} 

Proof. We already proved in (\ref{CK.flat.verss.imps}) that
formally K-flat $\Rightarrow$  K-flat $\Rightarrow$ stably C-flat.

Thus assume that $F$ is stably C-flat. Set $Y:=\ssupp F$ and 
pick a closed point $p\in Y$. We need to show that $F$ is formally K-flat at $p$. By (\ref{over.art.enough.ptop}) it is enough to prove this for Artin base schemes  and after a faithfully flat extension.
We may thus assume that $S=\spec A$ for a local Artin ring
 $(A,m,k)$ with $k$ infinite and   $p\in Y_s(k)$.

Let  $\hat\pi: \hat Y\to \hat\a^2_S=\spec A[[u,v]]$ be a finite morphism. 
After a linear coordinate change 
we may assume that the composite  $$
 \hat\rho:\hat Y\to \spec A[[u,v]]\to \spec A[[u]]\qtq{is also finite.}
$$
Thus $\hat\rho_*\hat F$ is a coherent sheaf on  $\spec A[[u]]$;   let 
 $\hat G$ denote its restriction  to the generic point.  Since $F$ is generically flat over $A$, $\hat G$ is  flat over $A$, hence we can write it as the sheafification of a  free  $A((u))$-module  $\oplus_j A((u))e_j$ of rank $n$ with basis $\{e_j\}$. Then multiplication by $v$ is given by a matrix  $M=\bigl(m_{ij}(u)\bigr)$ where
the $m_{ij}(u)\in A((u))$ are  Laurent series in $u$. By (\ref{dsupp.car.nd.dl})  $\dsupp(\hat \pi_*\hat F)$ is given by
$$
\bigl(\det (M-v{\mathbf 1}_n)=0\bigr)\subset \hat\a^2_S.
\eqno{(\ref{d=1.scf=kf=ff.prop}.1)}
$$
The finite map $\hat\pi$ is given by 2 power series $u, v$.  Fix some $m_0\in\n$, to be determined later. By (\ref{approx.proj.say}), for $m\gg m_0$ we can choose 
homogeneous polynomials  $g_1, g_2\in H^0\bigl(\p^n_A, \o_{\p^n_A}(m)\bigr)$
such that 
$$
\pi:Y\to \p^2_S\qtq{given by} (x_0^m{:}g_1{:}g_2)
\eqno{(\ref{d=1.scf=kf=ff.prop}.2)}
$$ is a finite morphism, 
$$
g_1/x_0^m\equiv u\mod \o_Y(-p)^{m_0}\qtq{and} g_2/x_0^m\equiv v\mod \o_Y(-p)^{m_0},
\eqno{(\ref{d=1.scf=kf=ff.prop}.3)}
$$
where $\o_Y(-p)$ is the ideal sheaf of $p\in Y$. 
Since the relative dimension is 1, $\o_Y(-p)^c\subset u\o_{Y_p}$ for some $c\in\n$, thus
$$
g_1/x_0^m\equiv u\mod u^{m_0/c}\o_{Y_p}\qtq{and} g_2/x_0^m\equiv v\mod u^{m_0/c}\o_{Y_p}
\eqno{(\ref{d=1.scf=kf=ff.prop}.4)}
$$
holds, whenever $c\mid m_0$.

We now compute $\dsupp (\pi_*F)$. This leads to  a  matrix $M'$ using
$u':=g_1/x_0^m, v':=g_2/x_0^m$, and  we see that
$$
M(u)\equiv M'(u')\mod u^{m_0/c}A[[u]].
\eqno{(\ref{d=1.scf=kf=ff.prop}.5)}
$$
Let $r$ be the maximal pole order
of the $m_{ij}(u)$.  Expanding the determinant we get terms whose  maximal pole order is $\leq nr$. Thus (\ref{d=1.scf=kf=ff.prop}.5) implies  that 
$$
\det (M-v{\mathbf 1}_n)\equiv \det (M'-v'{\mathbf 1}_n)\mod  u^{-nr+m_0/c}A[[u,v]].
\eqno{(\ref{d=1.scf=kf=ff.prop}.6)}
$$
If $m_0\gg 1$ (depending on $c,n,r$ and $M$) then, by (\ref{divs.sheaves.dim.2.cor.2}),  
$\det (M-v{\mathbf 1}_n)$ defines a Cartier divisor iff 
$\det (M'-v'{\mathbf 1}_n)$ does.
Since $F$ is stably C-flat,  $\det (M'-v'{\mathbf 1}_n)$ 
defines a Cartier divisor, hence so does $\det (M-v{\mathbf 1}_n)$.
Thus $F$ is  formally K-flat at $p$. \qed

\begin{say}[Approximation of formal projections]\label{approx.proj.say}
Let $(s, S)$ be a local scheme and $Y\subset \p^n_S$ a closed subset of pure relative dimension $d$. Let $p\in Y_s$ be a closed point with maximal ideal $m_p$ such that $x_0(p)\neq 0$. Let
$(\hat g_1{:}\cdots{:} \hat g_e):\hat Y_p\to  \hat{\a}^e_{\hat S}$ be a finite morphism.  Fix $m_0\in\n$.

Then, for every $m\gg m_0$ there are  $g_1,\dots, g_e\in H^0\bigl(\p^n_S, \o_{\p^n_S}(m)\bigr)$ such that
\begin{enumerate}
\item  $\pi: (x_0^m{:}g_1{:}\cdots{:} g_e):Y\to \p^e_S$ is a finite morphism, 
\item  $\pi^{-1}(\pi(p))\cap Y=\{p\}$ and
\item  $\hat g_i\equiv g_i/x_0^m \mod m_p^{m_0}$ for every $i$. \qed
\end{enumerate}
\end{say}

\begin{say}[Proof of (\ref{K-flat.=.sC-flat.thm})]\label{K-flat.=.sC-flat.thm.pf}
We already noted in (\ref{CK.flat.verss.imps}) that
K-flat $\Rightarrow$ stably C-flat.

To see the converse, assume that $F$ is stably C-flat and let 
$\pi:X\to \p^{d+1}_S$ be a finite projection. 
Set $L:=\pi^*\o_{\p^{d+1}_S}(1)$. By (\ref{bert.sec.lem.5})
$F$ is stably C-flat for $L$, hence 
$\dsupp(\pi_*F)$ is a relative Cartier divisor by
 (\ref{bert.sec.lem.1-1}). \qed
\end{say}

\section{Ideal of Chow equations}\label{chow.eq.ideal.sect}

The ideal of Chow equations was introduced in various forms in \cite{MR1174902, MR1368005, MR1714736}.

\begin{defn}[Ideal of Chow equations]\label{id.ch.eq.defn}
Let $S$ be a local scheme with infinite residue field and
 $Y\subset \p^n_S$  a generically flat family of pure subschemes of dimension $d$ (\ref{pure.over.defn}).   Let $\pi:\p^n_S\map \p^{d+1}_S$ be a projection that is finite on $Y$ and let $f(\pi, Y)$ be an equation of
$\dsupp_S(\pi_*\o_Y)$. The {\it ideal of Chow equations} of $Y$ is
$$
I^{\rm ch}(Y):=\bigl(\pi^*f(\pi, Y)\colon \mbox{ $\pi$ is finite on $Y$}\bigr)
\subset \o_S[x_0,\dots, x_n].
\eqno{(\ref{id.ch.eq.defn}.1)}
$$
We use 2 other versions of this concept.
If $I_Y$ is the ideal sheaf or the homogeneous ideal of $Y$, then sometimes we write $I^{\rm ch}(I_Y) $ instead of $I^{\rm ch}(Y)$.
If $S$ is reduced and $Z$ is well defined family of $d$-cycles,
then $\dsupp_S(\pi_*D)$ is defined in (\ref{id.ch.eq.defn}), hence
$I^{\rm ch}(Z)$ can be defined as in  (\ref{id.ch.eq.defn}.1). We see that if $Z(Y)$ denotes the  well defined family of $d$-cycles
associated to $Y$ then 
$$
I^{\rm ch}(Y)=I^{\rm ch}\bigl(Z(Y)\bigr).
\eqno{(\ref{id.ch.eq.defn}.2)}
$$
One can also define the obvious formal-local version of this notion.
In the examples that I know of, the global version is  compatible with the formal-local one, but this is not known in general. This is closely related to (\ref{CK.flat.verss.imps.4}). 

One difficulty is that  a typical space curve has non-linear projections that are not isomorphic to linear projections, not even locally analytically; see for example  (\ref{Cn.anal.proj.n34.exmp}). 
\end{defn}

We would like to compare the ideal of Chow equations with the ideal of $Y$.
A rather straightforward  result is the following.

\begin{prop}\cite[Sec.8]{MR1714736}
\label{k-null.sec.9.prop} Let $Y\subset \p^n_S \to S$ be C-flat. Then
$$
I^{\rm ch}\bigl(\pure (Y_s)\bigr)\subset I(Y_s). \qed
$$
\end{prop}

To get a more precise answer,  we need some definitions.

\begin{say}[Element-wise powers of ideals]\label{frob.power.defn}
Let $R$ be a ring, $I\subset R$ an ideal and $m\in \n$. 
Set
 $$
I^{[m]}:=(r^m: r\in I).
$$
These ideals have been studied mostly when $\chr k=p>0$ and $m=q$ is a power of $p$; one of the early occurrences is in   \cite{MR0432625}. In these cases  $I^{[q]}$ is called a {\it Frobenius power} of $I$. 
Other values of the exponent are also interesting, the following properties follow from (\ref{frob.power.defn}.6).
We assume for simplicity that $R$ is a $k$-algebra.
\begin{enumerate}
\item  If $I$ is principal then $I^{[m]}=I^m$.
\item  If $\chr k=0$ then  $I^{[m]}=I^m$.  
\item If $m<\chr k$ then $I^{[m]}=I^m$.
\item If $k$ is infinite then  $(r_1,\dots, r_n)^{[m]}=
\bigl((\sum c_i r_i)^m: c_i\in k\bigr)$.
\item If $k$ is infinite and $U\subset k^n$ is Zariski dense then  
$$(r_1,\dots, r_n)^{[m]}=
\bigl((\tsum c_i r_i)^m\colon (c_1,\dots, c_n)\in U\bigr).
$$
\end{enumerate}
Note that (2) is close to being optimal.  For example, if $I=(x,y)\subset k[x,y]$  and $\chr k=p\geq 3$ then 
$$
(x,y)^{[p+1]}=(x^{p+1}, x^py, xy^p, y^{p+1})\subsetneq (x,y)^{p+1}.
$$

{\it Claim \ref{frob.power.defn}.6.} Let $k$ be an infinite field. Then
$$
\bigl\langle (c_1x_1+\cdots+c_nx_n)^m: c_i\in k\bigr\rangle =
\bigl\langle x_1^{i_1}\cdots x_n^{i_n} : \tbinom{m}{i_1\dots i_n}\neq 0\bigr\rangle. 
$$
Here $\binom{m}{i_1\dots i_n}$ denotes the multinomial coefficient, that is, the
coefficient of $x_1^{i_1}\cdots x_n^{i_n} $ in  $(x_1+\cdots+x_n)^m$. 
\medskip

Proof. The containment $\subset $ is clear. If the  2 sides are not equal then the left hand side is contained in some hyperplane of the form $\sum \lambda_I x^I=0$, but this would give a nontrivial polynomial identity 
 $\sum \tbinom{m}{i_1\dots i_n} \lambda_I c^I=0$
for the $c_i$. \qed
\end{say}

\begin{prop} \label{chow.eq.bra.pow.thm} Let $k$ be an infinite field and 
 $Z_i\subset \p^n$ be  distinct, geometrically  irreducible and reduced $k$-cycles of dimension $d$. Then
$$
\pure \bigl( \o_{\p^n}/I^{\rm ch} \bigl(\tsum_i m_iZ_i\bigr)\bigr)= \pure \bigl(\o_{\p^n}/\cap_i I(Z_i)^{[m_i]}\bigr). 
\eqno{(\ref{chow.eq.bra.pow.thm}.1)}
$$
\end{prop}

Proof. Assume first that $d=0$, thus $Z_i=\{p_i\}$. Let $\ell_1$ be a linear form on $\p^n$ that vanishes  at $p_i$ but not at the other points.  Choosing a general $\ell_0$ we get a projection  $\pi:=(\ell_0, \ell_1):\p^n\map \p^1$.
At $p_i$, the resulting Chow equation  can be written as  $\ell_1^{m_i}(\mbox{unit})$.  By (\ref{frob.power.defn}.5) 
these generate $I(p_i)^{[m]}$. 

If $d\geq 0$ then it is clear that 
$$
 I^{\rm ch} \bigl(\tsum_i m_iZ_i\bigr)\subset \cap_i I(Z_i)^{[m_i]}. 
\eqno{(\ref{chow.eq.bra.pow.thm}.2)}
$$
In order  to check that they agree generically, choose a   projection   $\rho:\p^n\map  \p^d$ that is finite on $\cup_iZ_i$, and  a point $p\in \p^d$ such that  $\rho$ is \'etale over $p$.

Let $L\cong \p^{n-d}$ denote the closure of $\rho^{-1}(p)$ and set
$W_i:=L\cap Z_i$. It is clear that
$$
I(Z_i)^{[m_i]}|_L=I(W_i)^{[m_i]}\qtq{and}
 I^{\rm ch} \bigl(m_iW_i\bigr)\subset  I^{\rm ch} \bigl(m_iZ_i\bigr)|_L.
\eqno{(\ref{chow.eq.bra.pow.thm}.3)}
$$
(See (\ref{I.ch.d.n.rest.exmp}) for the second of these.)
Thus we get that
$$
\cap_i I(Z_i)^{[m_i]}|_L=\cap_i I(W_i)^{[m_i]}=
 I^{\rm ch} \bigl(\tsum_i m_iW_i\bigr)\subset  I^{\rm ch} \bigl(\tsum_i m_iZ_i\bigr)|_L.
\eqno{(\ref{chow.eq.bra.pow.thm}.4)}
$$
Therefore
$$
\cap_i I(Z_i)^{[m_i]}|_L= I^{\rm ch} \bigl(\tsum_i m_iZ_i\bigr)|_L.
\eqno{(\ref{chow.eq.bra.pow.thm}.5)}
$$
Note that $\cap_i I(Z_i)^{[m_i]}$ is flat over $p$, thus
(\ref{chow.eq.bra.pow.thm}.1) and (\ref{chow.eq.bra.pow.thm}.3) together
 imply that
$$
\cap_i I(Z_i)^{[m_i]}=  I^{\rm ch} \bigl(\tsum_i m_iZ_i\bigr)
\qtq{near  $L\cap \supp Z$;} 
$$
see for example  \cite[I.7.4.1]{rc-book}.\qed

\begin{defn}\label{chow.eq.bra.pow.thm.2} We call the  subscheme (\ref{chow.eq.bra.pow.thm}.1) the {\it Chow hull} of the cycle $Z=\sum m_iZ_i$, and denote it by $\chull(Z)$. 
\end{defn}

The following consequence is key to our study of Mumford divisors.

\begin{cor} \label{chow.eq.bra.pow.MDcor.0}
Let $k$ be an infinite field, $X\subset \p^N_k$  a reduced subscheme of pure dimension $n+1$ and $D\subset X$ a Mumford divisor of degree $d$. Then
$$
\pure \bigl(X\cap \chull(D)\bigr)= D. 
\eqno{(\ref{chow.eq.bra.pow.MDcor.0}.1)}
$$
\end{cor}

Proof. After a field extension we can write  $D=\sum m_iD_i$ where the $D_i$ are geometrically irreducible and reduced. Then
(\ref{chow.eq.bra.pow.thm}) says that
$$
\chull(D)=
\pure \bigl(\spec \o_{\p^N_k}/\cap_i I(D_i)^{[m_i]}\bigr).
$$
Let $g_i\in D_i$ be the generic point and $R_i$ its local ring in $\p^N_k$.
Let $J_i\subset R_i$ be the ideal defining $X$ and  $(J_i, h_i)$ 
the ideal defining $D_i$. The ideal defining the left hand side of
(\ref{chow.eq.bra.pow.MDcor.0}.1) is then
$\bigl(J_i+(J_i, h_i)^{[m_i]}\bigr)/J_i$. This is the same as
$(h_i)^{[m_i]}$, as an ideal in  $R_i/J_i$, which equals
$(h_i^{m_i})$ by (\ref{frob.power.defn}.1).  \qed

\medskip

Combining (\ref{k-null.sec.9.prop}) and (\ref{chow.eq.bra.pow.thm}), we get the  following partial answer to Question~\ref{bd.torrs.ques}.

\begin{cor}\cite[Sec.8]{MR1714736}\label{k.null.sec.8.cor}
 Let ${\mathbf C}\to S$ be a C-flat family whose pure fibers are geometrically reduced curves. Then
$$
\tors {\mathbf C}_s\subset I\bigl(\pure {\mathbf C}_s\bigr)/I^{\rm ch}\bigl(\pure {\mathbf C}_s\bigr). \qed
$$
\end{cor}

 It is not clear that  (\ref{k.null.sec.8.cor}) is optimal, 
but the next example shows that  
 is close to it  in some directions.

\begin{exmp} Consider the monomial curve $s\mapsto (s^a, s^{a+1})$ with equation
$x^{a+1}=y^a$.
The surface  $(s,t)\mapsto (s^a, s^{a+1}, st, t)$ defines a non-flat deformation of it in $\a^3_{xyz}\times \a^1_t$. 

If $t\neq 0$ then we get a complete intersection with equations
$x^{a+1}-y^a=xz-ty=0$.  Near $t=0$ we need more equations.  These are obtained by multiplying  $x^{a+1}=y^a$ with $(t/x)^m=(z/y)^m$ for $0<m\leq a$ to get
$x^{a+1-m}t^m=y^{a-m}z^m$. This shows that the central fiber is
$$
\spec k[x,y,z]/z\bigl((x)+(y,z)^{a-1}\bigr).
$$
Thus the torsion at the origin is isomorphic to $ k[y,z]/(y,z)^{a-1}$.

We compute in (\ref{codim.2.in.hypp.exmp}.7) that the ideal of Chow equations is
$$
\bigl(x^{a+1}-y^a, z^ix^{a+1-i},z^iy^{a-i}\colon  i=1,\dots, a\bigr). 
$$
Thus (\ref{k.null.sec.8.cor}) is close to being optimal, as far as the $y,z$ variables are concerned.

\end{exmp}

\begin{exmp} \cite[4.8]{MR1714736}\label{codim.2.in.hypp.exmp}
Using (\ref{proj.formulas.exmp}.1) we see that 
 the ideal of Chow equations of  the codimension 2 subvariety 
$\bigl(x_{n+1}=f(x_0,\dots, x_n)=0\bigr)\subset \p^{n+1}$
is generated by the forms
$$
f\bigl(x_0-a_0x_{n+1} : \cdots : 
x_n-a_nx_{n+1}\bigr) \qtq{for all} a_0,\dots,a_{n}.
\eqno{(\ref{codim.2.in.hypp.exmp}.1)}
$$
If the characteristic is 0 then Taylor's theorem  gives that
$$
f\bigl(x_0-a_0x_{n+1} : \cdots : 
x_n-a_nx_{n+1}\bigr)=
\sum_I  \frac{(-1)^I}{I!}a^I\frac{\partial^I f}{\partial x^I} x_{n+1}^{|I|},
\eqno{(\ref{codim.2.in.hypp.exmp}.2)}
$$
where $I=(i_0,\dots, i_n)\in \n^{n+1}$.
 The $a^{|I|}$ are linearly independent, hence we get that the  ideal of Chow equations is
$$
I^{\rm ch}\bigl(f(x_0,\dots, x_n), x_{n+1}\bigr)=\bigl(f, x_{n+1}D(f),  \dots, x_{n+1}^mD^m(f)\bigr),
\eqno{(\ref{codim.2.in.hypp.exmp}.3)}
$$
where we can stop at $m=\deg f$.  Here we use the usual notation for derivative ideals
$$
D(f):=\bigl(f, \tfrac{\partial f}{\partial x_0},\dots, \tfrac{\partial f}{\partial x_n}\bigr).
\eqno{(\ref{codim.2.in.hypp.exmp}.4)}
$$ 
(Note that if we work with the ideal $(f)$ and not just the polynomial $f$, then we must include
$f$ itself in its derivative ideal.)

If we want to work locally at the point $(x_1=\cdots=x_n=0)$, the we can set
$x_0=1$ to get the local version
$$
I^{\rm ch}\bigl(f(x_1,\dots, x_n), x_{n+1}\bigr)=\bigl(f, x_{n+1}D(f),  \dots, x_{n+1}^mD^m(f)\bigr),
\eqno{(\ref{codim.2.in.hypp.exmp}.7)}
$$
where we can now stop at $m=\mult f$. 
This also holds if $f$ is an analytic function, though this needs to be worked out using the more complicated   formulas (\ref{proj.formulas.exmp}.6) that for us become
$$
\pi: (x_1,\dots,x_{n+1}) \to \bigl(x_{1}-x_{n+1}\psi_1,\dots,
x_{n}-x_{n+1}\psi_n\bigr),
\eqno{(\ref{codim.2.in.hypp.exmp}.8)}
$$
where $\psi_i=\psi_i(x_0,\dots, x_{n+1})$ are analytic functions.
Expanding as in (\ref{codim.2.in.hypp.exmp}.2) we see that
$$
f\bigl(x_{1}-x_{n+1}\psi_1,\dots,
x_{n}-x_{n+1}\psi_n\bigr)\in
I^{\rm ch}\bigl(f(x_1,\dots, x_n), x_{n+1}\bigr).
\eqno{(\ref{codim.2.in.hypp.exmp}.9)}
$$
Thus we get the same ideal if we compute $I^{\rm ch} $ using
analytic or formal projections. 
\end{exmp}

The next example shows that taking the  ideal of Chow equations does not commute with general hyperplane sections; see \cite[5.1]{MR1714736} for a positive result.

\begin{exmp}\label{I.ch.d.n.rest.exmp}
 Start with  $F=f(x,y)+zg(x,y)$ where $\mult_0f=\mult_0g=m$.
So family is equimultiple along the $z$-axis.

We compute the ideal of Chow equations for   $\bigl(F(x,y,z), t\bigr)\subset k[x,y,z,t]$ and also for its restrictions
 $\bigl(F(x,y,c), t\bigr)\subset k[x,y,t]$.
The first one is
$$
\bigl(F, tD(F), t^2D^2(F), \dots, \bigr).
$$
Here $D(F)=(F, F_x, F_y, F_z)$, thus setting $z=c$ we get
$$
D\bigl(F(x,y,z)\bigr)|_{z=c}=\bigl(f +cg ,f_x +cg_x ,  f_y +cg_y , g \bigr).
\eqno{(\ref{I.ch.d.n.rest.exmp}.1)}
$$
By contrast 
$$
D\bigl(F(x,y,c)\bigr)=\bigl(f +cg ,f_x +cg_x ,  f_y +cg_y \bigr).
\eqno{(\ref{I.ch.d.n.rest.exmp}.2)}
$$
We have an extra term $g $ in (\ref{I.ch.d.n.rest.exmp}.1), which shows the following.

\medskip
{\it Claim \ref{I.ch.d.n.rest.exmp}.3.} If $g \notin \bigl( f +cg, f_x +cg_x ,  f_y +cg_y  \bigr)$ then taking $I^{\rm ch}$ does {\em not} commute with restriction to $(z=c)$. \qed

\medskip

To get a concrete example, take
$f(x,y)=x^4+y^4, g(x,y)=x^2y^2$. The ideal 
$ \bigl(f_x +cg_x ,  f_y +cg_y , f +cg \bigr)$
contains 5 obvious degree 4 elements, but Euler's equation
$$
x(f_x +cg_x)+y(f_y +cg_y)=4(f +cg)
$$ tells us that they are dependent.  An easy explicit computation shows that $x^2y^2$ is not in the ideal.

Note also that we  get the exact same computation if we restrict to
some other plane  $z=c+ax+by$.

\medskip
{\it Remark \ref{I.ch.d.n.rest.exmp}.4.} It is possible that this can be used to improve on (\ref{k.null.sec.8.cor}) and get a better 
 answer to (\ref{bd.torrs.ques}).
\end{exmp}

\section{Representability  Theorems}\label{C.K.rep.thms.sect}

\begin{defn}\label{c/k.flat.funct.repr.defn}
 Let $S$ be a   scheme  and  $F$ a  generically flat family of pure, coherent sheaves  of dimension $d$ on  $\p^n_S$. The {\it functor of C-flat pull-backs} is 
$$
CFlat_F(q:T\to S)=
\left\{
\begin{array}{l}
1\qtq{ if  $q^{[*]}F\to T$ is C-flat, and}\\
0  \qtq{otherwise,}
\end{array}
\right.
$$
where $q^{[*]}F :=\vpure(q^*F)$  is the divisorial pull-back as in (\ref{fam.divs.ch3}.5) or  (\ref{vpure.defn}).
We say that a monomorphism $j^{\rm cflat}_F:S^{\rm cflat}_F\to S$ {\it represents}  $ CFlat_F$, provided
$CFlat_F(q)=1$ iff $q$ factors as  $q:T\to S^{\rm cflat}_F\to S$.

If $Y\subset \p^n_S$ is a  generically flat family of pure subschemes  of dimension $d$ then we write $CFlat_Y $ instead of $CFlat_{\o_Y} $.

One defines  analogously the {\it functor of stably C-flat pull-backs}
 $ SCFlat_F$, and the {\it functor of K-flat pull-backs}
 $ KFlat_F$.
The monomorphisms representing them are denoted by  $j^{\rm scflat}_F:S^{\rm scflat}_F\to S$ and $ j^{\rm kflat}_F:S^{\rm kflat}_F\to S$.

Let $f:X\to S$ be a flat morphism and  $D$ a relative Mumford divisor on $X$.
The {\it functor of Cartier pull-backs} is 
$$
Cartier_D(q:T\to S)=
\left\{
\begin{array}{l}
1\qtq{ if  $q^{[*]}D\to T$ is Cartier, and}\\
0  \qtq{otherwise.}
\end{array}
\right.
$$
By \cite[4.35]{k-modbook},  the  functor of Cartier pull-backs
is represented by a monomorphism
$j^{\rm car}_D:S^{\rm car}_D\to S$.
(Note. Unfortunately \cite[4.35]{k-modbook} is about reduced schemes, but this assumption is not necessary in the proof. This will be fixed later.)
\end{defn}

An immediate consequence of (\ref{CC.flat.4.vers.thm}) is the following.

\begin{prop}\label{c.flat.funct.repr.thm}
 Let $S$ be a   scheme  and  $F$ a  generically flat family of pure, coherent sheaves  of dimension $d$ on  $\p^n_S$. Then the functor of C-flat pull-backs of $F$ is represented by a monomorphism
$j^{\rm cflat}_F:S^{\rm cflat}_F\to S$.
\end{prop}

Proof.  By (\ref{CC.flat.4.vers.thm}), 
$j^{\rm cflat}_F:S^{\rm cflat}_F\to S$ is the same as  $j^{\rm car}_{\chh_p(F)}:S^{\rm car}_{\chh_p(F)}\to S$,  with the Chow hypersurface  $\chh_p(F)$ as defined in (\ref{CC.defn.4.vers}.2). \qed

\begin{cor}\label{sc.flat.funct.repr.thm}
 Let $S$ be a   scheme  and  $F$ a  generically flat family of pure, coherent sheaves  of dimension $d$ on  $\p^n_S$. Then the functor of stably C-flat pull-backs of $F$ is represented by a monomorphism
$j^{\rm scflat}_F:S^{\rm scflat}_F\to S$.
\end{cor}

Proof. Let $F_m$ denote the push-forward of $F$ by the $m$th Veronese embedding $v_m$. Then $ S^{\rm scflat}_F$ should be the fiber product
$$
S^{\rm scflat}_{F_1}\times S^{\rm scflat}_{F_2}\times \cdots.
\eqno{(\ref{sc.flat.funct.repr.thm}.1)}
$$
However, this need not make sense if the product is truly infinite.
To understand this, we make a slight twist and as the first step we replace
$S$ by $S^{\rm cflat}_F$ and $F$ by $F^*:=(j^{\rm cflat}_F)^{[*]}F$. Then we consider
$$
j^{\rm cflat}_{F^*_m}:\bigl(S^{\rm cflat}_F\bigr)^{\rm cflat}_{F^*_m}\to S^{\rm cflat}_F.
\eqno{(\ref{sc.flat.funct.repr.thm}.2)}
$$
Note that  $j^{\rm cflat}_{F^*_m} $ is an isomorphism on the underlying reduced subschemes by (\ref{over.red.K.autom.ptop}).
Thus $j^{\rm cflat}_{F^*_m} $ is a proper monomorphism, hence a closed immersion \cite[3.47.1]{k-modbook}.

Now the infinite fiber product (\ref{sc.flat.funct.repr.thm}.1) is an infinite intersection of closed subschemes, which always exists. Thus we get that 
$$
S^{\rm scflat}_F=\cap_m \bigl(S^{\rm cflat}_F\bigr)^{\rm cflat}_{F^*_m}
\subset S^{\rm cflat}_F,
\eqno{(\ref{sc.flat.funct.repr.thm}.3)}
$$
and $j^{\rm scflat}_F$ is the restriction of $j^{\rm cflat}_F$ to it. 
 \qed

\medskip
A combination of (\ref{K-flat.=.sC-flat.thm}) and (\ref{sc.flat.funct.repr.thm}) gives the following.

\begin{cor}\label{k.flat.funct.repr.thm}
 Let $S$ be a   scheme  and  $F$ a  generically flat family of pure, coherent sheaves  of dimension $d$ on  $\p^n_S$. Then the functor of K-flat pull-backs of $F$ is represented by a monomorphism
$j^{\rm kflat}_F:S^{\rm kflat}_F\to S$. \qed
\end{cor}

\begin{say}[Construction of the Chow variety I]\label{Cvar.exits.thm}
In order to construct $\chow_{n,d}(\p^N_S)$, the Chow variety of degree $d$ cycles of dimension $n$ in $\p^N_S$, 
 we start with the diagram 
(\ref{CC.defn.4.vers}.2) 
$$
\p^N_S \stackrel{\sigma}{\longleftarrow}\icorr_S   \stackrel{\pi}{\longrightarrow}(\check{\p}^N)^{n+1}_S,
\eqno{(\ref{Cvar.exits.thm}.1)}
$$
where  the incidence variety $\icorr_S:=\icorr_S\bigl((\mbox{point}), (\check{\p}^N)^{n+1}\bigr)$
 parametrizes   $(n+2)$-tuples
$\bigl(p, H_0,\dots, H_n\bigr)$ satisfying
$p\in H_i$ for every $i$, where $p$ is a point in $\p^N_S$ and we view the $H_i$ either as hyperplanes in  $\p^N_S$ or points in $\check{\p}^N_S$.

Let ${\mathbf P}_{N,n,d}=|\o_{(\check{\p}^N)^{n+1}}(d,\dots, d)|$
be the linear system of hypersurfaces of multidegree $(d,\dots, d) $
in $(\check{\p}^N)^{n+1}$   with universal  hypersurface
$$
\mathbf{CH}_{N,n,d}\subset (\check{\p}^N)^{n+1}\times {\mathbf P}_{N,n,d}.
$$
Thus (\ref{Cvar.exits.thm}.1) extends to
$$
\begin{array}{rcl}
& \icorr_S \times_S {\mathbf P}_{N,n,d}&\\
\sigma_{N,n,d}\swarrow && \searrow\pi_{N,n,d}\\[1ex]
\p^N_S \times_S {\mathbf P}_{N,n,d}&& (\check{\p}^N)^{n+1}_S\times_S {\mathbf P}_{N,n,d}
\end{array}
\eqno{(\ref{Cvar.exits.thm}.2)}
$$
Consider now the restriction of the left hand projection
$$
\bar\sigma_{N,n,d}: \bigl(\icorr_S \times_S {\mathbf P}_{N,n,d}\bigr)\cap 
\pi_{N,n,d}^{-1}\mathbf{CH}_{N,n,d}
\ \to  \ \p^N_S \times_S {\mathbf P}_{N,n,d}.
\eqno{(\ref{Cvar.exits.thm}.3)}
$$
Note that the preimage of a pair  $$\bigl(p=(\mbox{point}), CH=(\mbox{Cayley-Chow-type hypersurface})\bigr)$$ consists of all $(d+1)$-tuples
$(H_0,\dots, H_n)$ such that
$p\in H_i$ for every $i$ and 
$(H_0,\dots, H_n)\in CH$.

In particular, if $Z$ is an $n$-cycle of  degree $d$ on $\p^N_S$ and 
$\chh_p(Z)$ is its Cayley-Chow  hypersurface then $\bar\sigma_{N,n,d} $ 
is a $(\check{\p}^{N-1}_S)^{n+1}$-bundle over $\supp Z$. The key observation is that  this property alone  is enough to construct the Chow variety.

By the Flattening Decomposition Theorem \cite[Lec.8]{mumf66}, there is a unique, largest, locally closed subscheme
$$
{\mathbf W}_{N,n,d}\into \p^N_S \times_S {\mathbf P}_{N,n,d}
\eqno{(\ref{Cvar.exits.thm}.4)}
$$
over which  $\bar\sigma_{N,n,d}$ is a   $(\check{\p}^{N-1})^{n+1}$-bundle.
Note that if $(p, CH)\in  {\mathbf W}_{N,n,d}$ then
 $\bar\sigma_{N,n,d}^{-1}(p, CH)$ is the product of $n+1$ copies of the dual hyperplane  
$H(\check{p})\subset \check{\p}^N$, that is, the set of all hyperplanes that contain $p$.

The set-theoretic behavior of the projection
$$\rho_{N,n,d}:{\mathbf W}_{N,n,d}\to  {\mathbf P}_{N,n,d}\eqno{(\ref{Cvar.exits.thm}.5)}
$$ is rather clear; see
 \cite[I.3.24.4]{rc-book}. The fiber dimension of $\rho_{N,n,d}$  is $\leq n$, and if $CH\in {\mathbf P}_{N,n,d}$ is irreducible then 
$\dim \rho_{N,n,d}^{-1}(CH)=n$  iff $CH$ is the Cayley-Chow hypersurface of
$Z:=\red \rho_{N,n,d}^{-1}(CH)$. 

In the reducible case one has to be more careful with the multiplicities; 
this was completed in \cite{MR1513117}. A scheme-theoretic version of this is done in (\ref{chow.eq.bra.pow.thm}). 
The end result is that  there is a closed subset
$\chow'_{n,d}(\p^N_S)\into {\mathbf P}_{N,n,d}$ that parametrizes 
Cayley-Chow hypersurfaces of $n$-cycles of degree $d$.

One would like this closed subset to be  $\chow_{n,d}(\p^N_S)$. Unfortunately,  its scheme structure may change if we apply a Veronese embedding of $\p^N_S$; see \cite{MR0096668} or \cite[I.4.2]{rc-book}.
For this reason \cite{rc-book} defines the Chow variety  $\chow_{n,d}(\p^N_S)$ as the seminormalization of  $\chow'_{n,d}(\p^N_S)$.
\end{say}

In order to get a scheme-theoretic version of $\chow_{n,d}(\p^N_S)$, one needs to understand the scheme-theoretic fibers of $\rho_{N,n,d}$.  We consider this next.

\begin{say}[Construction of the Chow variety II] \label{cc.inverse.general.say}  
Let $S$ be a local scheme with residue field $k$. 
Let $Y\subset \p^N_S$ be a generically flat family of subschemes of dimension $n$,  degree $d$ and  $\chh_p(Y)\subset  (\check{\p}^N)^{n+1}_S$ its Cayley-Chow hypersurface.

Choose coordinates  $(x_0{:}\dots{:}x_N)$ on $\p^N_S$ and dual coordinates
 $(x_{0j}^\vee{:}\dots{:}x_{Nj}^\vee)$ on the $j$th copy of $\check{\p}^N_S $
for $j=0,\dots, n$.  Let $F_Y(x_{ij}^\vee)=0$ be an equation of $\chh_p(Y)$.

 For notational simplicity we compute in the affine cart
$\a^N_S=\p^N_S\setminus (x_0=0)$.  Assume that none of the irreducible components of $Y_k$ is contained in the coordinate hyperplane  $(x_0=0)$. 

For $(p_1,\dots, p_N)\in  \a^N_S$, the hyperplanes $H$ in the   $j$th copy of $\check{\p}^N_S $ that pass through $(p_1,\dots, p_N)$  are all written in  the form
$$
\bigl(-\tsum_{i=1}^N p_i x_{ij}^\vee: x_{1j}^\vee:\dots: x_{Nj}^\vee\bigr).
$$
Thus $(p_1,\dots, p_N)\in Y\cap \a^N_S$
iff
$F$ identically vanishes after the substitutions
$$
x_{0j}\mapsto -\tsum_{i=1}^N p_i x_{ij}^\vee\qtq{for}  j=0,\dots, n.
\eqno{(\ref{cc.inverse.general.say}.1)}
$$
If $M(x_{ij}^\vee)$ are all the monomials in the $x_{ij}^\vee $
for $1\leq i\leq N, 0\leq j\leq n$ then, after the substitutions (\ref{cc.inverse.general.say}.1),  we can write $F_Y$ as
$$
\tsum_M  f_{Y,M}(p_1,\dots, p_N)M(x_{ij}^\vee).
\eqno{(\ref{cc.inverse.general.say}.2)}
$$
Since the monomials $M(x_{ij}^\vee)$ are linearly independent, this vanishes for all $x_{ij}^\vee$ iff  $f_{Y,M}(p_1,\dots, p_N)=0 $ for every $M$. Equivalently:

\medskip
{\it Claim \ref{cc.inverse.general.say}.3.}  The subscheme 
$$
{\mathbf W}_{N,n,d}\cap (\a^N_S \times_S {\mathbf P}_{N,n,d})\subset \a^N_S \times_S {\mathbf P}_{N,n,d}
$$
defined in (\ref{Cvar.exits.thm}.4)
 is given by the equations  
$f_{Y,M}(x_1,\dots, x_N)=0$ for all monomials $M$. \qed
\medskip

If we fix  $x_{ij}^\vee=c_{ij}$ then these give a linear projection
$\pi_{\mathbf c}:\a^N_S \to \a^{n+1}_S$, and  the corresponding Chow equation is
$$
\tsum_M  f_{Y,M}(x_1,\dots, x_N)M(c_{ij})=0.
\eqno{(\ref{cc.inverse.general.say}.4)}
$$
Thus we get the following.

\medskip
{\it Claim \ref{cc.inverse.general.say}.5.}  If the residue field of $S$  is infinite then $I^{\rm ch}(Y)|_{\a^N_S}$ is generated by the Chow equations  
of the linear projections  $\pi_{\mathbf c}:\a^N_S \to \a^{n+1}_S$. \qed
\medskip

Note that a priori we would need to use the more general projections
(\ref{proj.formulas.exmp}.4).
\end{say}

Using (\ref{cc.inverse.general.say}.3), 
we can reformulate (\ref{chow.eq.bra.pow.thm}) as follows.

\begin{cor} \label{chow.eq.bra.pow.cor}
Let $k$ be an infinite field and $Z=\sum m_iZ_i$  a sum of   distinct, geometrically irreducible and reduced $n$-cycles of dimension $d$ in $\p^N_S$. Then the purified fiber of 
$\rho_{N,n,d}:{\mathbf W}_{N,n,d}\to  {\mathbf P}_{N,n,d}$ over 
$\chh_p(Z)$ is 
$$
\pure\bigl(\rho_{N,n,d}^{-1}[\chh_pZ]\bigr)=\chull(Z).\qed
$$
\end{cor}

Combining  (\ref{chow.eq.bra.pow.cor}) and (\ref{chow.eq.bra.pow.MDcor.0}) gives the following.

\begin{cor} \label{chow.eq.bra.pow.MDcor}
Let $k$ be an infinite field, $X\subset \p^N_k$  a reduced subscheme of pure dimension $n+1$ and $D\subset X$ a Mumford divisor of degree $d$. Then
$$
\pure \bigl(X\cap  \rho_{N,n,d}^{-1} [\chh_pD]\bigr)= D. \qed
$$
\end{cor}

\begin{say}[Construction of the Chow variety III]\label{Cvar.exits.thm.2}
Let us now return to  (\ref{Cvar.exits.thm}.4) 
$$\rho_{N,n,d}:{\mathbf W}_{N,n,d}\to  {\mathbf P}_{N,n,d}.
\eqno{(\ref{Cvar.exits.thm.2}.1)}
$$
By (\ref{chow.eq.bra.pow.cor}) the (purified) fiber of $\rho_{N,n,d}$ over
$\chh_pZ$ is the Chow hull of $Z$, which is usually much larger than $Z$.
Therefore  $\rho_{N,n,d}$ is not even generically flat over $[\chh_pZ]\in {\mathbf P}_{N,n,d}$
if $Z$ is not geometrically reduced.

However, if $Z$ is  geometrically reduced then
its Chow hull equals $\o_Z$ by (\ref{chow.eq.bra.pow.thm}).  Thus there is a chance that $\rho_{N,n,d}$ is  generically flat over $[\chh_pZ]$.
 
Using  (\ref{generic.flattening.thm}) we get a generically flattening decomposition
$$
j^{\rm g-flat}_{{\mathbf W}_{N,n,d}}\colon \bigl({\mathbf P}_{N,n,d}\bigr)^{\rm g-flat}_{{\mathbf W}_{N,n,d}}\to {\mathbf P}_{N,n,d},
\eqno{(\ref{Cvar.exits.thm.2}.2)}
$$
and then (\ref{k.flat.funct.repr.thm}) gives
$$
j^{\rm kflat}_{{\mathbf W}_{N,n,d}}\colon \bigl({\mathbf P}_{N,n,d}\bigr)^{\rm kflat}_{{\mathbf W}_{N,n,d}}\to {\mathbf P}_{N,n,d},
\eqno{(\ref{Cvar.exits.thm.2}.3)}
$$
which represents K-flatness. Finally,  being  generically geometrically reduced is an open condition, hence we get 
 $$
\chow_{n,d}^{\rm g-red}(\p^N_S)\subset \bigl({\mathbf P}_{N,n,d}\bigr)^{\rm kflat}_{{\mathbf W}_{N,n,d}},
\eqno{(\ref{Cvar.exits.thm.2}.4)}
$$
which represents the functor of K-flat families of
geometrically reduced $n$-cycles of degree $d$ in $\p^N_S$.
\medskip

{\it Remark \ref{Cvar.exits.thm.2}.5.} It would be more in the spirit of classical Chow theory to use C-flatness in (\ref{Cvar.exits.thm.2}.3) instead of K-flatness. However, when one  defines  $\chow_{n,d}^{\rm g-red}(X/S)$ for some projective scheme $X\to S$, we would like the  result to be independent of the embedding $X\into \p^N_S$. Thus K-flatness is more natural. 
\end{say}

I do not know whether it is possible to push through the above approach for the whole Chow variety. Fortunately, the method works with minor changes for C-flat families of Mumford divisors. This then completes the proof of (\ref{KDIV.exits.thm}).

\begin{thm} \label{CMDIV.exits.thm} Let $X\subset \p^N_S$ be a closed subscheme that is flat over $S$  with $S_2$ fibers of pure dimension $n$.   Then the functor $KDiv_d(X/S)$ of K-flat,  relative Mumford divisors of degree $d$
 is representable by  a separated $S$-scheme of finite type  $\kdiv_d(X/S)$.
\end{thm}

Proof.  Let $D_s\subset X_s$ be a Mumford divisor of degree $d$. We can also view it as an
$(n-1)$-cycle  of degree $d$ in $\p^N_S$. 
We proceed as in  (\ref{Cvar.exits.thm}.1--4)  to get
$$
{\mathbf W}_{N,n-1,d}\into \p^N_S \times_S {\mathbf P}_{N,n-1,d}.
\eqno{(\ref{CMDIV.exits.thm}.1)}
$$
We are only interested in cycles that lie on $X$, hence we focus on the restriction of the coordinate projection 
$$
\bar\rho_{N,n-1,d}: {\mathbf W}_{N,n-1,d}\cap \bigl(X\times_S {\mathbf P}_{N,n-1,d}\bigr)\to {\mathbf P}_{N,n-1,d}.
\eqno{(\ref{CMDIV.exits.thm}.2)}
$$
Let $D_s\subset X_s$ be a Mumford divisor of degree $d$.
By (\ref{chow.eq.bra.pow.MDcor}),
$$
D_s=\pure\bigl(X_s\cap \chull(Z(D_s))\bigr)=\pure\bigl(\bar\rho^{-1}_{N,n-1,d}[\chh_pD_s]\bigr).
\eqno{(\ref{CMDIV.exits.thm}.3)}
$$
Let $F$ be the structure sheaf of ${\mathbf W}_{N,n-1,d}\cap \bigl(X\times_S {\mathbf P}_{N,n-1,d}\bigr) $.
By (\ref{generic.flattening.thm}) 
there is a locally closed partial decomposition 
$$
j^{\rm g-flat}_F:{\mathbf P}_{N,n-1,d,F}^{\rm g-flat}\to {\mathbf P}_{N,n-1,d},
\eqno{(\ref{CMDIV.exits.thm}.4)}
$$ such that 
$F_W$ is generically flat in dimension $n-1$ over an $S$-scheme $W$  iff $W\to S$ factors through   $j^{\rm g-flat}_F$. 

Thus ${\mathbf P}_{N,n-1,d,F}^{\rm g-flat}$ paramerizes generically flat families of divisorial subschemes of $X$ of degree $d$.  Applying  (\ref{k.flat.funct.repr.thm}) now gives
$$
j^{\rm kflat}_F:\bigl({\mathbf P}_{N,n-1,d,F}^{\rm g-flat}\bigr)^{\rm kflat}\to
{\mathbf P}_{N,n-1,d,F}^{\rm g-flat},
\eqno{(\ref{CMDIV.exits.thm}.5)}
$$
which  paramerizes K-flat Mumford divisors.  \qed
\medskip

We have used the following is a variant of the Flattening Decomposition Theorem of 
\cite[Lec.8]{mumf66}. 

\begin{thm} \label{generic.flattening.thm}
Let $f:X\to S$ be a projective morphisms and $F$ a coherent sheaf on $X$. Let  $n$ be the maximal fiber dimension of $\supp F\to S$.
There is a locally closed decomposition $j^{\rm g-flat}_F:S^{\rm g-flat}_F\to S$ such that 
$F_W$ is generically flat over $W$ in dimension $n$ iff $W\to S$ factors through   $j^{\rm g-flat}_F$.
\end{thm}

Proof. We may replace $X$ by $\supp F$. The question is local on $S$. By Noether normalization we  may  assume that there is a finite morphism $\pi:X\to  \p^n_S$. Note that $F_W$ is generically flat over $W$ in dimension $n$ iff the same holds for $\pi_*F_W$. We may thus assume that  $X=\p^n_S$. 

Applying \cite[Lec.8]{mumf66}  to the identity $X\to X$ and $F$, we get a   decomposition
$\amalg_i X_i\to X=\p^n_S$ where every 
$F|_{X_i}$ is flat, hence locally free of rank $i$.

Let  $Z\subset \p^n_S$ be a closed subscheme. 
 Applying \cite[Lec.8]{mumf66}  to the projection  $\p^n_S\to S$ and $\o_Z$, we see that there is a unique
 largest subscheme $S(Z)\subset S$ such that
$S(Z)\times_S\p^n_S\subset Z$.  For a locally closed subscheme  $Z\subset \p^n_S$
set $S(Z)=S(\bar Z)\setminus S(\bar Z\setminus Z)$, where $\bar Z$ denotes the scheme-theoretic closure of $Z\subset \p^n_S$. 

Note that $S(Z)$ is the largest subscheme  $T\subset S$ with the following property:
\begin{enumerate}
\item There is an open subscheme ${\mathbf P}^0_T\subset \p^n_T$
that  contains the
generic point of $\p^n_t$ for every $t\in T$ and such that
${\mathbf P}^0_T\subset Z$.
\end{enumerate}

We claim that $S^{\rm g-flat}_F=\amalg_i S(X_i)$.

First, $F|_{X_i}$ is locally free of rank $i$, so the
restriction of $F$ to $S(X_i)\times_S\p^n_S$ is locally free, hence flat, 
at every generic point of every fiber.

Conversely, let $W$ be a connected scheme and $q:W\to S$ a morphism  
such that $F_W$ is generically flat over $W$ in dimension $n$.
Since $F_w$ is generically free for every $w\in W$, this implies that
$F_W$ is locally free at the generic point of every fiber. 
Let ${\mathbf P}^0_W\subset \p^n_W$ be the open set where $F_W$ is locally free.
By assumption the closure of ${\mathbf P}^0_W$ equals $\p^n_W$. 

Since  ${\mathbf P}^0_W$ contains the generic point of every fiber $\p^n_w$, it is connected. Thus $F$ has constant rank, say $i$, on ${\mathbf P}^0_W$. Therefore, the restriction of $q$ to  ${\mathbf P}^0_W$ lifts to  $\tilde q: {\mathbf P}^0_W\to X_i$, which in turn extends to the closures
$\bar q: \p^n_W\to \bar X_i$. Thus  $\bar q$ gives $q_W: W\to S(X_i)$ in view of (1). \qed

\section{Hypersurface singularities}\label{hyp.surf.K.flat.sec}

In this section we give a detailed description of K-flat deformations of hypersurface singularities over $k[\epsilon]$.

\begin{say}[Non-flat deformations]\label{non.fl.def.gen.say}
Let $X\subset \a^n$ be a reduced subscheme of pure dimension $d$. 
We aim to describe non-flat deformations of $X$ that are flat outside a subset $W\subset X$. 

Choose equations $g_1,\dots, g_{n-d}$ such that
$$
(g_1=\cdots=g_{n-d}=0)=X\cup X', 
$$
where $Z:=X\cap X'$ has dimension $<d$.
Let $h$ be an equation of $X'\cup W$ that does not vanish on any irreducible component of $X$. Thus $X$ is a complete intersection in $\a^n\setminus (h=0)$
with equation $g_1=\cdots=g_{n-d}=0$.  Its flat deformations over an Artin ring $(A,m,k)$ are then given by
$$
g_i({\mathbf x})=\Psi_i({\mathbf x})\qtq{where} \Psi_i\in m[x_1,\dots, x_n, h^{-1}].
\eqno{(\ref{non.fl.def.gen.say}.1)}
$$
Note that we can freely change the $\Psi_i$ by any element of the ideal
$\bigl(g_1-\Psi_1,\dots, g_{n-d}-\Psi_{n-d}\bigr)$. 
We get especially simple normal forms  if $A=k[\epsilon]$, that is, we look at first order deformations. In this case the equations can be written as
$$
g_i({\mathbf x})=\Phi_i({\mathbf x})\epsilon \qtq{where} \Phi_i\in k[x_1,\dots, x_n, h^{-1}].
\eqno{(\ref{non.fl.def.gen.say}.2)}
$$
Now we can freely change the $\Phi_i$ by any element of the ideal
$(g_1,\dots, g_{n-d})$. Thus first order generically flat deformations can be given in the form
$$
g_i=\phi_i\epsilon\qtq{where} \phi_i\in H^0(X, \o_X)[h^{-1}].
\eqno{(\ref{non.fl.def.gen.say}.3)}
$$
Set $X^\circ:=X\setminus (Z\cup W)$. By varying $h$ we see that in fact
$$
g_i=\phi_i\epsilon\qtq{where} \phi_i\in H^0\bigl(X^\circ, \o_{X^\circ}\bigr).
\eqno{(\ref{non.fl.def.gen.say}.4)}
$$
This shows that the choice of $h$ is largely irrelevant.

If the deformation is flat then the equations defining $X$ lift, that is,
$\phi_i\in H^0\bigl(X, \o_{X}\bigr)$.  In some simple cases, for example if
$X$ is a complete intersection, this is equivalent to flatness. 
In the examples that we compute the most important information is carried by the polar parts
$$
\bar \phi_i\in H^0\bigl(X^\circ, \o_{X^\circ}\bigr)/H^0\bigl(X, \o_{X}\bigr).
\eqno{(\ref{non.fl.def.gen.say}.5)}
$$
\end{say}

We study first order non-flat  deformations of hypersurface singularities.
Plane curves turn out to be the most interesting ones.

\begin{say}\label{hypsurf.say.1}
Consider a  hypersurface  singularity $X:=(f=0)\subset \a^n_{\mathbf x}$
and a generically flat deformation of it 
$$
{\mathbf X}\subset \a^{n+1}_{\mathbf x,z}[\epsilon]\to \spec k[\epsilon].
\eqno{(\ref{hypsurf.say.1}.1)}
$$
Aiming to work inductively, we assume that the deformation is flat outside the origin.  Choose coordinates such that the $x_i$ do not divide $f$.

As in (\ref{non.fl.def.gen.say}.3)  any such deformation can be given as 
$$
f({\mathbf x})=\psi({\mathbf x})\epsilon\qtq{and}
z=\phi({\mathbf x})\epsilon,
\eqno{(\ref{hypsurf.say.1}.2)}
$$
where 
$\psi, \phi  \in  H^0(X, \o_X)[x_n^{-1}]$. 
Note that the choice of $x_n$ is not intrinsic, so in fact
$$
\psi, \phi  \in  \cap_i H^0(X, \o_X)[x_i^{-1}].
\eqno{(\ref{hypsurf.say.1}.3)}
$$
If $n\geq 3$ then $\cap_i\o_X[x_i^{-1}]=\o_X$ and we get the following
special case of \cite[10.68]{k-modbook}. 

\medskip

{\it Claim \ref{hypsurf.say.1}.4.}  Let ${\mathbf X}\subset \a^{n+1}$ be a 
first order deformation over $k[\epsilon]$ of a
 hypersurface  singularity $X:=(f=0)\subset \a^n$ 
that is flat outside the origin. If $n\geq 3$ then ${\mathbf X}$ is flat over $k[\epsilon]$. \qed
\medskip

If $n=2$ then we have a curve singularity $C=\bigl(f(x,y)=0\bigr)\subset \a^2$.
Set $C^\circ:=C\setminus\{(0,0)\}$. Then the deformation is given as
$$
f(x,y)=\psi\epsilon\qtq{and}
z=\phi\epsilon, 
\eqno{(\ref{hypsurf.say.1}.5)}
$$
where the relevant information about $\phi,\psi$ is carried by the polar parts
$$
\bar\psi,\bar \phi\in H^0(C^\circ, \o_{C^\circ})/ H^0(C, \o_{C}).
\eqno{(\ref{hypsurf.say.1}.6)}
$$
\end{say}

\begin{defn} We say that a (flat resp.\ generically flat) deformation over $k[\epsilon]$ {\it globalizes} if it is induced from a (flat resp.\ generically flat) deformation over $k[[t]]$ by base change.
\end{defn} 

\begin{thm} \label{planecurve.K.defns.thm} Consider a generically flat deformation ${\mathbf C}$ of the reduced plane curve  singularity $C:=(f=0)\subset \a^2_{xy}$ given in (\ref{hypsurf.say.1}.5--6).

\begin{enumerate}
\item If ${\mathbf C}$  is C-flat then $\psi \in H^0(C, \o_{C})$. 
\item If $\psi \in H^0(C, \o_{C})$ then
the deformation is
\begin{enumerate}
\item flat iff $\phi\in  H^0(C, \o_{C})$,
\item globalizes iff $\phi\in H^0(\bar C,\o_{\bar C})$ where $\bar C\to C$ is the normalization, and
\item C-flat iff  $f_x\phi, f_y\phi\in H^0(C, \o_{C})$. 
\end{enumerate}
\end{enumerate}
\end{thm}

Proof. If $\psi, \phi\in H^0(C, \o_{C})$ then we can rewrite
(\ref{hypsurf.say.1}.5) as
$$
f(x,y)-\tilde\psi(x,y)\epsilon=z-\tilde\phi(x,y)\epsilon=0
$$
 where $\tilde\psi, \tilde\phi$ are regular; this is a flat deformation and
the converse is clear.

If $\phi\in H^0(\bar C,\o_{\bar C})$  then it is integral over $k[x,y]_{(x,y)}$, so it satisfies an equation
$$
\phi^m+\tsum_{j=0}^{m-1} r_j(x, y)\phi^j=0.
$$ 
Consider now the surface $S$ given by the equations
$$
\begin{array}{l}
f(x,y)-\psi(x,y)s = 0\\
x^{a_1}z-\Phi_1(x,y)s= y^{c_1}z-\Phi_2(x,y)s=0\qtq{and}\\
z^m+\tsum_{j=0}^{m-1} r_j(x, y)z^js^{m-j}=0,
\end{array}
\eqno{(\ref{planecurve.K.defns.thm}.3)}
$$
where $a_1, c_1$ are chosen so that  $\Phi_1(x,y)=x^{a_1}\phi, \Phi_2(x,y)=y^{c_1}\phi$ are regular. The equations in line 2 of (\ref{planecurve.K.defns.thm}.3) guarantee that the projection 
$S\to  \bigl( f(x,y) =\psi(x,y)s\bigr)\subset \a^3_{xys}$ is birational and the
last equation shows  that  it is finite.  We also see that
$z-\phi s\in (s^2)$, thus $S$ is a globalization of ${\mathbf C}$. 
The converse assertion in (b) follows from (\ref{1par.def.S2.hull.say}).

As for (c), we  write down the image of  the projection 
$$
(x,y,z)\mapsto (\bar x, \bar y)=\bigl(x-\alpha(x,y,z) z, y-\gamma(x,y,z) z\bigr).
$$ 
where $\alpha, \gamma$ are constants for linear projections and 
power series that are nonzero at the origin in general. 

Since $z^2=0$ we get that
$$
f(x,y)=f(\bar x,\bar y)+\alpha(x,y,z) f_x(x,y) z +\gamma(x,y,z) f_y(x,y) z
$$
holds in $\o_{\mathbf C}$. 
Similarly, for any polynomial $F(x,y)$ we get that
$F(\bar x,\bar y)\equiv F(x,y)\mod \epsilon\o_{\mathbf C}$, hence
$F(\bar x,\bar y)z= F(x,y)z$ in $\o_{\mathbf C}$ since $z\epsilon=0$. 
(Here we use that $\epsilon^2=0$, we get other terms in (\ref{monomial.exmp.3}.3) otherwise.)

Thus the equation of the projection is
$$
f(\bar x,\bar y)+\bigl(\psi(\bar x,\bar y)+\alpha(\bar x,\bar y,0) f_x(\bar x,\bar y)\phi +\gamma(\bar x,\bar y,0) f_y(\bar x,\bar y)\phi\bigr)\cdot \epsilon =0. 
\eqno{(\ref{planecurve.K.defns.thm}.4)}
$$
By (\ref{divs.sheaves.dim.2.cor}) this defines a relative Cartier divisor 
for every $\alpha, \gamma$
iff
$\psi, f_x\phi, f_y\phi\in \o_C$.  (In particular, linear projections and formal projections give the same restrictions.) \qed

\begin{rem}  Note that  $\Omega^1_C$ is generated by $dx|_C, dy|_C$, while  $\omega_C$ is generated by  $f_y^{-1}dx=-f_x^{-1}dy$. 

Since $C$ is reduced, $\Omega^1_C$  and $\omega_C$  are naturally isomorphic over the smooth locus $C^\circ$. This gives a natural inclusion
 $\Hom(\Omega^1_C, \omega_C)\into \o_{C^\circ}$. Then ((\ref{planecurve.K.defns.thm}.2.c) says that
 C-flat deformations as in (\ref{hypsurf.say.1}.5) are parametrized by
$\Hom(\Omega^1_C, \omega_C)$. 
We describe this space  for monomial curves next. 
\end{rem}

\begin{exmp}[Monomial curves] \label{monomial.exmp.3}  We can be more explicit if $C$ is the irreducible monomial curve 
 $C:=(x^a=y^c)\subset \a^2$ where $(a,c)=1$.
It can be parametrized as  $t\mapsto (t^c, t^a)$. Thus
$\o_C=k[t^c, t^a]$. Let $E_C=\n a+\n c\subset \n$ denote the semigroup of exponents.  Then the condition  (\ref{planecurve.K.defns.thm}.2.c)  becomes
$$
t^{ac-c}\phi(t), t^{ac-a}\phi(t)\in k[t^a, t^c].
$$
This needs to be checked one monomial at a time. For $\phi=t^{-m}$ we get
$ac-c-m\in E_C$ and  $ac-a-m\in E_C$.
By (\ref{monomial.exmp.3}.1) these are  equivalent to  $ac-a-c-m\in E_C$.
The largest value of $m$ satisfying this condition gives the deformation
$$
\bigl(x^a-y^c=z-t^{-ac+a+c}\epsilon=0)\qtq{over} k[\epsilon].
$$
Note also that for $0\leq m\leq ac-a-c$, we have that
$ac-a-c-m\in E_C$ iff $m\notin E_C$.  Thus we see that the space of C-flat deformations that do not globalize has dimension $\frac12 (a-1)(c-1)$.
(This is an integer since one of $a,c$ must be odd.)

The following is left as an exercise.

\medskip
{\it Lemma \ref{monomial.exmp.3}.1.}  For  $(a,c)=1$ set  $E=\n a+\n c\subset \n$. Then
\begin{enumerate}
\item[(a)] If 
$0\leq m\leq \min\{ac-a, ac-c\}$ then 
$ac-a-m, ac-c-m\in E$ iff $ac-a-c-m\in E$.
\item[(b)] If  $0\leq m\leq ac-a-c$ then 
$ac-a-c-m\in E$ iff $m\notin E$. \qed
\end{enumerate}
\end{exmp}

\begin{say}[$S_2$ hull of a deformation]\label{1par.def.S2.hull.say} Let $T$ be the spectrum of a DVR with maximal ideal $(t)$ and residue field $k$.
Let  $g:X\to T$ be a flat morphism of pure relative dimension $d$ and  $Z:=\supp \tors (X_k)$.
Let $j:X\setminus Z\into X$ the  natural injection and
set $\bar X:=\spec_X j_*\o_{X\setminus Z}$. 
If $X$ is excellent then $\pi:\bar X\to X$ is finite and $\bar X$ is $S_2$.

By composition we get $\bar g:\bar X\to T$. Note that
$\pi_k:\bar X_kto X_k$ is an isomorphism over $X_k\setminus Z $ and
$\bar X_k$ is $S_1$. Thus if  $\pure (X_k)$ is reduced then 
$\bar X_k$ is dominated by the normalization $X_k^{\rm nor}\to X_k$.

Note that  $ t^n\o_X$ usually has some embedded primes contained in $Z$.
The intersection of its height 1 primary ideals (also called the $n$th symbolic power of $t\o_X$) is
$$
(t\o_X)^{(n)}=\o_X\cap t^n\o_{\bar X}=\ker\bigl[\o_X\to \pure\bigl(\o_X/t^n\o_X)\bigr]
\eqno{(\ref{1par.def.S2.hull.say}.1)}
$$
Multiplication by $t$ gives injections
$$
\pure (\o_{X_k})=\o_X/(t\o_X)^{(1)}\stackrel{t}{\into}  (t\o_X)^{(1)}/(t\o_X)^{(2)} \stackrel{t}{\into}  \cdots\into 
\eqno{(\ref{1par.def.S2.hull.say}.2)}
$$
Note that
$$
(t\o_X)^{(n)}/(t\o_X)^{(n+1)}\into  t^n\o_{\bar X}/t^{n+1}\o_{\bar X}\cong \o_{\bar X_k},
\eqno{(\ref{1par.def.S2.hull.say}.3)}
$$
thus the sequence (\ref{1par.def.S2.hull.say}.2) eventually stabilizes.
We can thus view the quotients
$$
(t\o_X)^{(n+1)}/t(t\o_X)^{(n)}
$$
as graded pieces of two  filtrations, one of $\tors (X_k)$ and one of
$\o_{\bar X_k}/\o_{X_k}$.

To formalize this, let us  write $M\preceq N$ to mean that there are  filtrations  $0=M_0\subset \cdots\subset M_m=M$,  $0=N_0\subset \cdots\subset N_n=N$ and an injection  $\sigma:[1,\dots,m]\into [1,\dots,n]$ such that, 
$M_{i}/M_{i-1}\cong N_{\sigma(i)}/N_{\sigma(i)-1}$ for every $i=1,\dots, m$.
If $M, N$ are artinian modules over a local ring then this holds iff $\len M\leq \len N$.

We have thus proved the following.
\end{say}

\begin{cor} \label{glob.cdef.partl.norm.c1} Using the notation of (\ref{1par.def.S2.hull.say}), assume that $\pure (X_k)$ is reduced with  normalization  $X_k^{\rm nor}\to X_k$. Then
$$
\tors \o_{X_k}\preceq \o_{X_k^{\rm nor}}/\o_{X_k}.
$$
In particular, if $\dim X_k=1$  then
$$
\len\bigl(\tors \o_{X_k}\bigr)\leq \len\bigl( \o_{X_k^{\rm nor}}/\o_{X_k}\bigr).\qed
$$
\end{cor}

\section{Seminormal curves}\label{s.n.curve.K.flat.sec}

Over an algebraically closed field $k$, every seminormal curve singularity is formally isomorphic to
$C_n\subset \a^n_{\mathbf x}$,   formed by the union of the $n$ coordinate axes. 
Equivalently,
$$
C_n=\spec k[x_1,\dots, x_n]/(x_ix_j: i\neq j).
$$
In this section we study deformations of $C_n$ over $k[\epsilon]$
that are flat outside the origin.

A normal form is worked out in (\ref{Cn.genflat.1st.say}.4), which shows that the space of these deformations is infinite dimensional.
Then we describe the flat deformations (\ref{Cn.flat.1st.exmp}) and their relationship to  smoothings  (\ref{Cn.smoothings.exmp}).

We compute C-flat and K-flat deformations in (\ref{Cn.CC.K.flat.thm}); these turn out to be quite close to flat deformations.

The ideal of Chow equations of $C_n$ is computed in (\ref{I.ch.axes.prop.1}).
 For $n=3$ these are close to C-flat deformations, but the difference between the two classes increases rapidly with $n$.

\begin{say}[Generically flat deformations of $C_n$]\label{Cn.genflat.1st.say}
Let  ${\mathbf C}_n\subset \a^m_{\mathbf x}[\epsilon]$ be  a generically flat deformation  of $C_n\subset \a^m_{\mathbf x}$ over $k[\epsilon]$.  

If ${\mathbf C}_n$ is flat over $k[\epsilon]$ then we can assume that $n=m$, but a priori we only know that $n\leq m$.

Following (\ref{non.fl.def.gen.say}), we can describe ${\mathbf C}_n$ as follows.

Along the $x_j$-axis and away from the origin, the deformation is flat
and the $x_j$-axis is a complete intersection. Thus, in the $(x_j\neq 0)$ open set,  ${\mathbf C}_n$ can be  given as
$$
x_i=\Phi_{ij}(x_1,\dots, x_n)\epsilon \qtq{where} i\neq j\qtq{and} \Phi_{ij}\in k[x_1,\dots, x_n, x_j^{-1}].
\eqno{(\ref{Cn.genflat.1st.say}.1)}
$$
Note that $(x_1,\dots, \hat x_j, \dots, x_n, \epsilon)^2$ is identically zero on ${\mathbf C}_n\cap (x_n\neq 0)$, so the terms in this ideal can be ignored. Thus along the $x_j$-axis we can change (\ref{Cn.genflat.1st.say}.1) to  the simpler form
$$
x_i=\phi_{ij}(x_j)\epsilon \qtq{where} i\neq j\qtq{and} \phi_{ij}\in k[x_j, x_j^{-1}].
\eqno{(\ref{Cn.genflat.1st.say}.2)}
$$
There is one more simplification that we can make.  Write  
$$
\phi_{ij}=\phi^p_{ij}+\gamma_{ij}\qtq{where}
\phi^p_{ij}\in k[x_j^{-1}], \gamma_{ij}\in (x_j)\subset k[x_j],
$$
and set $x'_i= x_i-\tsum_{j\neq i} \gamma_{ij}(x_j)$. Then we get the description
$$
x'_i=\phi^p_{ij}(x'_j)\epsilon \qtq{where} i\neq j\qtq{and} \phi^p_{ij}\in k[{x'_j}^{-1}].
\eqno{(\ref{Cn.genflat.1st.say}.3)}
$$
For most of our computations the latter coordinate change is not very important. 
Thus we write our deformations as
$$
{\mathbf C}_n: 
\bigl\{x_i=\phi_{ij}(x_j)\epsilon \quad\mbox{along the $x_j$-axis}\bigr\},
\eqno{(\ref{Cn.genflat.1st.say}.4)}
$$
where $\phi_{ij}(x_j)\in k[x_j, x_j^{-1}]$,  but we keep in mind that we can choose  $\phi_{ij}(x_j)\in k[x_j^{-1}]$  if it is convenient. 
In order to deal wit the cases when $m>n$, we make the following

\medskip
{\it Convention \ref{Cn.genflat.1st.say}.5.} We set $\phi_{ij}\equiv 0$ for $j>n$.
\medskip

Writing ${\mathbf C}_n $ as in (\ref{Cn.genflat.1st.say}.4) is almost  unique;
see  (\ref{Cn.flat.1st.rem}) for one more coordinate change that leads to a unique  normal form.

We get the same result (\ref{Cn.genflat.1st.say}.4) if we work with the analytic or formal local scheme of $C_n$: we still end up with $\phi_{ij}(x_j)\in k[x_j^{-1}]$. 

\end{say}

\begin{prop} \label{Cn.flat.1st.exmp}
For $n\geq 3$ the generically  flat deformation ${\mathbf C}_n\subset \a^n_{\mathbf x}[\epsilon]$  as in (\ref{Cn.genflat.1st.say}.4)
is flat  iff the $\phi_{ij}(x_j) $ have no poles.
(See (\ref{Cn.gflat.sp.exmp}.5) for the $n=2$ case.)
\end{prop}

 Proof.   If the $\phi_{ij}(x_j) $ are regular then 
$$
x_ix_j+\bigl(x_j\phi_{ij}(x_j)+x_i\phi_{ji}(x_i)\bigr)\epsilon=0
\eqno{(\ref{Cn.flat.1st.exmp}.1)}
$$
is an equation for ${\mathbf C}_n$. Thus every equation of $C_n$ lifts to an 
equation of ${\mathbf C}_n$, hence  ${\mathbf C}_n$ is flat over $k[\epsilon]$ by (\ref{flat.over.art.lem}).

Conversely, if the deformation is flat then  the equations defining $C_n$ lift, so we have
a set of defining equations for ${\mathbf C}_n$ of the form
$$
x_ix_j=\Psi_{ij}(x_1,\dots, x_n)\epsilon.
\eqno{(\ref{Cn.flat.1st.exmp}.2)}
$$
As in (\ref{Cn.genflat.1st.say}.2), this simplifies to
$$
x_ix_j=\psi_{ij}(x_j)\epsilon \qtq{along the $x_j$-axis.}
$$
Note that $x_ix_j$ vanishes along the other $n-2$ axes, so we must have
$\psi_{ij}(0)=0$.  (Here we use that $n\geq 3$.)  Thus $\phi_{ij}:=x_j^{-1}\psi_{ij}$ is regular as needed. \qed

\begin{rem}  Choosing $r\leq n$ of the coordinate axes we get an embedding
$\tau_r:C_r\into C_n$ and any generically flat deformation ${\mathbf C}_n$ of $C_n$
induces a generically flat deformation ${\mathbf C}_r:=\tau_r^*{\mathbf C}_n$ of $C_r$.

From (\ref{Cn.flat.1st.exmp}) we conclude that ${\mathbf C}_n$ is flat iff 
$\tau_3^*{\mathbf C}_n$ is flat for every $\tau_3:C_3\into C_n$. Neither direction of this claim  seems to follow from general principles.
For example, if $\tau_2^*{\mathbf C}_n$ is flat for every $\tau_2:C_2\into C_n$
then ${\mathbf C}_n$ need not be flat; see (\ref{Cn.gflat.sp.exmp}.5)
and (\ref{Cn.CC.K.flat.thm}.5). 
\end{rem}

\begin{rem} \label{Cn.flat.1st.rem}
Putting (\ref{Cn.genflat.1st.say}.3) and (\ref{Cn.flat.1st.exmp}) together we get that
flat deformations can be given as
$$
{\mathbf C}_n: 
\bigl\{x_i=e_{ij}\epsilon \quad\mbox{along the $x_j$-axis, where $e_{ij}\in k$}\bigr\}.
\eqno{(\ref{Cn.flat.1st.rem}.1)}
$$
The constants $e_{ij}$  are not yet unique, translations 
$$
x_i\mapsto x_i- a_i\epsilon\qtq{change} e_{ij}\mapsto e_{ij}-a_j.
 \eqno{(\ref{Cn.flat.1st.rem}.2)}
$$
So we get a  first order deformation  space of dimension   $n(n-1)-n=n(n-2)$.

We can also think of $\o_{{\mathbf C}_n}$ as a subring of
$\oplus_j  k[X_j, \epsilon_j]$ given by
$$
x_i\mapsto \bigl(e_{i1}\epsilon_1, \dots, e_{i,i-1}\epsilon_{i-1}, X_i, e_{i,i+1}\epsilon_{i+1},\dots, e_{in}\epsilon_n\bigr).
$$

Strangely, (\ref{Cn.flat.1st.rem}.1) says that every flat first order deformation of $C_n$ is obtained by translating the axes independently of each other.
These deformations all globalize in the obvious way, but the globalization is not a flat deformation of $C_n$  unless the translated axes all pass through the same point. 
If this point is   $(a_1\epsilon,\dots, a_n\epsilon)$ then
$e_{ij}=a_j$ and applying (\ref{Cn.flat.1st.exmp}.3) we get the trivial deformation. 

If $n=2$ then the  universal deformation is
$x_1x_2+\epsilon=0$. One may ask why this deformation does not lift to a deformation of $C_3$: smooth 2 of the axes to a hyperbola and just move the 3rd axis along. If we use $x_1x_2+t=0$, then the $x_3$-axis should move to
the line  $(x_1-\sqrt{t}=x_2-\sqrt{t}=0)$. This gives the flat deformation
given by equations
$$
x_1x_2+t=x_3(x_1-\sqrt{t})=x_3(x_2-\sqrt{t})=0.
$$
Of course this only makes sense if $t$ is a square. Thus setting
$\epsilon=\sqrt{t}\mod t$ the $t=\epsilon^2\mod t$ term becomes 0 and we get
$$
x_1x_2=x_3x_1-x_3\epsilon=x_3x_2-x_3\epsilon=0,
$$
which is of the form given in (\ref{Cn.flat.1st.exmp}.1).
\end{rem}

\begin{exmp}[Smoothing $C_n$]\label{Cn.smoothings.exmp}
Rational normal curves  $R_n\subset \p^n$  have a moduli space of dimension $(n+1)(n+1)-1-3=n^2+2n-3$. 
The $C_n \subset \p^n$  have a moduli space of dimension
$n+n(n-1)=n^2$. Thus the smoothings of $C_n$ have a  
moduli space of dimension $n^2+2n-3-n^2=2n-3$. 
We can construct these smoothings explicitly as follows.

Fix distinct $p_1,\dots, p_n\in k$ and consider the map
$$
(t,z)\mapsto  \bigl(\tfrac{t}{z-p_1}, \dots, \tfrac{t}{z-p_n}\bigr).
$$
Eliminating $z$ gives the equations
$$
(p_i-p_j)x_ix_j+(x_i-x_j)t=0 \colon 1\leq i\neq j\leq n
\eqno{(\ref{Cn.smoothings.exmp}.1)}
$$
for the closure of the image, which
 is an affine cone over a degree $n$ rational normal curve
$R_{n}\subset \p^n_{t,\mathbf x}$.  So far this is an $(n-1)$-dimensional space of smoothings.

Applying the  torus action  $x_i\mapsto \lambda_i^{-1} x_i$, we get new smoothings given by the equations
$$
(p_i-p_j)x_ix_j+(\lambda_jx_i-\lambda_ix_j)t=0 \colon 1\leq i\neq j\leq n.
\eqno{\ref{Cn.smoothings.exmp}.2}
$$
Writing it in the form (\ref{Cn.genflat.1st.say}.4) we get
$$
x_i=\tfrac{\lambda_i}{p_i-p_j}\epsilon \quad\mbox{along the $x_j$-axis.}
\eqno{\ref{Cn.smoothings.exmp}.3}
$$
This looks like a $2n$-dimensional family,  but  $\aut(\p^1)$ acts on it,
reducing the dimension to the expected $2n-3$.
The action is clear for $z\mapsto \alpha z+\beta$, but $z\mapsto z^{-1}$ also works out using (\ref{Cn.flat.1st.rem}.2) since
$$
\tfrac{\lambda_i}{p_i^{-1}-p_j^{-1}}=
\tfrac{-\lambda_ip_i^2}{p_i-p_j}+\lambda_ip_i.
$$

{\it Claim \ref{Cn.smoothings.exmp}.4.}  For distinct $p_i\in k$  and $\lambda_j\in k^*$, the vectors
 $$
\bigl(\tfrac{\lambda_j}{p_i-p_j}\colon i\neq j\bigr)
\qtq{span} \langle e_{ij}\rangle \cong k^{\binom{n}{2}}.
$$
So the flat infinitesimal deformations determined in (\ref{Cn.flat.1st.rem}.1) form the Zariski  tangent space of the smoothings.
\medskip

Proof. Assume that there is a linear relation
$$
\tsum_{ij} m_{ij}\tfrac{\lambda_j}{p_i-p_j}=0.
$$
If we let $p_i\to p_j$ but keep  the others fixed, we get that $m_{ij}=0$. \qed

\medskip

{\it Remark \ref{Cn.smoothings.exmp}.5.}  If $n=3$ then $2n-3=n(n-2)$ and
the Hilbert scheme of degree 3 reduced space curves with  $p_a=0$ is smooth, see \cite{MR796901}. 
\end{exmp}

\begin{exmp}[Simple poles]\label{Cn.gflat.sp.exmp} 
Among non-flat deformations, the simplest ones are given by 
$\phi_{ij}(x_j)=c_{ij}x_j^{-1}+e_{ij}$.  Then we have
$$
q_{ij}:=x_ix_j-(e_{ij}x_i+e_{ji}x_j)\epsilon =
\left\{
\begin{array}{cl}
c_{ij}\epsilon &\mbox{along the $x_j$-axis,}\\
c_{ji}\epsilon &\mbox{along the $x_i$-axis,}\\
0 &\mbox{along the other axes.}
\end{array}
\right.
\eqno{(\ref{Cn.gflat.sp.exmp}.1)}
$$
Thus we see that
$\sum_{ij}\gamma_{ij}q_{ij}$ vanishes on   ${\mathbf C}_n$ iff
$$
\tsum_i\gamma_{ij}c_{ij}  \qtq{is independent of $j$.}
\eqno{(\ref{Cn.gflat.sp.exmp}.2)}
$$
These impose  $n-1$ linear conditions on the $\gamma_{ij}$, which are in general independent.  Thus we get the following.
\medskip

{\it Claim \ref{Cn.gflat.sp.exmp}.3.} For general $c_{ij}$, 
the torsion in the central fiber has length $n-1$. \qed
\medskip

In special cases the torsion can be smaller, but
 if the $c_{ij}$ are not identically 0, then we get at least 1 nontrivial condition. This is 
in accordance with (\ref{Cn.flat.1st.exmp}). 

The $n=2$ case is exceptional and is worth discussing separately.
We get that
$$
q_{12}:=x_1x_2-(e_{12}x_1+e_{21}x_2)\epsilon =
\left\{
\begin{array}{cl}
c_{12}\epsilon &\mbox{along the $x_2$-axis,}\\
c_{21}\epsilon &\mbox{along the $x_1$-axis.}\\
\end{array}
\right.
\eqno{(\ref{Cn.gflat.sp.exmp}.4)}
$$
This gives the following.
\medskip

{\it Claim \ref{Cn.gflat.sp.exmp}.5.} For $n=2$ the deformation
as in (\ref{Cn.genflat.1st.say}.4) is flat iff
  $\phi_{12}, \phi_{21}$ have only simple poles and with  the same residue. \qed
\end{exmp}

The main result is the following.

\begin{thm} \label{Cn.CC.K.flat.thm} For a first order deformation of $C_n\subset \a^m$ specified as in (\ref{Cn.genflat.1st.say}.4) by 
$$
{\mathbf C}_n: 
\bigl\{x_i=\phi_{ij}(x_j)\epsilon \quad\mbox{along the $x_j$-axis}\bigr\}
\eqno{(\ref{Cn.CC.K.flat.thm}.1)}
$$
the following are equivalent.
\begin{enumerate}\setcounter{enumi}{1}
\item ${\mathbf C}_n $ is C-flat.
\item ${\mathbf C}_n $ is K-flat.
\item The $\phi_{ij}$ have only simple poles and 
$\phi_{ij},\phi_{ji}$ have the same residue. 
\item ${\mathbf C}_n $ induces a flat deformation on any pair of lines
$C_2\into C_n$.
\end{enumerate}
\end{thm}

Proof. The proof consist of 2 parts. First we show in 
(\ref{Cn.CC.flat.ord.1.exmp}) that (2) and (4) are equivalent by
explicitly computing the equations of linear projections.

 We see in 
(\ref{Cn.K.flat.nl.exmp}) that if the $\phi_{ij}$ have only simple poles then 
there is only 1 term of the equation of a non-linear projection that could have a pole, and this term is the same for the linearization of the projection. Hence it vanishes iff it vanishes for linear projections. 
This  shows that (4)  $\Rightarrow$ (3).

Finally (4)  $\Leftrightarrow$ (5)  follows from (\ref{Cn.gflat.sp.exmp}.5). \qed

\begin{rem} If $j>n$ then $\phi_{ij}\equiv 0$ by (\ref{Cn.genflat.1st.say}.5),
so $\phi_{ji}$ is regular by   (\ref{Cn.CC.K.flat.thm}.4).
Thus
$$
x_j-\tsum_{\ell=1}^n \phi_{j\ell}x_\ell \epsilon
$$
is identically on ${\mathbf C}_n $. We can thus eliminate
the $x_j$ for $j>n$.  Hence we see that allowing $m>n$ did not result in more deformations. This is in contrast with (\ref{planecurve.K.defns.thm}). 
\end{rem}

\begin{say}[Linear projections]\label{Cn.CC.flat.ord.1.exmp}   
Recall that by our convention (\ref{Cn.genflat.1st.say}.5),  $\phi_{ij}\equiv 0$ for $j>n$.  Extending this, in the following proof  all sums/products involving $i$ go from $1$ to $m$  and sums/products involving $j$ go from $1$ to $n$.

With ${\mathbf C}_n $ as in 
(\ref{Cn.CC.K.flat.thm}.1)
 consider the special projections
$$
\pi_{\mathbf a}:\a^n_{\mathbf x}[\epsilon]\to \a^2_{uv}[\epsilon]
\qtq{given by} u=\tsum x_i, v=\tsum a_ix_i,
\eqno{(\ref{Cn.CC.flat.ord.1.exmp}.1)}
$$
 where $a_i\in k[\epsilon]$. Write $a_i=\bar a_i+ a'_i\epsilon$. 
(One should think that $a'_i=\partial a_i/\partial\epsilon$.)

In order to compute the projection, we follow the method of (\ref{dsupp.car.1d.dl}). 
Since we compute over $k[u,u^{-1}, \epsilon]$, we may as well work with 
 the $k[u, \epsilon]$-module  $M:=\oplus_j k[x_j, \epsilon]$ and   write $1_j\in  k[x_j, \epsilon]$ for the $j$th unit. Then multiplication by $u$ and $v$ are given by
$$
\begin{array}{rcrcl}
u\cdot 1_j &=& (\tsum_i x_i)1_j &=& \hphantom{a_j}  x_j+\tsum_i \phi_{ij}\epsilon\qtq{and}\\
v\cdot 1_j &=& (\sum_i a_ix_i)1_j &=&  a_jx_j+\tsum_i a_i\phi_{ij}\epsilon.\\
\end{array}
\eqno{(\ref{Cn.CC.flat.ord.1.exmp}.2)}
$$
Thus
$$
v\cdot 1_j = \bigl(a_ju+\tsum_i(a_i-a_j)\phi_{ij}(u)\epsilon\bigr)\cdot 1_j.
$$
Thus the  $v$-action on $M$ is given by the diagonal matrix
$$
\diag \bigl(a_ju+\tsum_i(a_i-a_j)\phi_{ij}(u)\epsilon\bigr),
$$
and by (\ref{dsupp.car.1d.dl}) the equation of the projection is its characteristic polynomial
$$
\tprod_j\bigl(v-a_j u -\tsum_i (a_i-a_j) \phi_{ij}(u)\epsilon\bigr)=0.
\eqno{(\ref{Cn.CC.flat.ord.1.exmp}.3)}
$$
Expanding it we get an equation of the form
$$
\tprod_j\bigl(v-\bar{a}_j u\bigr)-
B(u,v, a, \phi)\epsilon = 0,
\eqno{(\ref{Cn.CC.flat.ord.1.exmp}.4)}
$$
where
$$
B(u,v, a, \phi)=\tsum_j\bigl(\tprod_{i\neq j}(v-\bar{a}_j u)\bigr)\cdot
\bigl( a'_ju+\tsum_i(\bar{a}_i-\bar{a}_j) \phi_{ij}(u)\bigr).
\eqno{(\ref{Cn.CC.flat.ord.1.exmp}.5)}
$$
This is a polynomial of degree $\leq n-1$ in $v$, hence by (\ref{ci.bais.lem.5})
its restriction to the curve $\bigl(\tprod_j\bigl(v-\bar{a}_j u\bigr)=0\bigr) $
is regular iff  $ B(u,v, a, \phi)$ is a polynomial in $u$ as well.  
Let now $r$ be the highest pole order of the $\phi_{ij}$ and write
$$
\phi_{ij}(u)=c_{ij}u^{-r}+(\mbox{higher terms}).
$$
Then the leading part of the coefficient of $v^{n-1}$ is
$$
 \tsum_j
\tsum_i(\bar{a}_i-\bar{a}_j) c_{ij}u^{-r}=
u^{-r}\tsum_i\ \bar{a}_i\bigl(\tsum_j (c_{ij}-c_{ji})\bigr).
\eqno{(\ref{Cn.CC.flat.ord.1.exmp}.6)}
$$
Since the $\bar{a}_i$ are arbitrary, 
 we get that  
$$
\tsum_j (c_{ij}-c_{ji})=0 \qtq{for every $i$.}
\eqno{(\ref{Cn.CC.flat.ord.1.exmp}.7)}
$$
Next we use a linear reparametrization of  the lines $x_i=\lambda_i^{-1}y_i$ and then apply a  projection  $\pi_{\mathbf a}$ as in (\ref{Cn.CC.flat.ord.1.exmp}.1).  The equations 
$x_i=\phi_{ij}(x_j)\epsilon$  become
$$
y_i=\lambda_i\phi_{ij}(\lambda_j^{-1}y_j)\epsilon
$$
and $c_{ij}$ changes to $\lambda_i\lambda_j^rc_{ij}$. Thus the equations (\ref{Cn.CC.flat.ord.1.exmp}.7)
 become
$$
\tsum_j (\lambda_i\lambda_j^rc_{ij}-\lambda_j\lambda_i^rc_{ji})=0
\quad \forall i.
\eqno{(\ref{Cn.CC.flat.ord.1.exmp}.8)}
$$
If $r\geq 2$ this implies that $c_{ij}=0$ and if $r=1$ then we get that
$c_{ij}=c_{ji}$.

This completes the proof of (\ref{Cn.CC.K.flat.thm}.2)
$\Leftrightarrow$ (\ref{Cn.CC.K.flat.thm}.4).
\medskip

{\it Remark \ref{Cn.CC.flat.ord.1.exmp}.9.} Note that if we work over $\f_2$ then necessarily $\lambda_i=1$, hence  (\ref{Cn.CC.flat.ord.1.exmp}.8) does not exclude the $r\geq 2$ cases. 
\end{say}

\begin{say}[Non-linear projections]\label{Cn.K.flat.nl.exmp} 
Consider a general non-linear projection
$$
(x_1,\dots, x_n)\mapsto  \bigl (\Phi_1( x_1,\dots, x_n), \Phi_2(x_1,\dots, x_n)\bigr).
$$
After a formal coordinate change we may assume that
$\Phi_1=\tsum_ix_i$.
Note that the monomials of the form  $x_ix_jx_k, x_i^2x_j^2, x_ix_j\epsilon$ vanish on ${\mathbf C}_n$, so we can discard these terms from $\Phi_2$. 
Thus, in suitable local coordinates
a general non-linear projection can be written as
$$
u=\tsum_i x_i,\quad  v=\tsum_i\alpha_i(x_i)+\tsum_{i\neq j}x_i\beta_{ij}(x_j),
\eqno{(\ref{Cn.K.flat.nl.exmp}.1)}
$$
where $\alpha_i(0)=\beta_{ij}(0)=0$. 
Note that  $\alpha'_i(0)=a_i$ in the notation of (\ref{Cn.CC.flat.ord.1.exmp}). 
Now we get that
$$
\begin{array}{l}
u\cdot 1_j=x_j+\tsum_i \phi_{ij}(x_j)\epsilon\qtq{and}\\
v\cdot 1_j=\alpha_j(x_j)+\tsum_{i\neq j}\alpha_i\bigl(\phi_{ij}(x_j)\epsilon\bigr)
+ \tsum_{i\neq j}\phi_{ij}(x_j)\beta_{ij}(x_j)\epsilon.
\end{array}
\eqno{(\ref{Cn.K.flat.nl.exmp}.2)}
$$
Note further that 
$\alpha_i\bigl(\phi_{ij}(x_j)\epsilon\bigr)=\alpha'_i(0)\phi_{ij}(x_j)\epsilon$ and
$$
\alpha_j(x_j)=\alpha_j\bigl(u-\tsum_i \phi_{ij}(x_j)\epsilon\bigr)=
\alpha_j(u)-\alpha'_j(u)\tsum_i \phi_{ij}(x_j)\epsilon.
$$
Thus, as in (\ref{Cn.CC.flat.ord.1.exmp}.4),  the projection  is
defined by the vanishing of
$$
\begin{array}{l}
\prod_j\Bigl(
v-\alpha_j(u)-\tsum_i\bigl(\beta_{ij}(u)+\alpha'_i(0)-\alpha'_j(u)\bigr)\phi_{ij}(u)\epsilon\Bigr)\\
\quad =: \prod_j\bigl(v-\bar\alpha_j(u)\bigr)-B(u,v,\alpha,\beta, \phi)\epsilon.
\end{array}
\eqno{(\ref{Cn.K.flat.nl.exmp}.3)}
$$
Let  $\bar\beta_{ij}, \bar\alpha'_j$ denote the residue of $\beta_{ij}, \alpha'_j$  modulo $\epsilon$   and write
$\alpha_j(u)=\bar\alpha_j(u)+\partial_{\epsilon}\alpha_j(u)\epsilon$. 
As in (\ref{Cn.CC.flat.ord.1.exmp}.5), expanding the product
 gives that $B(u,v,\alpha,\beta, \phi)$ equals
$$
\tsum_j
\bigl(\tprod_{i\neq j}(v-\bar\alpha_i(u))\bigr)\cdot \bigl(\partial_{\epsilon}\alpha_j(u)+\tsum_i\bigl(\bar\beta_{ij}(u)+\bar\alpha'_i(0)-\bar\alpha'_j(u)\bigr)\phi_{ij}\bigr).
\eqno{(\ref{Cn.K.flat.nl.exmp}.4)}
$$
We already know that  $\phi_{ij}(u)=c_{ij}u^{-1}+(\mbox{higher terms})$, hence
 $B(u,v,\alpha,\beta, \phi)$ has at most simple pole along $(u=0)$. Computing its residue along $u=0$ we get
$$
v^{n-1}\tsum_j\tsum_i\bigl(\bar\beta_{ij}(0)+\bar\alpha'_i(0)-\bar\alpha'_j(0)\bigr)c_{ij}=
v^{n-1}\tsum_{ij}(\bar{a}_i-\bar{a}_j)c_{ij}.
\eqno{(\ref{Cn.K.flat.nl.exmp}.5)}
$$
These are the same  as in (\ref{Cn.CC.flat.ord.1.exmp}.6). 
Thus $B(u,v,\alpha,\beta, \phi)$ is regular iff it is regular for the linearization of the projection.
This completes the proof of (\ref{Cn.CC.K.flat.thm}.4)
$\Rightarrow$ (\ref{Cn.CC.K.flat.thm}.3).
\end{say}

\begin{exmp}\label{Cn.anal.proj.n34.exmp}
The image of a  general linear projection of $C_n\subset \a^n$  to $\a^2$ is $n$ distinct lines through the origin. Their equation is 
$g_n(x,y)=0$ where $g_n$ is homogeneous of degree $n$ with simple roots only. 
A typical example is $g_n=x^n+y^n$.

A general non-linear projection to $\a^2$ gives  $n$ smooth curve germs with distinct tangent lines through the origin. The equation of the image is 
$g_n(x,y)+(\mbox{higher terms})=0$ where $g_n$ is homogeneous of degree $n$ with simple roots only. 

The miniversal deformation  of 
 $(x^n+y^n=0)$ is
$$
\bigl(x^n+y^n+\tsum_{i,j\leq n-2} t_{ij}x^iy^j=0\bigr)\subset \a^2_{xy}\times \a^{(n-1)^2}_{\mathbf t}.
\eqno{(\ref{Cn.anal.proj.n34.exmp}.1)}
$$
A general deformation is a smoothing, but deformations that have $n$ smooth branches with the same tangents as  $(x^n+y^n=0)$ form the subfamily
$$
\bigl(x^n+y^n+\tsum_{ij} t_{ij}x^iy^j=0\bigr)\subset \a^2_{xy}\times \a^{\binom{n-3}{2}}_{\mathbf t},
\eqno{(\ref{Cn.anal.proj.n34.exmp}.2)}
$$
where summation is over those pairs $(i,j)$ that satisfy
$i,j\leq n-2$ and $n< i+j$.   For $n\leq 4$ there is no such pair $(i,j)$, which gives the following.

\medskip
{\it Claim  \ref{Cn.anal.proj.n34.exmp}.3.} For $n\leq 4$ every analytic projection $\hat C_n\to \hat\a^2$ is obtained as the composite of an automorphism of $\hat C_n$, followed by a linear projection and then by an automorphism of $\hat \a^2$. \qed
\medskip

For $n=5$ we get the deformations
$$
(x^5+y^5+t x^3y^3=0)\subset \a^2_{xy}\times \a_t.
\eqno{(\ref{Cn.anal.proj.n34.exmp}.4)}
$$
For $t\neq 0$ we get curves that are  images of $\hat C_n$ by a nonlinear projection, but not as a linear projection pre-composed/composed with
automorphisms. 

\end{exmp}

The following strengthens \cite[4.11]{MR1714736}.

\begin{prop}\label{I.ch.axes.prop.1} The ideal of Chow equations of $C_n$ is generated by
\begin{enumerate}
\item all degree $n$ monomials, save the $x_i^n$, if   $n$ is even, and
\item all degree $n$ monomials, save the $x_i^n$ and $x_1\cdots x_n$, if   $n$ is odd.
\end{enumerate} 
These hold both for linear, polynomial and analytic projections.
\end{prop}

Note that 
we can write the even case as $I^{\rm ch}_{C_n}=I_{C_n}\cap  (x_1,\dots, x_n)^n$.
\medskip

Proof.  Every Chow equation has multiplicity $\geq n$, and we get the
same equations modulo  $(x_1,\dots, x_n)^{n+1}$, whether we use linear, polynomial and analytic projections (\ref{proj.formulas.exmp}).

In both of our cases, $I_{C_n}\cap  (x_1,\dots, x_n)^{n+1}\subset I^{\rm ch}_{C_n}$, so the ideal of Chow equations coming from linear projections already contains every possible higher order monomial. Thus it is sufficient to prove (1--2) for linear projections.

The linear projections of $C_n$ to $\a^2_{uv}$ are given by
$u=\sum_i a_ix_i, v=\sum b_ix_i$. The image of the $x_j$-axis  is
$(b_ju-a_jv=0)$. So the pull-back of their product is
$$
\tprod_j \tsum_i (a_ib_j-a_jb_i)x_i.
\eqno{(\ref{I.ch.axes.prop.1}.3)}
$$
Since $C_n$ is toric, $I^{\rm ch}_{C_n}$ is a monomial ideal. 
Thus we need to understand which degree $n$ monomials in the $x_i$ have a
nonzero coefficient in  (\ref{I.ch.axes.prop.1}.3).

First, the coefficient of $x_j$ in $\tsum_i (a_ib_j-a_jb_i)x_i $ is 0, 
so we never get $x_j^n$. 
Next consider  $x_1\cdots x_n$. Its coefficient is
$$
\tsum_{\sigma\in S_n} \prod_i\bigl(a_ib_{\sigma(i)}-a_{\sigma(i)}b_i\bigr).
\eqno{(\ref{I.ch.axes.prop.1}.4)}
$$
Note that the product is 0 if $\sigma(i)=i$ for some $i$ and changes by
$(-1)^n$ when $\sigma$ is replaced by $\sigma^{-1}$. Thus if $n$ is odd then 
(\ref{I.ch.axes.prop.1}.4) is identically zero.  (More generally, the permanent of a skew-symmetric matrix of odd size is 0.) We have thus proved the following.
\medskip

{\it Claim \ref{I.ch.axes.prop.1}.5.} If $n$ is odd then the coefficient of  $x_1\cdots x_n$ in  (\ref{I.ch.axes.prop.1}.3) is 0.\qed
\medskip

It remains to show that all other degree $n$ monomials appear in  (\ref{I.ch.axes.prop.1}.3) with nonzero coefficient.

To show this we choose specific values of the $a_i, b_i$ and hope to get 
enough nonzero terms. Thus fix $1\leq r\leq n$, choose  $a_1=\cdots=a_r=1, a_{r+1}=\cdots=a_n=0$ and $b_1=\cdots=b_r=0, b_{r+1}=\cdots=b_n=1$. 
Then 
$$
a_ib_j-a_jb_i=
\left\{
\begin{array}{cl}
1 &\mbox{ if } i\leq r<j,\\
-1 &\mbox{ if } j\leq r<i, \qtq{and}\\
0 &\mbox{ otherwise.}
\end{array}
\right.
\eqno{(\ref{I.ch.axes.prop.1}.6)}
$$
 Thus (\ref{I.ch.axes.prop.1}.3) becomes
$$
(-1)^{r}(x_1+\cdots+x_r)^{n-r}(x_{r+1}+\cdots+x_n)^r
\eqno{(\ref{I.ch.axes.prop.1}.7)}
$$
Applying this to various permutations of the $x_i$ and choices of $r$  we get the following.
\medskip

{\it Claim \ref{I.ch.axes.prop.1}.8.} Let $M=\prod x_i^{w_i}$ be a degree $n$ monomial. Then $M\in  I^{\rm ch}_{C_n}$ if the following holds.
\begin{enumerate}
\item[(a)] There is a subset $I\subset \{1,\dots, n\}$ such that
$\tsum_{i\in I}w_i=n-|I|$. \qed
\end{enumerate}

While this is only a sufficient condition,
we check in  (\ref{I.ch.axes.prop.lem}) that it applies to every monomial other than $x_i^n$ and $x_1\cdots x_n$ for $n$ odd.
This completes the proof of (\ref{I.ch.axes.prop.1}). \qed

\begin{lem} \label{I.ch.axes.prop.lem} Let $M=\prod x_i^{w_i}$ be a degree $n$ monomial other than  $x_i^n$ or $x_1\cdots x_n$ for odd $n$. 
Then there is a subset $I\subset \{1,\dots, n\}$ such that
$\tsum_{i\in I}w_i=n-|I|$.
\end{lem}

Proof. We use induction on $n$, the case $n=1$ is empty and $n=2$ is obvious. 

Assume first that $w_{n-1}=w_n=1$.  
If $M=x_1^{n-2}x_{n-1}x_n$ then $I=\{1,2\}$ works. 
Otherwise, by induction, 
there is a subset $J\subset \{1,\dots, n-2\}$ such that
$\tsum_{i\in J}w_i=n-2-|J|$. Set $I=J\cup\{n\}$. Then 
$\tsum_{i\in J}w_i=n-2-|J|+1=n-|I|$ and we are done.

If the inductive step does not apply, then there is at most one $w_i=1$, hence at least $\tfrac{n-1}{2}$ of the $w_i=0$.

Reorder the $x_i$ such that $w_i$ is a  decreasing function and take $r$ such that
$w_1+\cdots+w_{r-1}<n/2$ but $w_1+\cdots+w_{r}\geq n/2$. 
If $w\leq n-r$ then we take
$I=\{1,\dots, r, n-s,\dots, n\}$ where $s=n-r-w-1$. 
Since $w_i=0$ for $i\geq \tfrac{n+1}{2}$, 
$$
 \tsum_{i\in J}w_i=w_1+\cdots+w_{r}=w\qtq{and} |I|=r+s+1=n-w.
$$
What happens if $w>n-r$?. Note that then $r\geq 2$ and  $w_1\geq\cdots\geq w_r\geq 2$ so
$(r-1)w_r<n/2$ and $2(r-1)<n/2$. On the other hand, 
$w_1+\cdots+w_{r}< n/2+w_r<n/2+n/(2r-2)$. 
One checks that $n/2+n/(2r-2)>n-r$ and $2(r-1)<n/2$ both hold only for
$r=2$.
Furthermore,  the only monomial for which the above choice of $I$ does not work is $x_1^{(n-1)/2}x_2^{(n-1)/2}x_3$ for $n$ odd. In this case we can take
$I=\{1, 3, n-s,\dots, n\}$ where $s=\tfrac{n-5}{2}$. \qed

\begin{prop} \label{Cn.CC.ideal.defs.prop} 
Consider the deformation ${\mathbf C}_n $ as in 
(\ref{Cn.CC.K.flat.thm}.1). Assume that  $n\geq 3$.
Then $I^{\rm ch}(C_n)$ vanishes on the central fiber $({\mathbf C}_n)_k $ 
 iff $\ord \phi_{ij}\leq n-2$ for every $i\neq j$. 
\end{prop}

Proof. Note that $x_ix_jx_k, x_i^2x_j^2\in I_{C_n}^{(2)}$ hence
the only condition is the liftability of $x_ix_j^{n-1}$.

If $\ord \phi_{ij}\leq n-2$ then
$x_ix_j^{n-1}-x_j^{n-1}\phi_{ij}(x_j)\epsilon$ vanishes along the $x_j$-axis
and everywhere else.

Conversely, assume that we have equations
$$
x_ix_j^{n-1}-\Psi_{ij}(x_1,\dots, x_n)\epsilon=0.
$$
As in (\ref{Cn.genflat.1st.say}.2), they simplify to
$$
x_ix_j^{n-1}-\psi_{ij}(x_j)\epsilon=0 \qtq{along the $x_j$-axis.}
$$
Note that $x_ix_j^{n-1}$ vanishes along the other $n-2$ axes, so we must have
$\psi_{ij}(0)=0$.    Thus $\phi_{ij}:=x_j^{1-n}\psi_{ij}$ has a pole of order at most $n-2$. \qed


\def\cprime{$'$} \def\cprime{$'$} \def\cprime{$'$} \def\cprime{$'$}
  \def\cprime{$'$} \def\dbar{\leavevmode\hbox to 0pt{\hskip.2ex
  \accent"16\hss}d} \def\cprime{$'$} \def\cprime{$'$}
  \def\polhk#1{\setbox0=\hbox{#1}{\ooalign{\hidewidth
  \lower1.5ex\hbox{`}\hidewidth\crcr\unhbox0}}} \def\cprime{$'$}
  \def\cprime{$'$} \def\cprime{$'$} \def\cprime{$'$}
  \def\polhk#1{\setbox0=\hbox{#1}{\ooalign{\hidewidth
  \lower1.5ex\hbox{`}\hidewidth\crcr\unhbox0}}} \def\cdprime{$''$}
  \def\cprime{$'$} \def\cprime{$'$} \def\cprime{$'$} \def\cprime{$'$}
\providecommand{\bysame}{\leavevmode\hbox to3em{\hrulefill}\thinspace}
\providecommand{\MR}{\relax\ifhmode\unskip\space\fi MR }
\providecommand{\MRhref}[2]{%
  \href{http://www.ams.org/mathscinet-getitem?mr=#1}{#2}
}
\providecommand{\href}[2]{#2}

\bigskip

  Princeton University, Princeton NJ 08544-1000

\email{kollar@math.princeton.edu}

\end{document}